\documentclass[10pt,a4paper,draft]{article}
\usepackage{bbm}
\usepackage{pifont}
\usepackage{bbding}
\usepackage{amsmath, esint}
\usepackage{mathrsfs}
\usepackage{amsfonts}
\usepackage{amssymb}
\usepackage{amsfonts,amssymb,amsmath,indentfirst,cases,amsthm}
\usepackage{epsfig}
\usepackage[numbers,sort&compress]{natbib}
\allowdisplaybreaks
\def\ess~sup{\mathop{\rm ess~sup}}

\numberwithin{equation}{section}

\newenvironment{key words}{\emph{\texttt{Keywords}}\mbox{  }}{ }
\newtheorem{theorem}{Theorem}[section]
\newtheorem{lemma}[theorem]{Lemma}
\newtheorem{remark}[theorem]{Remark}

\theoremstyle{remark}

\theoremstyle{plain}

\makeatletter

\newcommand{\Rmnum}[1]{\expandafter\@slowromancap\romannumeral #1@}
\makeatother
\begin{document}

 \title{\textbf{Weak Differentiability to Nonuniform Nonlinear Degenerate Elliptic Systems under $p,q$-growth Condition
 on the Heisenberg Group
}}
\author{Junli Zhang$^1$\thanks{Corresponding
author's E-mail: jlzhang2020@163.com}, Zhouyu Li$^2$\\
\small{1. School of Mathematics and Data Science, Shaanxi University of Science and Technology,}\\
\small{ Xi'an, Shaanxi, 710021, P. R. China}\\
\small{2. School of Sciences, Xi'an University of Technology,}\\
\small{ Xi'an, Shaanxi, 710054, P. R. China}}

\date{} \maketitle

\maketitle {\bf Abstract}\
The paper concerns the weak differentiability of weak solutions to two kinds of nonuniform nonlinear degenerate elliptic systems under the $p,q$-growth condition on the Heisenberg Group. We use the iteration to fractional difference quotients on the Heisenberg Group to get the weak differentiability of weak solution $u$ in the vertical direction (i.e., ${L^p}$($1<p<4$) integrability of $Tu$) and then the second order weak differentiability of weak solution in the horizontal directions (i.e., ${L^2}$ integrability of $\nabla _H^2u$) and weak differentiability of gradient of weak solution in the vertical direction (i.e., ${L^2}$ integrability of $T{\nabla _H}u$)).


\textbf {AMS} Primary: 35H20, 35B65; Secondary: 35J60, 35J70.
\def\Xint#1{\mathchoice
    {\XXint\displaystyle\textstyle{#1}}%
    {\XXint\textstyle\scriptstyle{#1}}%
    {\XXint\scriptstyle\scriptscriptstyle{#1}}%
    {\XXint\scriptscriptstyle\scriptscriptstyle{#1}}%
    \!\int}
\def\XXint#1#2#3{{\setbox0=\hbox{$#1{#2#3}{\int}$}
    \vcenter{\hbox{$#2#3$}}\kern-.5\wd0}}
\def\dashint{\Xint-}
\section{Introduction}

In this paper, we attempt to establish the weak differentiability of weak solutions (i.e., integrabilities of $Tu$
, $\nabla _H^2u$ and $T{\nabla _H}u$) to the nonuniform nonlinear degenerate elliptic system
\begin{equation}\label{eq11}
  {div}_H\left( {A\left( {x,{\nabla _H}u} \right)} \right) = 0,\;x \in \Omega
\end{equation}
and the nonuniform nonlinear degenerate elliptic system with drift $Tu$
\begin{equation}\label{eq12}
  {div}_H\left( {A\left( {x,{\nabla _H}u} \right)} \right) + Tu = 0,\;x \in \Omega
\end{equation}
under the $p,q$-growth condition on the Heisenberg Group ${\mathbb{H}^n}$, where $\Omega  \subset {\mathbb{H}^n}$ is a bounded domain, $u = \left( {{u_1}, \cdots ,{u_N}} \right)$ is a vector valued function, ${\nabla _H}u$ and ${{div}_H} A$ denote
\[{\nabla _H}u{\rm{ = }}\left[ {\begin{array}{*{20}{c}}
   {{X_1}{u_1}} &  \cdots  & {{X_{2n}}{u_1}}  \\
    \vdots  &  \cdots  &  \vdots   \\
   {{X_1}{u_N}} &  \cdots  & {{X_{2n}}{u_N}}  \\
\end{array}} \right]\]
and
\[di{v_H}A = \sum\limits_{j = 1}^{2n} {{X_j}{A_j}} ,\;A = A(x,{\nabla _H}u),\]
respectively, $A = ({A_1}, \cdots ,{A_{2n}}),{A_j} = \left( {{A_{1,j}},{A_{2,j}}, \cdots ,{A_{N,j}}} \right),j = 1, \cdots ,2n$. For any $P,U \in {\mathbb{R}^{ N \times 2n }}$ and almost everywhere $x,{x_0} \in \Omega $, the matrix function $A = A(x,{\nabla _H}u):\Omega  \times {\mathbb{R}^{N \times 2n}} \to {\mathbb{R}^{N \times 2n}}$ satisfies the $p,q$-growth condition
\begin{equation}\label{eq13}
\langle {A\left( {x,P} \right) - A\left( {x,U} \right)} , {P - U} \rangle \ge a\left( x \right){\left( {1 + {{\left| P \right|}^2} + {{\left| U \right|}^2}} \right)^{\frac{{p - 2}}{2}}}{\left| {P - U} \right|^2},
\end{equation}
\begin{equation}\label{eq14}
\left| {A\left( {x,P} \right) - A\left( {x,U} \right)} \right| \le b\left( x \right){\left( {1 + {{\left| P \right|}^2} + {{\left| U \right|}^2}} \right)^{\frac{{q - 2}}{2}}}\left| {P - U} \right|
\end{equation}
and the continuity condition
\begin{equation}\label{eq15}
  \left| {A\left( {x,P} \right) - A\left( {{x_0},P} \right)} \right| \le cd\left( {x,{x_0}} \right)k\left( x \right){\left( {1 + {{\left| P \right|}^2}} \right)^{\frac{{q - 1}}{2}}},
\end{equation}
where $d\left( {x,{x_0}} \right)$ denote the Kor\`{a}nyi metric between $x$ and $x_0$, and nonnegative measurable functions $a\left( x \right),\;b\left( x \right)$ and $k\left( x \right)$ satisfy
\begin{equation}\label{eq16}
{a(x)^{ - 1}} \in L_{loc}^\infty \left( \Omega  \right),\;b(x),k(x) \in L_{loc}^m\left( \Omega  \right),
\end{equation}
the exponents $p,\;q$ and $m$ satisfy
\begin{equation}\label{eq17}
{\rm{1}} < p \le q,\;m > \frac{{2p}}{{p - 2}},\frac{q}{p} < 1 + \frac{1}{Q} - \frac{1}{m},
\end{equation}
in which $Q = 2n + 2$ is the homogeneous dimension of $\mathbb{H}^n$. It should be emphasized that the matrix $A = A(x,{\nabla _H}u)$ of coefficients to system \eqref{eq11} does not need to satisfy the uniform elliptic condition for $P,U \in {\mathbb{R}^{ N \times 2n }}$, because $a\left( x \right)$ in \eqref{eq13} does not need to be greater than zero, but can be zero on a zero measure set, and $b\left( x \right)$ in \eqref{eq14} can be unbounded on a zero measure set. We note
\[Tu = {X_i}{X_{n + i}}u - {X_{n + i}}{X_i}u\]
and write
\[\nabla _H^2u = {\left( {{X_i}{X_j}u} \right)_{1 \le i,j \le 2n}}.\]

System \eqref{eq11} considered by us includes system
\begin{equation*}
di{v_H}\left( {a\left( x \right){{\left( {1 + {{\left| {{\nabla _H}u} \right|}^2}} \right)}^{\frac{{p - 2}}{2}}}{\nabla _H}u + b\left( x \right){{\left( {1 + {{\left| {{\nabla _H}u} \right|}^2}} \right)}^{\frac{{q - 2}}{2}}}{\nabla _H}u} \right) = 0
\end{equation*}
(where $a(x)$ and $b(x)$ satisfy \eqref{eq16}) and system
\begin{equation*}
 di{v_H}\left( {a\left( x \right){{\left( {1 + {{\left| {{\nabla _H}u} \right|}^2}} \right)}^{\frac{{p\left( x \right) - 2}}{2}}}{\nabla _H}u + b\left( x \right){{\left( {1 + {{\left| {{\nabla _H}u} \right|}^2}} \right)}^{\frac{{q\left( x \right) - 2}}{2}}}{\nabla _H}u} \right) = 0,
\end{equation*}
where $a(x)$ and $b(x)$ satisfy \eqref{eq16}, $p(x)$ and $q(x)$  meet
\[p \le p\left( x \right) \le q\left( x \right) \le q,\;p\left( x \right),q\left( x \right) \in {L^m}\left( \Omega  \right).\]

Many scholars have studied regularity of solutions to nonlinear elliptic equations (systems) under the $p,q$-growth condition on the Euclidean space, such as ${C^\alpha }$ regularity, ${C^{1,\alpha }}$ regularity and Lipschitz continuity and so on. It is noted that many results were obtained for the equations (systems) satisfying the uniform elliptic condition (i.e., $0 < {c_1} \le a\left( x \right) \le b\left( x \right) \le {c_2}$), see \cite{BCM20}, \cite{CM15}, \cite{CM151}, \cite{CGM03}, \cite{CMM14}, \cite{DM20}, \cite{ELP19}, \cite{ELM02}, \cite{ELM04}, \cite{EMM16}, \cite{EMM20}, \cite{M91} and \cite{M20}. Some results were derived for the equations (systems) not satisfying the uniform elliptic condition (i.e., $0 \le a\left( x \right) \le b\left( x \right) \le  + \infty $), see \cite{BCM18}, \cite{CMM18} and \cite{CMMP21}. For viscosity solutions of the double-phase equation see \cite{FZ20} and \cite{FZ201}. For more detailed introduction for regularity of solutions to nonlinear elliptic equations (systems) under the $p,q$-growth condition, see \cite{M21}.

An interesting problem is whether these regularity results hold to nonlinear degenerate elliptic equations (systems) on the Heisenberg group $\mathbb{H}^n$? Wang, Manfredi, Zhu et al consider the partial regularity of weak solutions to nonlinear uniform degenerate elliptic equations satisfying $p$-growth condition, see \cite{WM18}, \cite{WLZ19} and \cite{WZGL20}. Zhang and Niu in \cite{ZN20} have studied ${C^\alpha }$ regularity of weak solutions to nonlinear nonuniform degenerate elliptic equations satisfying generalized Orlicz growth conditions (including $p$-growth, $p,q$-growth and variable exponential growth conditions). Since $\mathbb{H}^n$ is a noncommutative two step nilpotent Lie group and the equations (systems) on the group is degenerate elliptic, in order to obtain ${C^{1,\alpha }}$ regularity and Lipschitz continuity of weak solutions, it is usually necessary to establish first the weak differentiabilities of weak solutions.

${C^\alpha }\left( {0 < \alpha  < 1} \right)$ regularity of weak solutions to nonlinear degenerate elliptic equations composed of vector fields ${X_i}$ satisfying the H\"{o}rmander condition has been studied by some scholars, see \cite{CDG93}, \cite{L96} and \cite{X92}. Specifically, in \cite{X92}, ${C^\infty }$ smoothness of weak solutions to the non-divergence equation
\[\sum\limits_{i,j = 1}^m {{A_{ij}}\left( {x,u,Xu} \right){X_i}{X_j}u}  + B\left( {x,u,Xu} \right) = 0\]
was obtained under the assumption that second order weak derivatives ${X_i}{X_j}u$ of weak solutions satisfy ${X_i}{X_j}u \in {C^\alpha }$. The authors in \cite{CDG93} investigated the divergence equation
\[\sum\limits_{j = 1}^m {X_j^*{A_j}\left( {x,u,Xu} \right)}  = f\left( {x,u,Xu} \right)\]
under the $p$-growth condition and proved ${C^\alpha }\left( {0 < \alpha  < 1} \right)$ regularity
by using the Sobolev embedding theorem, the Poincar\'{e} inequality and the Moser iteration. Lu
in \cite{L96} concluded ${C^\alpha }\left( {0 < \alpha  < 1} \right)$ regularity to the divergence equation
\[\sum\limits_{j = 1}^m {X_j^*\left( {{{\left| {Xu} \right|}^{p - 2}}{X_j}u} \right)}  = A{\left| {Xu} \right|^p} + B{\left| {Xu} \right|^{p - 1}} + C{\left| u \right|^{p - 1}} + D,\;1 < p < \infty \]
by applying various embedding theorems into nonisotropic Lipschitz spaces and BMO spaces.

${C^{1,\alpha }}\left( {0 < \alpha  < 1} \right)$ regularity of weak solutions
to the uniform nonlinear divergence equation on $\mathbb{H}^n$
\begin{equation}\label{eq112}
di{v_H}\left( {Dg\left( {x,{\nabla _H}u} \right)} \right) = h\left( x \right)
\end{equation}
under $p$-growth condition has also been concerned by many scholars. For example, Capogna in \cite{CL97} obtained H\"{o}lder regularity of Euclidean gradients of weak solutions to \eqref{eq112}, where $Dg$ is differentiable and satisfies
\begin{equation*}\label{eqb15}
  \left\{ \begin{array}{l}
{L^{ - 1}}{\left| \xi  \right|^2} \le \left\langle {{D^2}g(z)\xi ,\;\xi } \right\rangle  \le L{\left| \xi  \right|^2},\\
\left| {Dg(z)} \right| \le L{\left( {\delta  + {{\left| z \right|}^2}} \right)^{\frac{1}{2}}},
\end{array} \right.
\end{equation*}
where $\delta  \ge 0,\:L \ge 1$  are constants. To get this result, the author first proved that Heisenberg gradients and vertical gradients of weak solutions belong to the local Sobolev space $HW_{loc}^{1,2}\left( \Omega  \right)$ by using the fractional difference quotients, and then gained ${C^\infty }$ smoothness of weak solutions to the uniform nonlinear divergence equation
\begin{equation*}\label{eq119}
di{v_H}\left( {Dg\left( {x,{\nabla _H}u} \right)} \right) = 0.
\end{equation*}
 These results in \cite{CL97} were extended to the Carnot group, see \cite{CL99}. The weak differentiability of weak solutions (i.e., $Tu \in L_{loc}^p\left( \Omega  \right)\left( {1 < p < 4} \right)$, $u \in HW_{loc}^{2,p}\left( \Omega  \right)\left( {\frac{{\sqrt {17}  - 1}}{2} \le p \le 2} \right)$, $u \in HW_{loc}^{2,2}\left( \Omega  \right)\left( {2 \le p < 4} \right)$, ${\nabla _H}Tu \in L_{loc}^2\left( \Omega  \right)\left( {\frac{{\sqrt {17}  - 1}}{2} \le p \le 2} \right)$) to subelliptic $p$-Laplacian equations
\begin{equation}\label{eq119}
 - \sum\limits_{i = 1}^{2n} {{X_i}\left( {{{\left( {\mu  + {{\left| {{\nabla _H}u} \right|}^2}} \right)}^{\frac{{p - 2}}{2}}}{X_i}u} \right)}  = 0,\;\mu  > 0
\end{equation}
were deduced by Domokos as $1 < p < 4$ in \cite{DA04} by the fractional difference quotients, which generalized the results in \cite{MS03} and \cite{MS031}. Manfredi and Mingione in \cite{MM07} derived Lipschitz continuity of weak solutions to \eqref{eq119}, where $p$ is not far from 2, and ${C^\infty }$ smoothness by using the method in \cite{CL97}. Domokos and Manfredi in \cite{DM05} used Calder\'{o}n-Zygmund theory to treat ${C^{1,\alpha }}$ regularity of weak solutions to \eqref{eq119} and
 \begin{equation}\label{eq113}
\sum\limits_{i = 1}^{2n} {{X_i}\left( {{{\left| {{\nabla _H}u} \right|}^{p - 2}}{X_i}u} \right)}  = 0,
 \end{equation}
 where $p$ is not far from 2. Based on \cite{DA04}, Mingione, Zatorska-Goldstein and Zhong in \cite{MZZ07} concluded ${C^{1,\alpha }}$-regularity of weak solutions to equation \eqref{eq119} and the Lipschitz continuity of weak solutions to equation \eqref{eq113} as $2 \le p < 4$ by using the energy estimate and interpolation inequality, and applying a double-bootstrap method to expand the scope of $p$ in $Tu \in L_{loc}^p\left( \Omega  \right)$
 in \cite{DA04} from $1 < p < 4$ to $1 < p < \infty $. The result in \cite{MZZ07} breaks through the limit of topological dimension $2n+1$ of $\mathbb{H}^n$ to $p$. Ricciotti in \cite{R18} presented a proof of the local H\"older regularity of the horizontal derivatives of weak solutions to \eqref{eq119} and \eqref{eq113} in the Heisenberg group $\mathbb{H}^1$ for $p>4 $. For more results, also see Ricciotti's book \cite{R15}. ${C^{1,\alpha }}$ regularity to \eqref{eq119} and \eqref{eq113} as $p\ge 2$ on Lie groups was proved by Domokos and Manfredi in \cite{DM20}. Mukherjee and Zhong in \cite{MZ17} got ${C^{1,\alpha }}$-regularity of weak solutions to \eqref{eq119} and \eqref{eq113} as $1 < p < \infty $ by using the energy estimate, Moser iteration and oscillation estimate.

It is worth noting that the above results are all obtained for the equation satisfying the uniform elliptic condition under the $p$-growth condition. In this paper, the coefficient matrix in \eqref{eq11} satisfies the nonuniform elliptic condition under the $p,q$-growth condition. Clearly, when $p=q$, the $p,q$-growth condition becomes $p$-growth condition, and the nonuniform elliptic condition includes the uniform elliptic condition. Inspired by \cite{CMM18} and \cite{CMMP21}, we deal with the case that $a\left( x \right)$ in \eqref{eq13} allows 0 on a zero measure set and $b\left( x \right)$ in \eqref{eq14} and $k\left( x \right)$ in \eqref{eq15} don't need to be bounded on the zero measure set, but only integrable.

We study the weak differentiabilities of weak solutions $u \in H{W^{1,{\tilde q}}}\left( \Omega  \right)$ to systems \eqref{eq11} and \eqref{eq12} under conditions \eqref{eq13}-\eqref{eq17}, where
\[\tilde q = \left\{ \begin{array}{l}
 \frac{{mp\left( {2q - 2 - p} \right)}}{{{m\left( {p - 2} \right) - 2p} }},\;2 < p < \infty , \\
 \frac{{mp\left( {q - 1} \right)}}{{ {m\left( {p - 1} \right) - p} }},\;1 < p \le 2. \\
 \end{array} \right.\]
Obviously, when $m = \infty ,p = q$, we have $\tilde q = p$, which returns to the situation in \cite{DA04}.

The main results are as follows:
\begin{theorem}\label{Th11}
 Let system \eqref{eq11} satisfy the conditions \eqref{eq13}-\eqref{eq17}. If $u \in H{W^{1,{\tilde q}}}\left( \Omega  \right)$ (the definitions of $H{W^{1,{\tilde q}}}\left( \Omega  \right)$ and $HW_0^{1,{\tilde q}}\left( \Omega  \right)$ see Section 2.1) is a weak solution to \eqref{eq11}, i.e., for any $\varphi  \in HW_0^{1,{\tilde q}}\left( \Omega  \right)$,
\begin{equation}\label{eq212}
\int_\Omega  {\langle A\left( {x,{\nabla _H}u} \right),{\nabla _H}\varphi\rangle dx}  = 0,
\end{equation}
 then for ${x_0} \in \Omega ,\;r > 0$ with $B\left( {{x_0},3r} \right) \subset \Omega $, there exist $l \in \mathbb{N}$ depending only on $p$ and the constant $c$ depending on $p,\;m$ and ${\left\| {{a^{ - 1}}} \right\|_{{L^\infty }\left( \Omega  \right)}}$ such that for $2 < p < 4,$ it holds
 \begin{equation}\label{eq18}
 \int_{B\left( {{x_0},\frac{r}{{{2^{l + 1}}}}} \right)} {{{\left| {Tu} \right|}^p}dx}  \le c\kappa^{l+1} \left[ {\int_{B\left( {{x_0},2r} \right)} {\left( {{{\left( {1 + {{\left| {{\nabla _H}u} \right|}^2}} \right)}^{\frac{{mp\left( {2q - 2 - p} \right)}}{{2\left[ {m\left( {p - 2} \right) - 2p} \right]}}}} + {{\left| u \right|}^p}} \right)dx}  + 1} \right];
 \end{equation}
for $1<p\le2$,
\begin{equation}\label{eq19}
 \int_{B\left( {{x_0},\frac{r}{{{2^{l + 1}}}}} \right)} {{{\left| {Tu} \right|}^p}dx}  \le c\kappa^{l+1} \left[ {\int_{B\left( {{x_0},2r} \right)} {\left( {{{\left( {1 + {{\left| {{\nabla _H}u} \right|}^2}} \right)}^{\frac{{mp(q-1)}}{{2\left[ {m\left( {p - 1} \right) - p} \right]}}}} + {{\left| u \right|}^p}} \right)dx}  + 1} \right],
\end{equation}
where $\kappa  = \left( {1 + {{\left\| b \right\|}_{{L^m}\left( \Omega  \right)}} + {{\left\| k \right\|}_{{L^m}\left( \Omega  \right)}} + {{\left\| {{a^{ - 1}}} \right\|}_{{L^\infty }\left( \Omega  \right)}}\left( {\left\| k \right\|_{{L^m}\left( \Omega  \right)}^2 + \left\| b \right\|_{{L^m}\left( \Omega  \right)}^2} \right)} \right)$. Furthermore, \eqref{eq18} and \eqref{eq19} imply $Tu \in L_{loc}^p\left( \Omega  \right)$ for $1<p<4$.
\end{theorem}

We note that when $4 \le p < \infty $, it follows
\[\mathop {\sup }\limits_{0 < \left| s \right| \le h} \int_{B\left( {{x_0},\frac{r}{{{2^l}}}} \right)} {{{\left| {{D_{T,s,\frac{1}{2} + \alpha '}}u} \right|}^p}dx}  \le c{\kappa ^l}\left[ {\int_{B\left( {{x_0},2r} \right)} {\left( {{{\left( {1 + {{\left| {{\nabla _H}u} \right|}^2}} \right)}^{\frac{{mp\left( {2q - 2 - p} \right)}}{{2\left[ {m\left( {p - 2} \right) - 2p} \right]}}}} + {{\left| u \right|}^p}} \right)dx}  + 1} \right]\]
for the weak solution $u$ to system \eqref{eq11}, where $\alpha '$  is smaller than $\frac{1}{{p - 2}}$ but arbitrarily close to $\frac{1}{{p - 2}}$  and $l \in \mathbb{N}$ is a constant depending only on $p$.

\begin{theorem}\label{Th12}
 Let system \eqref{eq11} satisfy the conditions \eqref{eq13}-\eqref{eq17}. If $u \in H{W^{1,{\tilde q}}}\left( \Omega  \right)$ is a weak solution to \eqref{eq11}, then for ${x_0} \in \Omega ,\;r > 0$ with $B\left( {{x_0},3r} \right) \subset \Omega $, there exist $l \in \mathbb{N}$ depending only on $p$ and the constant $c$ depending on $p,\;m$ and ${\left\| {{a^{ - 1}}} \right\|_{{L^\infty }\left( \Omega  \right)}}$ such that for $2 < p < 4,$
 \begin{align}\label{eq110}
 &\int_{B\left( {{x_0},\frac{r}{{{2^{l + 2}}}}} \right)} {{{\left( {1 + {{\left| {{\nabla _H}u} \right|}^2}} \right)}^{\frac{{p - 2}}{2}}}{{\left| {\nabla _H^2u} \right|}^2}dx} \nonumber\\
 \le &c{\kappa ^{l+2}}\left[ {\int_{B\left( {{x_0},2r} \right)} {\left( {{{\left( {1 + {{\left| {{\nabla _H}u} \right|}^2}} \right)}^{\frac{{mp\left( {2q - 2 - p} \right)}}{{2\left[ {m\left( {p - 2} \right) - 2p} \right]}}}} + {{\left| u \right|}^p}} \right)dx}  + 1} \right],
 \end{align}
where $\kappa  = \left( {1 + {{\left\| b \right\|}_{{L^m}\left( \Omega  \right)}} + {{\left\| k \right\|}_{{L^m}\left( \Omega  \right)}} + {{\left\| {{a^{ - 1}}} \right\|}_{{L^\infty }\left( \Omega  \right)}}\left( {\left\| k \right\|_{{L^m}\left( \Omega  \right)}^2 + \left\| b \right\|_{{L^m}\left( \Omega  \right)}^2} \right)} \right)$. Thus it implies $u \in
HW_{loc}^{2,2}\left( \Omega  \right)$ (its definition sees Section 2.1).
\end{theorem}

\begin{theorem}\label{Th13}
 Let system \eqref{eq11} satisfy the conditions \eqref{eq13}-\eqref{eq17}. If $u \in H{W^{1,{\tilde q}}}\left( \Omega  \right)$ is a weak solution to \eqref{eq11}, then for ${x_0} \in \Omega ,\;r > 0$ with $B\left( {{x_0},3r} \right) \subset \Omega $, there exist $l \in \mathbb{N}$ depending only on $p$ and the constant $c$ depending on $p,\;m$ and ${\left\| {{a^{ - 1}}} \right\|_{{L^\infty }\left( \Omega  \right)}}$ such that for $2 < p < 4,$
 \begin{align}\label{eq111}
 & \int_{B\left( {{x_0},\frac{r}{{{2^{l + 2}}}}} \right)} {{{\left( {1 + {{\left| {{\nabla _H}u} \right|}^2}} \right)}^{\frac{{p - 2}}{2}}}{{\left| {T{\nabla _H}u} \right|}^2}dx}\nonumber\\
 \le & c{\kappa ^{l+2}}\left[ {\int_{B\left( {{x_0},2r} \right)} {\left( {{{\left( {1 + {{\left| {{\nabla _H}u} \right|}^2}} \right)}^{\frac{{mp\left( {2q - 2 - p} \right)}}{{2\left[ {m\left( {p - 2} \right) - 2p} \right]}}}} + {{\left| u \right|}^p}} \right)dx}  + 1} \right],
 \end{align}
where $\kappa  = \left( {1 + {{\left\| b \right\|}_{{L^m}\left( \Omega  \right)}} + {{\left\| k \right\|}_{{L^m}\left( \Omega  \right)}} + {{\left\| {{a^{ - 1}}} \right\|}_{{L^\infty }\left( \Omega  \right)}}\left( {\left\| k \right\|_{{L^m}\left( \Omega  \right)}^2 + \left\| b \right\|_{{L^m}\left( \Omega  \right)}^2} \right)} \right)$. Thus it follows $T{\nabla _H}u = {\nabla _H}Tu \in L_{loc}^2\left( \Omega  \right)$.
\end{theorem}

\begin{remark}\label{Re14}
Let system \eqref{eq12} satisfy the conditions \eqref{eq13}-\eqref{eq17}. If $u \in H{W^{1,{\tilde q}}}\left( \Omega  \right)$ is a weak solution to \eqref{eq12}, i.e., for any $\varphi  \in HW_0^{1,{\tilde q}}\left( \Omega  \right)$ and fixed $i \in \{1,2,\cdots, n\}$,
\begin{equation}\label{eq213}
 \int_\Omega  {\langle A\left( {x,{\nabla _H}u} \right),{\nabla _H}\varphi \rangle dx}  + \int_\Omega  {{X_{n+i}}u \cdot {X_i}\varphi dx}  - \int_\Omega  {{X_{i}}u \cdot {X_{n+i}}\varphi dx}  = 0,
\end{equation}
then Theorems \ref{Th11}-\ref{Th13} are true for the weak solution $u \in H{W^{1,{\tilde q}}}\left( \Omega  \right)$ to \eqref{eq12}.
\end{remark}

This paper is organized as follows. In Section 2, we introduce the related knowledge of the Heisenberg Group $\mathbb{H}^n$ and some known lemmas for difference quotients. In Section 3, we first obtain an iterative relation for $2 <p<\infty$ of fractional difference quotients in the vertical direction by using the energy estimate, then prove \eqref{eq18} by this iteration relation, and \eqref{eq19} by giving another iterative relation for $1<p\le2$ of fractional difference quotients in the vertical direction. In Section 4, we use the energy estimate, first-order difference quotients in the horizontal direction and Theorem \ref{Th11} to prove Theorem \ref{Th12}. Section 5 is devoted to the proof of Theorem \ref{Th13} by using the Nirenberg difference on $\mathbb{H}^n$ and Theorem \ref{Th11}.

\section{Preliminaries}
\subsection{The Heisenberg group $\mathbb{H}^n$}

The Euclidean space ${\mathbb{R}^{2n + 1}},\;n \ge 1$ with the group multiplication
\begin{equation}\label{eq21}
x \circ y = \left( {{x_1} + {y_1},{x_2} + {y_2}, \cdots ,{x_{2n}} + {y_{2n}},t + s + \frac{1}{2}\sum\limits_{i = 1}^n {\left( {{x_i}{y_{n + i}} - {x_{n + i}}{y_i}} \right)} } \right),
\end{equation}
where $x = \left( {{x_1},{x_2}, \cdots ,{x_{2n}},t} \right),\;y = \left( {{y_1},{y_2}, \cdots ,{y_{2n}},s} \right) \in {\mathbb{R}^{2n + 1}},$ leads to the Heisenberg group $\mathbb{H}^n$. The left invariant vector fields on $\mathbb{H}^n$ are of the form
\begin{equation}\label{eq22}
{X_i} = {\partial _{{x_i}}} - \frac{{{x_{n + i}}}}{2}{\partial _t},\;{X_{n + i}} = {\partial _{{x_{n + i}}}} + \frac{{{x_i}}}{2}{\partial _t},\;\;1 \le i \le n,
\end{equation}
and a non-trivial commutator on $\mathbb{H}^n$ is
\[T = {\partial _t} = \left[ {{X_i},{X_{n + i}}} \right] = {X_i}{X_{n + i}} - {X_{n + i}}{X_i},\;\;1 \le i \le n.\]
We call that ${X_1},{X_2}, \cdots ,{X_{2n}}$ are the horizontal vector fields and $T$ the vertical vector field on $\mathbb{H}^n$. Denote the horizontal gradient of a smooth function $u$ on $\mathbb{H}^n$ by
\[{\nabla _H}u = \left( {{X_1}u,{X_2}u, \cdots ,{X_{2n}}u} \right).\]
The Haar measure in $\mathbb{H}^n$ is equivalent to the Lebesgue measure in ${\mathbb{R}^{2n + 1}}$. The Lebesgue measure of a measurable set $E \subset {\mathbb{H}^n}$ is denoted by $\left| E \right|.$

The Carnot-Carath\`{e}odary metric (C-C metric) between two points in $\mathbb{H}^n$ is the shortest length of the horizontal curve joining them, denoted by $d$. The ball induced by the C-C metric is
\[{B_R}\left( x \right) = \left\{ {y \in {\mathbb{H}^n}:d\left( {y,x} \right) < R} \right\}.\]
For $x = \left( {{x_1},{x_2}, \cdots ,{x_{2n}},t} \right)\in \mathbb{H}^n,$ its module is defined as
\[{\left\| x \right\|_{{\mathbb{H}^n}}} = {\left( {{{\left( {\sum\limits_{i = 1}^{2n} {{x_i}^2} } \right)}^2} + {t^2}} \right)^{\frac{1}{4}}}.\]
The C-C metric $d$ is equivalent to the Kor\`{a}nyi metric
\[d\left( {x,y} \right) = {\left\| {{x^{ - 1}} \circ y} \right\|_{{\mathbb{H}^n}}}.\]

If $Z$ is a left invariant vector field on $\mathbb{H}^n$, then for some $z = \left( {{z_1},{z_2}, \cdots ,{z_{2n + 1}}} \right) = \left( {z',{z_{2n + 1}}} \right)\in \mathbb{H}^n,$  we write
\[Z = \sum\limits_{l = 1}^{2n} {{z_l}{X_l}}  + {z_{2n + 1}}T.\]
The exponential mapping in canonical coordinates is defined as
\[{{\mathop{\rm e}\nolimits} ^Z} = z.\]
By \eqref{eq21}, it follows the Baker-Campbell-Hausdorff formula: if $Z$ and $Y$ are left invariant vector fields with components $z$ and $y$, then
\begin{equation}\label{eq23}
{{\mathop{\rm e}\nolimits} ^Z}{{\mathop{\rm e}\nolimits} ^Y} = \left( {z',{z_{2n + 1}}} \right) \circ \left( {y',{y_{2n + 1}}} \right) = {{\mathop{\rm e}\nolimits} ^{Z + Y + \frac{1}{2}\left[ {Z,Y} \right]}}.
\end{equation}

For $s \in \mathbb{R}\backslash \left\{ 0 \right\}$, the first order Nirenberg difference of the function $v$ along the left invariant vector field $Z$-direction is defined as
\begin{equation}\label{eq24}
{\Delta _{Z,s}}v\left( x \right) = v\left( {x{e^{sZ}}} \right) - v\left( x \right),
\end{equation}
and the second order Nirenberg difference is defined as
\begin{equation}\label{eq25}
 \Delta _{Z,s}^2v\left( x \right) = v\left( {x{e^{sZ}}} \right) + v\left( {x{e^{ - sZ}}} \right) - 2v\left( x \right).
\end{equation}
Then $\alpha $ ($0 < \alpha  \le 1$) order right difference quotient and left difference quotient (it is called the fractional right difference quotient and fractional left difference quotient if $0 < \alpha  < 1$) of the function $v$ along the left invariant vector field $Z$-direction are defined as
\begin{equation}\label{eq26}
 {D_{Z,s,\alpha }}v\left( x \right) = \frac{{v\left( {x{e^{sZ}}} \right) - v\left( x \right)}}{{{{\left| s \right|}^{ \alpha }}}}
\end{equation}
and
\begin{equation}\label{eq27}
{D_{Z, - s,\alpha }}v\left( x \right) = \frac{{v\left( {x{e^{ - sZ}}} \right) - v\left( x \right)}}{{ - {{\left| s \right|}^\alpha }}},
\end{equation}
respectively, where $s>0$. The quadratic difference quotient of $v$ is defined as
\begin{equation}\label{eq28}
{D_{Z, - s,\alpha }}{D_{Z,s,\beta }}v\left( x \right) = {D_{Z,s,\beta }}{D_{Z, - s,\alpha }}v\left( x \right) = \frac{{\Delta _{Z,s}^2v\left( x \right)}}{{{{\left| s \right|}^{\alpha  + \beta }}}}.
\end{equation}
For any $f,g \in L_{loc}^1\left( \Omega  \right)$, it implies
\begin{equation}\label{eq116}
\int_\Omega  {f\left( x \right){D_{Z,s,\alpha }}g\left( x \right)dx}  = -\int_\Omega  {{D_{Z, - s,\alpha }}f\left( x \right)g\left( x \right)dx} .
\end{equation}
For more details, see \cite{CL97} and \cite{DA04}.

For $1 \le p < \infty $ and $\Omega  \subset {\mathbb{H}^n},$ the horizontal Sobolev space $H{W^{k,p}}\left( \Omega  \right)$ is defined as
\[H{W^{k,p}}\left( \Omega  \right) = \left\{ {u \in {L^p}\left( \Omega  \right):{\nabla _H}u \in {L^p}\left( \Omega  \right),\nabla _H^2u \in {L^p}\left( \Omega  \right), \cdots ,\nabla _H^ku \in {L^p}\left( \Omega  \right)} \right\},\]
which is a Banach space under the norm
\[{\left\| u \right\|_{H{W^{k,p}}\left( \Omega  \right)}} = {\left\| u \right\|_{{L^p}\left( \Omega  \right)}} + \sum\limits_{m = 1}^k {{{\left\| {\nabla _H^mu} \right\|}_{{L^p}\left( \Omega  \right)}}} .\]
The local horizontal Sobolev space $HW_{loc}^{k,p}\left( \Omega  \right)$ is defined as
\[HW_{loc}^{k,p}\left( \Omega  \right): = \left\{ {u:u \in H{W^{k,p}}\left( {\Omega '} \right),\forall \Omega ' \subset  \subset \Omega } \right\}\]
and the space $HW_0^{k,p}\left( \Omega  \right)$ is the closure of $C_0^\infty \left( \Omega  \right)$ in $H{W^{k,p}}\left( \Omega  \right)$.

\subsection{Some known lemmas}

\begin{lemma}[Campbell-Hausdorff formula, \cite{MZZ07}]\label{Le21}
 Let $Z$ and $Y$ be two left invariant vector fields, $v \in H{W^{1,p}}\left( \Omega  \right)$ and $\left[ {Y,Z} \right]v \in L_{loc}^p\left( \Omega  \right),\;1 \le p < \infty $. Denote $\tilde v: = v\left( {x{e^{sZ}}} \right)$, then $Y\tilde v \in L_{loc}^p\left( \Omega  \right)$ and for $x,x{e^Z} \in \Omega $, it holds
\begin{equation}\label{eq29}
Y\left( {v\left( {x{e^{sZ}}} \right)} \right) = Y\tilde v\left( x \right) = Yv\left( {x{e^{sZ}}} \right) + s\left[ {Y,Z} \right]v\left( {x{e^{sZ}}} \right).
\end{equation}
Moreover, it follows for $s \ne 0$,
\begin{equation}\label{eq210}
Y\left( {{\Delta _{Z,s}}v\left( x \right)} \right) = {\Delta _{Z,s}}\left( {Yv} \right)\left( x \right) + s\left[ {Y,Z} \right]v\left( {x{e^{sZ}}} \right).
\end{equation}
\end{lemma}

\begin{lemma}[\cite{CL97}]\label{Le22}
Let $\Omega  \subset {\mathbb{H}^n}$ be an open set, $K \subset \Omega $ a
compact set, $Z$ a left invariant vector field and $v \in L_{loc}^p\left( \Omega  \right)$ for $1 \le p < \infty $. If there
exist the positive constants $h$ and $c$ such that
\[\mathop {\sup }\limits_{0 < \left| s \right| < h} \int_K {{{\left| {{D_{Z,s,1}}v\left( x \right)} \right|}^p}dx}  \le {c^p},\]
then
\[Zv \in {L^p}\left( K \right)\; \hbox{and}\;{\left\| {Zv} \right\|_{{L^p}\left( K \right)}} \le c.\]
Conversely, if $Zv \in {L^p}\left( K \right)$, then for some $h > 0$,
\[\mathop {\sup }\limits_{0 < \left| s \right| < h} \int_K {{{\left| {{D_{Z,s,1}}v\left( x \right)} \right|}^p}dx}  \le {\left( {2{{\left\| {Zv} \right\|}_{{L^p}\left( K \right)}}} \right)^p}.\]
\end{lemma}

\begin{lemma}[\cite{CL97}]\label{Le23}
 Let $\Omega  \subset {\mathbb{H}^n}$ be an open set, and $v \in HW_{loc}^{1,p}\left( \Omega  \right)$, $1 \le p < \infty $. For ${x_0} \in \Omega ,\;r > 0$ with $B\left( {{x_0},3r} \right) \subset \Omega $, there exists a positive constant $c$, such that
 \[\int_{B\left( {{x_0},r} \right)} {{{\left| {{D_{T,s,\frac{1}{2}}}v\left( x \right)} \right|}^p}dx}  \le c\int_{B\left( {{x_0},2r} \right)} {\left( {{{\left| v \right|}^p} + {{\left| {{\nabla _H}v} \right|}^p}} \right)dx} .\]
\end{lemma}

\begin{lemma}[\cite{DA04}]\label{Le24}
Let $v \in {L^p}\left( \Omega  \right),\;{\rm{1}} < p < \infty ,\;\;\alpha  > 0,\;h > 0,\;0 \le M < \infty $, and $Z$ be a left invariant vector field. If
\begin{equation}\label{eq211}
\mathop {\sup }\limits_{0 < \left| s \right| \le h} \frac{{{{\left\| {\Delta _{Z,s}^2v\left( x \right)} \right\|}_{{L^p}\left( \Omega  \right)}}}}{{{{\left| s \right|}^\alpha }}} \le M,
\end{equation}
then there exist positive constants $c$ and $h'$ such that
\[\mathop {\sup }\limits_{0 < \left| s \right| \le h'} \frac{{{{\left\| {{\Delta _{Z,s}}v\left( x \right)} \right\|}_{{L^p}\left( \Omega  \right)}}}}{{{{\left| s \right|}^\beta }}} \le c\left( {{{\left\| v \right\|}_{{L^p}\left( \Omega  \right)}} + M} \right),\]
where
\[\beta  = \alpha \;\hbox{as}\;0 < \alpha  < 1;\;\;\beta  \in \left( {0,1} \right)\;\hbox{as}\;\alpha  = 1;\;\;\beta=1\;\hbox{as}\;\alpha  > 1.\]
\end{lemma}

\begin{remark}\label{Re25}
 Lemma \ref{Le22} and Lemma \ref{Le24} imply that if $v$ has a compact support $K$ and \eqref{eq211} holds with $\alpha  > 1$, then $Zv \in {L^p}\left( K \right)$, $1 < p < \infty $.
\end{remark}

For vectors $\nu  = \left( {{\nu _1}, \cdots ,{\nu _{2n}}} \right)$ and $\mu  = \left( {{\mu _1}, \cdots ,{\mu _N}} \right)$, we define
\[\nu  \otimes \mu  = {\left( {{\nu _i}{\mu _j}} \right)_{N \times 2n}}.\]

\section{Proof of Theorem \ref{Th11}}

To prove Theorem \ref{Th11}, we first derive an iterative relation (Lemma \ref{Le31} below) of the fractional difference quotients of the weak solution $u$ to system \eqref{eq11} in the vertical direction as $2 < p < \infty $, and then ${L^p}\left( {2 < p < 4} \right)$ integrability of the vertical derivative $Tu$ of $u$ by this iterative relation. As a by-product, we also give ${L^p}(4 \le p < \infty)$ integrability of the fractional difference quotients of $u$ in the vertical direction, see Remark \ref{Re32}. Next, we prove ${L^p}(1 < p \le 2)$ integrability of $Tu$ by deriving another iterative relations (\eqref{eq324} and \eqref{eq324} below) of fractional difference quotients of $u$ in the vertical direction. By combining the results for $2 < p < 4$ and $1 < p \le 2$, the proof of Theorem \ref{Th11} is completed.

\begin{lemma}\label{Le31}
 Let $2 < p < \infty ,$ $u \in H{W^{1,{\tilde q}}}\left( \Omega  \right)$ be a weak solution to system \eqref{eq11}
with the conditions \eqref{eq13}-\eqref{eq17}. For ${x_0} \in \Omega ,$ $r > 0$ with $B\left( {{x_0},3r} \right) \subset \Omega $, assume that there exist the constants $c > 0$, $0 < h < 1$ and $\alpha  \in \left[ {0,\frac{1}{2}} \right)$ such that
\begin{align}\label{eq31}
 & \mathop {\sup }\limits_{0 < \left| s \right| \le h} \int_{B\left( {{x_0},r} \right)} {{{\left| {{D_{T,s,\frac{1}{2} + \alpha }}u\left( x \right)} \right|}^p}dx}\nonumber\\
   \le&  c \left[ {\int_{B\left( {{x_0},2r} \right)} {\left( {{{\left( {1 + {{\left| {{\nabla _H}u} \right|}^2}} \right)}^{\frac{{mp\left( {2q - 2-p} \right)}}{{2[{m\left( {p - 2} \right) - 2p}]}}}} + {{\left| u \right|}^p}} \right)dx}  + 1} \right].
\end{align}
If $\frac{{1 + 2\alpha }}{p} < \frac{1}{2}$, then for $c > 0$ and $h > 0$, it holds
\begin{align}\label{eq32}
 & \mathop {\sup }\limits_{0 < \left| s \right| \le h} \int_{B\left( {{x_0},\frac{r}{2}} \right)} {{{\left| {{D_{T,s,\frac{1}{2} + \frac{1}{p} + \frac{{2\alpha }}{p}}}u\left( x \right)} \right|}^p}dx}  \nonumber\\
  \le &c\kappa \left[ {\int_{B\left( {{x_0},2r} \right)} {\left( {{{\left( {1 + {{\left| {{\nabla _H}u} \right|}^2}} \right)}^{\frac{{mp\left( {2q - 2-p} \right)}}{{2[{m\left( {p - 2} \right) - 2p}]}}}} + {{\left| u \right|}^p}} \right)dx}  + 1} \right].
\end{align}
If $\frac{{1 + 2\alpha }}{p} \ge \frac{1}{2}$, then it yields
\begin{equation}\label{eq33}
  \int_{B\left( {{x_0},\frac{r}{4}} \right)} {{{\left| {Tu\left( x \right)} \right|}^p}dx}  \le c\kappa^2 \left[ {\int_{B\left( {{x_0},2r} \right)} {\left( {{{\left( {1 + {{\left| {{\nabla _H}u} \right|}^2}} \right)}^{\frac{{mp\left( {2q - 2-p} \right)}}{{2[{m\left( {p - 2} \right) - 2p}]}}}} + {{\left| u \right|}^p}} \right)dx}  + 1} \right],
\end{equation}
where $\kappa  = \left( {1 + {{\left\| b \right\|}_{{L^m}\left( \Omega  \right)}} + {{\left\| k \right\|}_{{L^m}\left( \Omega  \right)}} + {{\left\| {{a^{ - 1}}} \right\|}_{{L^\infty }\left( \Omega  \right)}}\left( {\left\| k \right\|_{{L^m}\left( \Omega  \right)}^2 + \left\| b \right\|_{{L^m}\left( \Omega  \right)}^2} \right)} \right)$.
\end{lemma}

\textbf{Proof.} Write $\gamma  = \frac{1}{2} + \alpha $ and let $\eta $ be a cut-off function between $B\left( {{x_0},\frac{r}{2}} \right)$ and $B\left( {{x_0},r} \right)$ with $\left| {{\nabla _H}\eta } \right| \le c$ and $\left| {T\eta } \right| \le c$. We take the test function in \eqref{eq212}
\[\varphi  = {D_{T, - s,\gamma }}\left( {{\eta ^2}{D_{T,s,\gamma }}u\left( x \right)} \right)\]
to gain
\begin{equation}\label{eq34}
  \int_\Omega  {\langle{D_{T,s,\gamma }}\left( {A\left( {x,{\nabla _H}u} \right)} \right), {\nabla _H}\left( {{\eta ^2}{D_{T,s,\gamma }}u\left( x \right)} \right)\rangle dx}  = 0.
\end{equation}
For the sake of brevity, let us denote
\[{{\rm A}_1}: = \frac{{A\left( {x,{\nabla _H}u\left( {x{e^{sT}}} \right)} \right) - A\left( {x,{\nabla _H}u\left( x \right)} \right)}}{{{{\left| s \right|}^\gamma }}}\]
and
\[{{\rm A}_2}: = \frac{{A\left( {x{e^{sT}},{\nabla _H}u\left( {x{e^{sT}}} \right)} \right) - A\left( {x,{\nabla _H}u\left( {x{e^{sT}}} \right)} \right)}}{{{{\left| s \right|}^\gamma }}},\]
then
\[{D_{T,s,\gamma }}\left( {A\left( {x,{\nabla _H}u\left( x \right)} \right)} \right) = {{\rm A}_1} + {{\rm A}_2}.\]
Since
\[{\nabla _H}\left( {{\eta ^2}{D_{T,s,\gamma }}u\left( x \right)} \right) = {\eta ^2}{D_{T,s,\gamma }}{\nabla _H}u\left( x \right) + 2\eta {\nabla _H}\eta  \otimes {D_{T,s,\gamma }}u\left( x \right),\]
it sees that \eqref{eq34} becomes
\begin{align}\label{eq35}
   0 =& \int_\Omega  {{\eta ^2}\langle{{\rm A}_1} , {D_{T,s,\gamma }}{\nabla _H}u\left( x \right)\rangle dx}  + 2\int_\Omega  {\eta \langle{{\rm A}_1}, {{\nabla _H}\eta  \otimes {D_{T,s,\gamma }}u\left( x \right)} \rangle dx} \nonumber \\
   &  + \int_\Omega  {{\eta ^2}\langle{{\rm A}_2}, {D_{T,s,\gamma }}{\nabla _H}u\left( x \right)\rangle dx}  + 2\int_\Omega  {\eta \langle{{\rm A}_2} , {{\nabla _H}\eta  \otimes {D_{T,s,\gamma }}u\left( x \right)} \rangle dx} \nonumber\\
    = :&{I_1} + {I_2} + {I_3} + {I_4}.
\end{align}
Let us estimate ${I_1},\;{I_2},\;{I_3}$ and ${I_4}$, respectively.

Using \eqref{eq13}, it derives
\begin{align}\label{eq36}
   {I_1} &= \int_\Omega  {{\eta ^2}\langle{{\rm A}_1} , {D_{T,s,\gamma }}{\nabla _H}u\left( x \right)\rangle dx} \nonumber\\
  & = {\left| s \right|^{ - 2\gamma }}\int_\Omega  {{\eta ^2}\left\langle {A\left( {x,{\nabla _H}u\left( {x{e^{sT}}} \right)} \right) - A\left( {x,{\nabla _H}u\left( x \right)} \right),{\Delta _{T,s}}{\nabla _H}u\left( x \right)} \right\rangle dx}  \nonumber\\
&\ge {\left| s \right|^{ - 2\gamma }}\int_{B\left( {{x_0},r} \right)} {a\left( x \right){\eta ^2}{{\left( {1 + {{\left| {{\nabla _H}u} \right|}^2} + {{\left| {{\nabla _H}u\left( {x{e^{sT}}} \right)} \right|}^2}} \right)}^{\frac{{p - 2}}{2}}}{{\left| {{\Delta _{T,s}}{\nabla _H}u} \right|}^2}dx} \nonumber\\
&= \int_{B\left( {{x_0},r} \right)} {a\left( x \right){\eta ^2}{{\left( {1 + {{\left| {{\nabla _H}u} \right|}^2} + {{\left| {{\nabla _H}u\left( {x{e^{sT}}} \right)} \right|}^2}} \right)}^{\frac{{p - 2}}{2}}}{{\left| {{D_{T,s,\gamma }}{\nabla _H}u} \right|}^2}dx} .
\end{align}

We apply \eqref{eq14}, Young's inequality and H\"{o}lder's inequality($\frac{2}{m} + \frac{{m - 2}}{m} = 1$) to get
\begin{align*}
  {I_2} = &  2\int_\Omega  {\eta \langle{{\rm A}_1} , ,{{\nabla _H}\eta  \otimes {D_{T,s,\gamma }}u\left( x \right)} \rangle dx}\\
   \le&  2\int_\Omega  {\eta \left| {{\nabla _H}\eta } \right|b\left( x \right){{\left( {1 + {{\left| {{\nabla _H}u} \right|}^2} + {{\left| {{\nabla _H}u\left( {x{e^{sT}}} \right)} \right|}^2}} \right)}^{\frac{{q - 2}}{2}}}\left| {{D_{T,s,\gamma }}{\nabla _H}u} \right|\left| {{D_{T,s,\gamma }}u} \right|dx}\\
     = &2\int_\Omega  {\left[ {a{{\left( x \right)}^{\frac{1}{2}}}\eta {{\left( {1 + {{\left| {{\nabla _H}u} \right|}^2} + {{\left| {{\nabla _H}u\left( {x{e^{sT}}} \right)} \right|}^2}} \right)}^{\frac{{p - 2}}{4}}}\left| {{D_{T,s,\gamma }}{\nabla _H}u} \right|} \right] \cdot } \\
&\;\;\;\;\left[ {a{{\left( x \right)}^{ - \frac{1}{2}}}b\left( x \right)\left| {{\nabla _H}\eta } \right|{{\left( {1 + {{\left| {{\nabla _H}u} \right|}^2} + {{\left| {{\nabla _H}u\left( {x{e^{sT}}} \right)} \right|}^2}} \right)}^{\frac{{2q - 2 - p}}{4}}}\left| {{D_{T,s,\gamma }}u} \right|} \right]dx\\
    \le &\varepsilon \int_{B\left( {{x_0},r} \right)} {a\left( x \right){\eta ^2}{{\left( {1 + {{\left| {{\nabla _H}u} \right|}^2} + {{\left| {{\nabla _H}u\left( {x{e^{sT}}} \right)} \right|}^2}} \right)}^{\frac{{p - 2}}{2}}}{{\left| {{D_{T,s,\gamma }}{\nabla _H}u} \right|}^2}dx} \\
& + {c_\varepsilon }\int_{B\left( {{x_0},r} \right)} {\frac{{b{{\left( x \right)}^2}}}{{a\left( x \right)}}{{\left| {{\nabla _H}\eta } \right|}^2}{{\left( {1 + {{\left| {{\nabla _H}u} \right|}^2} + {{\left| {{\nabla _H}u\left( {x{e^{sT}}} \right)} \right|}^2}} \right)}^{\frac{{2q - 2 - p}}{2}}}{{\left| {{D_{T,s,\gamma }}u} \right|}^2}dx} \\
 \le& \varepsilon \int_{B\left( {{x_0},r} \right)} {a\left( x \right){\eta ^2}{{\left( {1 + {{\left| {{\nabla _H}u} \right|}^2} + {{\left| {{\nabla _H}u\left( {x{e^{sT}}} \right)} \right|}^2}} \right)}^{\frac{{p - 2}}{2}}}{{\left| {{D_{T,s,\gamma }}{\nabla _H}u} \right|}^2}dx} \\
 & + {c_\varepsilon }{\left\| {{a^{ - 1}}} \right\|_{{L^\infty }\left( \Omega  \right)}}\left\| b \right\|_{{L^m}\left( \Omega  \right)}^2{\left( {\int_{B\left( {{x_0},r} \right)} {{{\left( {1 + {{\left| {{\nabla _H}u} \right|}^2} + \left| {{\nabla _H}u\left( {x{e^{sT}}} \right)} \right|} \right)}^{\frac{{m\left( {2q - 2 - p} \right)}}{{2\left( {m - 2} \right)}}}}{{\left| {{D_{T,s,\gamma }}u} \right|}^{\frac{{2m}}{{m - 2}}}}dx} } \right)^{\frac{{m - 2}}{m}}},
\end{align*}
and then use Young's inequality twice(the exponents satisfy $\frac{m-2}{m} + \frac{{2}}{m} = 1$ and $\frac{{m\left( {p - 2} \right) - 2p}}{{p\left( {m - 2} \right)}} + \frac{{2m}}{{p\left( {m - 2} \right)}} = 1$, respectively) and \eqref{eq31} to conclude
\begin{align}\label{eq37}
   {I_2} \le& \varepsilon \int_{B\left( {{x_0},r} \right)} {a\left( x \right){\eta ^2}{{\left( {1 + {{\left| {{\nabla _H}u} \right|}^2} + {{\left| {{\nabla _H}u\left( {x{e^{sT}}} \right)} \right|}^2}} \right)}^{\frac{{p - 2}}{2}}}{{\left| {{D_{T,s,\gamma }}{\nabla _H}u} \right|}^2}dx} \nonumber \\
   &  + {c_\varepsilon }{\left\| {{a^{ - 1}}} \right\|_{{L^\infty }\left( \Omega  \right)}}\left\| b \right\|_{{L^m}\left( \Omega  \right)}^2\left( {\int_{B\left( {{x_0},r} \right)} {{{\left( {1 + {{\left| {{\nabla _H}u} \right|}^2} + \left| {{\nabla _H}u\left( {x{e^{sT}}} \right)} \right|} \right)}^{\frac{{m\left( {2q - 2 - p} \right)}}{{2\left( {m - 2} \right)}}}}{{\left| {{D_{T,s,\gamma }}u} \right|}^{\frac{{2m}}{{m - 2}}}}dx}  + 1} \right)\nonumber \\
   \le& \varepsilon \int_{B\left( {{x_0},r} \right)} {a\left( x \right){\eta ^2}{{\left( {1 + {{\left| {{\nabla _H}u} \right|}^2} + {{\left| {{\nabla _H}u\left( {x{e^{sT}}} \right)} \right|}^2}} \right)}^{\frac{{p - 2}}{2}}}{{\left| {{D_{T,s,\gamma }}{\nabla _H}u} \right|}^2}dx}  \nonumber \\
   &+ {c_\varepsilon }{\left\| {{a^{ - 1}}} \right\|_{{L^\infty }\left( \Omega  \right)}}\left\| b \right\|_{{L^m}\left( \Omega  \right)}^2\cdot\nonumber \\
   &\;\;\;\;\left[ {\int_{B\left( {{x_0},r} \right)} {\left( {{{\left( {1 + {{\left| {{\nabla _H}u} \right|}^2} + {{\left| {{\nabla _H}u\left( {x{e^{sT}}} \right)} \right|}^2}} \right)}^{\frac{{pm\left( {2q - 2 - p} \right)}}{{2\left[ {m\left( {p - 2} \right) - 2p} \right]}}}} + {{\left| {{D_{T,s,\gamma }}u} \right|}^p}} \right)dx}  + 1} \right]\nonumber \\
    \le& \varepsilon \int_{B\left( {{x_0},r} \right)} {a\left( x \right){\eta ^2}{{\left( {1 + {{\left| {{\nabla _H}u} \right|}^2} + {{\left| {{\nabla _H}u\left( {x{e^{sT}}} \right)} \right|}^2}} \right)}^{\frac{{p - 2}}{2}}}{{\left| {{D_{T,s,\gamma }}{\nabla _H}u} \right|}^2}dx} \nonumber \\
    &+ {c_\varepsilon }{\left\| {{a^{ - 1}}} \right\|_{{L^\infty }\left( \Omega  \right)}}\left\| b \right\|_{{L^m}\left( \Omega  \right)}^2\left[ {\int_{B\left( {{x_0},2r} \right)} {\left( {{{\left( {1 + {{\left| {{\nabla _H}u} \right|}^2}} \right)}^{\frac{{pm\left( {2q - 2 - p} \right)}}{{2\left[ {m\left( {p - 2} \right) - 2p} \right]}}}} + {{\left| u \right|}^p}} \right)dx}  + 1} \right].
\end{align}

By using \eqref{eq15}, $\gamma  < 1$, Young's inequality twice(the exponents satisfy $\frac{1}{2} + \frac{1}{2} = 1$ and $\frac{m-2}{m} + \frac{2}{m} = 1$, respectively) and H\"{o}lder's inequality($\frac{2}{m} + \frac{{m - 2}}{m} = 1$), it infers
\begin{align}\label{eq38}
   {I_3} =& \int_\Omega  {{\eta ^2}\langle{{\rm A}_2}, {D_{T,s,\gamma }}{\nabla _H}u\left( x \right)\rangle dx} \nonumber \\
   \le&  {\left| s \right|^{1 - \gamma }}\int_{B\left( {{x_0},r} \right)} {k\left( x \right){\eta ^2}{{\left( {1 + {{\left| {{\nabla _H}u\left( {x{e^{sT}}} \right)} \right|}^2}} \right)}^{\frac{{q - 1}}{2}}}\left| {{D_{T,s,\gamma }}{\nabla _H}u} \right|dx}\nonumber \\
   \le& c\int_{B\left( {{x_0},r} \right)} {k\left( x \right){\eta ^2}{{\left( {1 + {{\left| {{\nabla _H}u\left( {x{e^{sT}}} \right)} \right|}^2}} \right)}^{\frac{{q - 1}}{2}}}\left| {{D_{T,s,\gamma }}{\nabla _H}u} \right|dx} \nonumber \\
=& c\int_\Omega  {\left[ {a{{\left( x \right)}^{\frac{1}{2}}}\eta {{\left( {1 + {{\left| {{\nabla _H}u} \right|}^2} + {{\left| {{\nabla _H}u\left( {x{e^{sT}}} \right)} \right|}^2}} \right)}^{\frac{{p - 2}}{4}}}\left| {{D_{T,s,\gamma }}{\nabla _H}u} \right|} \right] \cdot } \nonumber \\
&\;\;\;\;\left[ {a{{\left( x \right)}^{ - \frac{1}{2}}}k\left( x \right)\eta {{\left( {1 + {{\left| {{\nabla _H}u} \right|}^2} + {{\left| {{\nabla _H}u\left( {x{e^{sT}}} \right)} \right|}^2}} \right)}^{\frac{{2q - p}}{4}}}} \right]dx\nonumber \\
   \le& \varepsilon \int_{B\left( {{x_0},r} \right)} {a\left( x \right){\eta ^2}{{\left( {1 + {{\left| {{\nabla _H}u} \right|}^2} + {{\left| {{\nabla _H}u\left( {x{e^{sT}}} \right)} \right|}^2}} \right)}^{\frac{{p - 2}}{2}}}{{\left| {{D_{T,s,\gamma }}{\nabla _H}u} \right|}^2}dx} \nonumber \\
   & + {c_\varepsilon }\int_{B\left( {{x_0},r} \right)} {\frac{{k{{\left( x \right)}^2}}}{{a\left( x \right)}}{\eta ^2}{{\left( {1 + {{\left| {{\nabla _H}u} \right|}^2} + {{\left| {{\nabla _H}u\left( {x{e^{sT}}} \right)} \right|}^2}} \right)}^{\frac{{2q - p}}{2}}}dx} \nonumber \\
   \le &\varepsilon \int_{B\left( {{x_0},r} \right)} {a\left( x \right){\eta ^2}{{\left( {1 + {{\left| {{\nabla _H}u} \right|}^2} + {{\left| {{\nabla _H}u\left( {x{e^{sT}}} \right)} \right|}^2}} \right)}^{\frac{{p - 2}}{2}}}{{\left| {{D_{T,s,\gamma }}{\nabla _H}u} \right|}^2}dx} \nonumber \\
   & + {c_\varepsilon }{\left\| {{a^{ - 1}}} \right\|_{{L^\infty }\left( \Omega  \right)}}\left\| k \right\|_{{L^m}\left( \Omega  \right)}^2{\left( {\int_{B\left( {{x_0},2r} \right)} {{{\left( {1 + {{\left| {{\nabla _H}u} \right|}^2}} \right)}^{\frac{{m\left( {2q - p} \right)}}{{2\left( {m - 2} \right)}}}}dx} } \right)^{\frac{{m - 2}}{m}}}\nonumber \\
   \le &\varepsilon \int_{B\left( {{x_0},r} \right)} {a\left( x \right){\eta ^2}{{\left( {1 + {{\left| {{\nabla _H}u} \right|}^2} + {{\left| {{\nabla _H}u\left( {x{e^{sT}}} \right)} \right|}^2}} \right)}^{\frac{{p - 2}}{2}}}{{\left| {{D_{T,s,\gamma }}{\nabla _H}u} \right|}^2}dx} \nonumber \\
   & + {c_\varepsilon }{\left\| {{a^{ - 1}}} \right\|_{{L^\infty }\left( \Omega  \right)}}\left\| k \right\|_{{L^m}\left( \Omega  \right)}^2\left[ {\int_{B\left( {{x_0},2r} \right)} {{{\left( {1 + {{\left| {{\nabla _H}u} \right|}^2}} \right)}^{\frac{{m\left( {2q - p} \right)}}{{2\left( {m - 2} \right)}}}}dx}  + 1} \right].
\end{align}

Applying \eqref{eq15}, $\gamma  < 1$, H\"{o}lder's inequality($\frac{1}{m} + \frac{{m - 1}}{m} = 1$), Young's inequality twice(the exponents satisfy $\frac{m-1}{m} + \frac{1}{m} = 1$ and $\frac{{m\left( {p - 1} \right) - p}}{{p\left( {m - 1} \right)}} + \frac{m}{{p\left( {m - 1} \right)}} = 1$, respectively) and \eqref{eq31}, it derives
\begin{align}\label{eq39}
   {I_4} =& 2\int_\Omega  {\eta \langle{{\rm A}_2} ,  {{\nabla _H}\eta  \otimes {D_{T,s,\gamma }}u\left( x \right)} \rangle dx} \nonumber\\
   \le&  {\left| s \right|^{1 - \gamma }}\int_{B\left( {{x_0},r} \right)} {k\left( x \right)\eta \left| {{\nabla _H}\eta } \right|{{\left( {1 + {{\left| {{\nabla _H}u\left( {x{e^{sT}}} \right)} \right|}^2}} \right)}^{\frac{{q - 1}}{2}}}\left| {{D_{T,s,\gamma }}u} \right|dx} \nonumber\\
   \le &c{\left\| k \right\|_{{L^m}\left( \Omega  \right)}}{\left( {\int_{B\left( {{x_0},r} \right)} {{{\left( {{{\left( {1 + {{\left| {{\nabla _H}u\left( {x{e^{sT}}} \right)} \right|}^2}} \right)}^{\frac{{q - 1}}{2}}}\left| {{D_{T,s,\gamma }}u} \right|} \right)}^{\frac{m}{{m - 1}}}}dx} } \right)^{\frac{{m - 1}}{m}}} \nonumber\\
    \le& c{\left\| k \right\|_{{L^m}\left( \Omega  \right)}}\left[ {\int_{B\left( {{x_0},r} \right)} {{{\left( {{{\left( {1 + {{\left| {{\nabla _H}u\left( {x{e^{sT}}} \right)} \right|}^2}} \right)}^{\frac{{q - 1}}{2}}}\left| {{D_{T,s,\gamma }}u} \right|} \right)}^{\frac{m}{{m - 1}}}}dx}  + 1} \right]\nonumber\\
    \le& c{\left\| k \right\|_{{L^m}\left( \Omega  \right)}}\left[ {\int_{B\left( {{x_0},r} \right)} {\left( {{{\left( {1 + {{\left| {{\nabla _H}u\left( {x{e^{sT}}} \right)} \right|}^2}} \right)}^{\frac{{pm\left( {q - 1} \right)}}{{2\left[ {m\left( {p - 1} \right) - p} \right]}}}} + {{\left| {{D_{T,s,\gamma }}u} \right|}^p}} \right)dx}  + 1} \right]\nonumber\\
     \le&c{\left\| k \right\|_{{L^m}\left( \Omega  \right)}}\left[ {\int_{B\left( {{x_0},2r} \right)} {\left( {{{\left( {1 + {{\left| {{\nabla _H}u} \right|}^2}} \right)}^{\frac{{pm\left( {q - 1} \right)}}{{2\left[ {m\left( {p - 1} \right) - p} \right]}}}} + {{\left( {1 + {{\left| {{\nabla _H}u} \right|}^2}} \right)}^{\frac{{pm\left( {2q - 2 - p} \right)}}{{2\left[ {m\left( {p - 2} \right) - 2p} \right]}}}} + {{\left| u \right|}^p}} \right)dx}  + 1} \right]\nonumber\\
  \le& c{\left\| k \right\|_{{L^m}\left( \Omega  \right)}}\left[ {\int_{B\left( {{x_0},2r} \right)} {\left( {{{\left( {1 + {{\left| {{\nabla _H}u} \right|}^2}} \right)}^{\frac{{pm\left( {2q - 2 - p} \right)}}{{2\left[ {m\left( {p - 2} \right) - 2p} \right]}}}} + {{\left| u \right|}^p}} \right)dx}  + 1} \right].
\end{align}
Substituting \eqref{eq36}-\eqref{eq39} into \eqref{eq35} and taking $\varepsilon $ small enough, we have
\begin{align}\label{eq310}
   & \int_{B\left( {{x_0},r} \right)} {a\left( x \right){\eta ^2}{{\left( {1 + {{\left| {{\nabla _H}u} \right|}^2} + {{\left| {{\nabla _H}u\left( {x{e^{sT}}} \right)} \right|}^2}} \right)}^{\frac{{p - 2}}{2}}}{{\left| {{D_{T,s,\gamma }}{\nabla _H}u} \right|}^2}dx}\nonumber\\
  \le & c\left( {{{\left\| {{a^{ - 1}}} \right\|}_{{L^\infty }\left( \Omega  \right)}}\left\| b \right\|_{{L^m}\left( \Omega  \right)}^2 + {{\left\| k \right\|}_{{L^m}\left( \Omega  \right)}}} \right)\left[ {\int_{B\left( {{x_0},2r} \right)} {\left( {{{\left( {1 + {{\left| {{\nabla _H}u} \right|}^2}} \right)}^{\frac{{pm\left( {2q - 2 - p} \right)}}{{2\left[ {m\left( {p - 2} \right) - 2p} \right]}}}} + {{\left| u \right|}^p}} \right)dx}  + 1} \right]\nonumber\\
  & + c{\left\| {{a^{ - 1}}} \right\|_{{L^\infty }\left( \Omega  \right)}}\left\| k \right\|_{{L^m}\left( \Omega  \right)}^2\left[ {\int_{B\left( {{x_0},2r} \right)} {{{\left( {1 + {{\left| {{\nabla _H}u} \right|}^2}} \right)}^{\frac{{m\left( {2q  - p} \right)}}{{2\left( {m - 2} \right)}}}}dx}  + 1} \right]\nonumber\\
   \le& c\kappa \left[ {\int_{B\left( {{x_0},2r} \right)} {\left( {{{\left( {1 + {{\left| {{\nabla _H}u} \right|}^2}} \right)}^{\frac{{mp\left( {2q - 2-p} \right)}}{2[{m\left( {p - 2} \right) - 2p}]}}} + {{\left| u \right|}^p}} \right)dx}  + 1} \right]\nonumber\\
   : =& \tilde M.
\end{align}

For convenience, we simply write $c\tilde M$ as $\tilde M$ in the sequel. Using \eqref{eq26}, $${\left| {{D_{T,s,\gamma }}{\nabla _H}u} \right|^p} = {\left| {{D_{T,s,\gamma }}{\nabla _H}u} \right|^{p - 2}}\cdot{\left| {{D_{T,s,\gamma }}{\nabla _H}u} \right|^2}$$
and
\[\left| {{{\left| s \right|}^\gamma }{D_{T,s,\gamma }}{\nabla _H}u} \right| \le \sqrt 2 {\left( {1 + {{\left| {{\nabla _H}u} \right|}^2} + {{\left| {{\nabla _H}u\left( {x{e^{sT}}} \right)} \right|}^2}} \right)^{\frac{1}{2}}},\]
it deduces from \eqref{eq310} that
\begin{align*}
   & \int_{B\left( {{x_0},r} \right)} {a\left( x \right){\eta ^2}{{\left| {{D_{T,s,\frac{2}{p}\gamma }}{\nabla _H}u} \right|}^p}dx}  = \int_{B\left( {{x_0},r} \right)} {a\left( x \right){\eta ^2}{{\left| s \right|}^{\left( {p - 2} \right)\gamma }}{{\left| {{D_{T,s,\gamma }}{\nabla _H}u} \right|}^p}dx}  \\
   \le& c \int_{B\left( {{x_0},r} \right)} {a\left( x \right){\eta ^2}{{\left( {1 + {{\left| {{\nabla _H}u} \right|}^2} + {{\left| {{\nabla _H}u\left( {x{e^{sT}}} \right)} \right|}^2}} \right)}^{\frac{{p - 2}}{2}}}{{\left| {{D_{T,s,\gamma }}{\nabla _H}u} \right|}^2}dx}  \le \tilde M,
\end{align*}
so
\begin{align}\label{eq311}
   \int_{B\left( {{x_0},r} \right)} {{\eta ^2}{{\left| {{D_{T,s,\frac{2}{p}\gamma }}{\nabla _H}u} \right|}^p}dx}  &= \int_{B\left( {{x_0},r} \right)} {a{{\left( x \right)}^{ - 1}}a\left( x \right){\eta ^2}{{\left| {{D_{T,s,\frac{2}{p}\gamma }}{\nabla _H}u} \right|}^p}dx} \nonumber \\
   &  \le c\int_{B\left( {{x_0},r} \right)} {a\left( x \right){\eta ^2}{{\left| {{D_{T,s,\frac{2}{p}\gamma }}{\nabla _H}u} \right|}^p}dx}  \le \tilde M.
\end{align}
We note the relation
\begin{align*}
   {D_{T,s,\gamma }}{\nabla _H}\left( {{\eta ^2}u} \right) =& {D_{T,s,\gamma }}{\nabla _H}\left( {{\eta ^2}} \right)u\left( {x{e^{sT}}} \right) + {\nabla _H}\left( {{\eta ^2}} \right){D_{T,s,\gamma }}u \\
   & + {D_{T,s,\gamma }}{\eta ^2}{\nabla _H}\left( u \right)\left( {x{e^{sT}}} \right) + {\eta ^2}{D_{T,s,\gamma }}{\nabla _H}u,
\end{align*}
\eqref{eq311} and \eqref{eq31} to show
\begin{align}\label{eq312}
   & \int_{B\left( {{x_0},r} \right)} {{{\left| {{D_{T,s,\frac{2}{p}\gamma }}{\nabla _H}\left( {{\eta ^2}u} \right)} \right|}^p}dx} \nonumber \\
  \le &  c\int_{B\left( {{x_0},r} \right)} {{{\left| {{D_{T,s,\frac{2}{p}\gamma }}{\nabla _H}\left( {{\eta ^2}} \right)u\left( {x{e^{sT}}} \right)} \right|}^p}dx}  + c\int_{B\left( {{x_0},r} \right)} {{{\left| {{\nabla _H}\left( {{\eta ^2}} \right){D_{T,s,\frac{2}{p}\gamma }}u} \right|}^p}dx}\nonumber \\
  & + c\int_{B\left( {{x_0},r} \right)} {{{\left| {{D_{T,s,\frac{2}{p}\gamma }}{\eta ^2}{\nabla _H}\left( u \right)\left( {x{e^{sT}}} \right)} \right|}^p}dx}  + c\int_{B\left( {{x_0},r} \right)} {{{\left| {{\eta ^2}{D_{T,s,\frac{2}{p}\gamma }}{\nabla _H}u} \right|}^p}dx} \nonumber \\
   = &c\int_{B\left( {{x_0},r} \right)} {{{\left| {{D_{T,s,\frac{2}{p}\gamma }}{\nabla _H}\left( {{\eta ^2}} \right)u\left( {x{e^{sT}}} \right)} \right|}^p}dx}  + c{\left| s \right|^{\left( {p - 2} \right)\gamma }}\int_{B\left( {{x_0},r} \right)} {{{\left| {{\nabla _H}\left( {{\eta ^2}} \right){D_{T,s,\gamma }}u} \right|}^p}dx} \nonumber \\
& + c\int_{B\left( {{x_0},r} \right)} {{{\left| {{D_{T,s,\frac{2}{p}\gamma }}{\eta ^2}{\nabla _H}\left( u \right)\left( {x{e^{sT}}} \right)} \right|}^p}dx}  + c\int_{B\left( {{x_0},r} \right)} {{{\left| {{\eta ^2}{D_{T,s,\frac{2}{p}\gamma }}{\nabla _H}u} \right|}^p}dx} \nonumber \\
 \le& c\int_{B\left( {{x_0},2r} \right)} {{{\left| u \right|}^p}dx}  + c\int_{B\left( {{x_0},r} \right)} {{{\left| {{D_{T,s,\gamma }}u} \right|}^p}dx}  + c\int_{B\left( {{x_0},2r} \right)} {{{\left| {{\nabla _H}u} \right|}^p}dx}  + \tilde M \nonumber \\
  \le& \tilde M.
\end{align}
Then it gets from Lemma \ref{Le23}, \eqref{eq312} and \eqref{eq31} that
\begin{align*}
   & \int_{B\left( {{x_0},r} \right)} {{{\left| {{D_{T, - s,\frac{1}{2}}}{D_{T,s,\frac{2}{p}\gamma }}\left( {{\eta ^2}u} \right)} \right|}^p}dx}  \\
   \le&  c\int_{B\left( {{x_0},r} \right)} {{{\left| {{D_{T,s,\frac{2}{p}\gamma }}\left( {{\eta ^2}u} \right)} \right|}^p}dx}  + c\int_{B\left( {{x_0},r} \right)} {{{\left| {{D_{T,s,\frac{2}{p}\gamma }}{\nabla _H}\left( {{\eta ^2}u} \right)} \right|}^p}dx} \\
    =& c{\left| s \right|^{\left( {p - 2} \right)\gamma }}\int_{B\left( {{x_0},r} \right)} {{{\left| {{D_{T,s,\gamma }}\left( {{\eta ^2}u} \right)} \right|}^p}dx}  + c\int_{B\left( {{x_0},r} \right)} {{{\left| {{D_{T,s,\frac{2}{p}\gamma }}{\nabla _H}\left( {{\eta ^2}u} \right)} \right|}^p}dx} \\
    \le& c\int_{B\left( {{x_0},2r} \right)} {\left( {{{\left| {{D_{T,s,\gamma }}u} \right|}^p} + {{\left| u \right|}^p}} \right)dx}  + \tilde M\\
     \le &\tilde M.
\end{align*}
Therefore for $s$ small enough, we know that \eqref{eq28} and the above formula imply
\[\frac{{{{\left\| {\Delta _{Z,s}^2\left( {{\eta ^2}u} \right)} \right\|}_{{L^p}\left( {B\left( {{x_0},r} \right)} \right)}}}}{{{{\left| s \right|}^{\frac{1}{2} + \frac{2}{p}\gamma }}}} \le {\tilde M^{\frac{1}{p}}},\]
where $\gamma  = \frac{1}{2} + \alpha $, and then there exists $0 < h < 1$ such that
\begin{equation}\label{eq313}
 \mathop {\sup }\limits_{0 < \left| s \right| \le h} \frac{{{{\left\| {\Delta _{Z,s}^2\left( {{\eta ^2}u} \right)} \right\|}_{{L^p}\left( {B\left( {{x_0},r} \right)} \right)}}}}{{{{\left| s \right|}^{\frac{1}{2} + \frac{{1 + 2\alpha }}{p}}}}} \le {\tilde M^{\frac{1}{p}}}.
\end{equation}
If $\frac{{1 + 2\alpha }}{p} < \frac{1}{2}$, then \eqref{eq32} holds by \eqref{eq313} and Lemma \ref{Le24}; If $\frac{{1 + 2\alpha }}{p} > \frac{1}{2}$, then Remark \ref{Re25} yields $Tu \in L_{loc}^p\left( \Omega  \right)$ and \eqref{eq33} holds; If $\frac{{1 + 2\alpha }}{p} = \frac{1}{2}$, then $2 \le p < 4$ from $\alpha  \in \left[ {0,\frac{1}{2}} \right)$ and it holds
\[\mathop {\sup }\limits_{0 < \left| s \right| \le h} \int_{B\left( {{x_0},\frac{r}{2}} \right)} {{{\left| {{D_{T,s,\frac{1}{2} + \alpha '}}u\left( x \right)} \right|}^p}dx}  \le \tilde M\]
for any $\alpha '\in \left[ {0,\frac{1}{2}} \right)$ by Lemma \ref{Le24}. In particular, we can choose $\alpha ' > \frac{{p - 2}}{4}$. We take again that $\eta $ is a cut-off function between $B\left( {{x_0},\frac{r}{4}} \right)$ and $B\left( {{x_0},\frac{r}{2}} \right)$, and get by the discussion similar to \eqref{eq313} that
\[\mathop {\sup }\limits_{0 < \left| s \right| \le h} \frac{{{{\left\| {\Delta _{Z,s}^2\left( {{\eta ^2}u} \right)} \right\|}_{{L^p}\left( {B\left( {{x_0},\frac{r}{2}} \right)} \right)}}}}{{{{\left| s \right|}^{\frac{1}{2} + \frac{{1 + 2\alpha '}}{p}}}}} \le {\left( {\kappa \tilde M} \right)^{\frac{1}{p}}}\]
and $\frac{{1 + 2\alpha '}}{p} > \frac{1}{2},$ so it follows $Tu \in L_{loc}^p\left( \Omega  \right)$ from Remark \ref{Re25}, and \eqref{eq33} holds.

\begin{remark}\label{Re31}
It observes from the proof of Lemma \ref{Le31} that we can have a larger integral region on the left hand side of \eqref{eq33} as $\frac{{1 + 2\alpha }}{p} > \frac{1}{2}$, namely, we have
\[\int_{B\left( {{x_0},\frac{r}{2}} \right)} {{{\left| {Tu\left( x \right)} \right|}^p}dx}  \le c\kappa \left[ {\int_{B\left( {{x_0},2r} \right)} {\left( {{{\left( {1 + {{\left| {{\nabla _H}u} \right|}^2}} \right)}^{\frac{{mp\left( {2q - 2-p} \right)}}{2[{m\left( {p - 2} \right) - 2p}]}}} + {{\left| u \right|}^p}} \right)dx}  + 1} \right].\]
\end{remark}

\textbf{Proof of Theorem \ref{Th11}($2 < p < 4$).}

From Lemma \ref{Le23}, we see that \eqref{eq31} holds with ${\alpha _0} = 0$, so \eqref{eq32} is true with ${\alpha _1} = \frac{1}{p}$ by Lemma \ref{Le31}. Then we use ${\alpha _1} = \frac{1}{p}$ in \eqref{eq31} to gain ${\alpha _2} = \frac{1}{p} + \frac{2}{p}{\alpha _1}$ in \eqref{eq32} by Lemma \ref{Le31}. In general, if we have already found ${\alpha _1}, \cdots ,{\alpha _l}$, then
\[{\alpha _{l + 1}} = \frac{1}{p} + \frac{2}{p}{\alpha _l} =  \cdots  = \frac{1}{p}\sum\limits_{i = 0}^{l - 1} {{{\left( {\frac{2}{p}} \right)}^i}}  = \frac{1}{p}\frac{{1 - {{\left( {\frac{2}{p}} \right)}^l}}}{{1 - \frac{2}{p}}}.\]
Thus for a fixed $p > 2$, it shows $\mathop {\sup }\limits_l {\alpha _l} = \frac{1}{{p - 2}}$. Hence, for $2 < p < 4$, after $l$ iterations ($l$ may be sufficiently large), we always have ${\alpha _l} > \frac{1}{2}$, and it follows
\[\mathop {\sup }\limits_{0 < \left| s \right| \le h} \int_{B\left( {{x_0},\frac{r}{{{2^l}}}} \right)} {{{\left| {{D_{T,s,\frac{1}{2} + {\alpha _l}}}u} \right|}^p}dx}  \le c{\kappa ^l}\left[ {\int_{B\left( {{x_0},2r} \right)} {\left( {{{\left( {1 + {{\left| {{\nabla _H}u} \right|}^2}} \right)}^{\frac{{mp\left( {2q - 2 - p} \right)}}{{2\left[ {m\left( {p - 2} \right) - 2p} \right]}}}} + {{\left| u \right|}^p}} \right)dx}  + 1} \right],\]
so \eqref{eq18} holds by Remark \ref{Re31}, i.e.,
\[\int_{B\left( {{x_0},\frac{r}{{{2^{l + 1}}}}} \right)} {{{\left| {Tu\left( x \right)} \right|}^p}dx}  \le c{\kappa ^{l + 1}}\left[ {\int_{B\left( {{x_0},2r} \right)} {\left( {{{\left( {1 + {{\left| {{\nabla _H}u} \right|}^2}} \right)}^{\frac{{mp\left( {2q - 2 - p} \right)}}{{2\left[ {m\left( {p - 2} \right) - 2p} \right]}}}} + {{\left| u \right|}^p}} \right)dx}  + 1} \right],\]
which implies $Tu \in L_{loc}^p\left( \Omega  \right)$.

\begin{remark}\label{Re32}
For $4\le p < \infty$, we can see from the above iteration that
\[\mathop {\sup }\limits_{0 < \left| s \right| \le h} \int_{B\left( {{x_0},\frac{r}{{{2^l}}}} \right)} {{{\left| {{D_{T,s,\frac{1}{2} + \alpha '}}u} \right|}^p}dx}  \le c\kappa^l \left[ {\int_{B\left( {{x_0},2r} \right)} {\left( {{{\left( {1 + {{\left| {{\nabla _H}u} \right|}^2}} \right)}^{\frac{{mp\left( {2q - 2-p} \right)}}{2[{m\left( {p - 2} \right) - 2p}]}}} + {{\left| u \right|}^p}} \right)dx}  + 1} \right],\]
where $\alpha '$  is smaller than $\frac{1}{{p - 2}}$ but arbitrarily close to $\frac{1}{{p - 2}}$, and $l$ is an iteration number.
\end{remark}

\begin{remark}\label{Re33}
From the proof process of Lemma \ref{Le31}, we have that a statement similarly to Lemma \ref{Le31} also holds for the weak solution $u \in H{W^{1,{\tilde q}}}\left( \Omega  \right)$ to system \eqref{eq12}. In fact, in considering system \eqref{eq12}, an additional term $\int_{B\left( {{x_{0}},r}\right) }{{X_{i}}u\ \cdot {X_{n+i}}\varphi dx}-\int_{B\left( {{x_{0}},r}\right) \ }{{X_{n+i}}u\ \cdot {X_{i}}\varphi dx}$ in \eqref{eq35} appears, denoted by ${I_5}$. Using the integration by parts, Young's inequality and \eqref{eq31}, it shows
\begin{align*}
{I_{5}}=& -\int_{B\left( {{x_{0}},r}\right) }{{D_{T,s,\gamma }X_{i}}u\cdot {%
X_{n+i}}\left( {{\eta ^{2}}{D_{T,s,\gamma }}u}\right) dx}+\int_{B\left( {{%
x_{0}},r}\right) }{{D_{T,s,\gamma }X_{n+i}}u\cdot {X_{i}}\left( {{\eta ^{2}}{%
D_{T,s,\gamma }}u}\right) dx} \\
=& -\int_{B\left( {{x_{0}},r}\right) }{{\eta ^{2}{X_{i}}}}\left( {{%
D_{T,s,\gamma }}u}\right) {\cdot {X_{n+i}}\left( {{D_{T,s,\gamma }}u}\right)
dx-2}\int_{B\left( {{x_{0}},r}\right) }{{\eta {X_{n+i}\eta \cdot X_{i}}}}%
\left( {{D_{T,s,\gamma }}u}\right) {\cdot {{D_{T,s,\gamma }}u}dx} \\
& +\int_{B\left( {{x_{0}},r}\right) }{{\eta ^{2}{X_{n+i}}}}\left( {{%
D_{T,s,\gamma }}u}\right) {\cdot {{X_{i}}}\left( {{D_{T,s,\gamma }}u}\right)
dx+2}\int_{B\left( {{x_{0}},r}\right) }{{\eta {{X_{i}}\eta \cdot {{X_{n+i}}}}%
}}\left( {{D_{T,s,\gamma }}u}\right) {\cdot {{D_{T,s,\gamma }}u}dx} \\
=&- {2}\int_{B\left( {{x_{0}},r}\right) }{{\eta {X_{n+i}\eta \cdot X_{i}}}}%
\left( {{D_{T,s,\gamma }}u}\right) {\cdot {{D_{T,s,\gamma }}u}dx+2}%
\int_{B\left( {{x_{0}},r}\right) }{{\eta {{X_{i}}\eta \cdot {{X_{n+i}}}}}}%
\left( {{D_{T,s,\gamma }}u}\right) {\cdot {{D_{T,s,\gamma }}u}dx} \\
=& -\int_{B\left( {{x_{0}},r}\right) }{{\eta {X_{n+i}\eta \cdot X_{i}}}}%
\left[ \left( {{D_{T,s,\gamma }}u}\right) ^{2}\right] {dx+}\int_{B\left( {{%
x_{0}},r}\right) }{{\eta {{X_{i}}\eta \cdot {{X_{n+i}}}}}}\left[ \left( {{%
D_{T,s,\gamma }}u}\right) ^{2}\right] {dx} \\
=& \int_{B\left( {{x_{0}},r}\right) }\left[ {{{X_{i}}}}\left( {{\eta {%
X_{n+i}\eta }}}\right) {{{{-{X_{n+i}}}}}}\left( {{\eta {{X_{i}}\eta }}}%
\right) \right] \left( {{D_{T,s,\gamma }}u}\right) ^{2}{dx=}\int_{B\left( {{%
x_{0}},r}\right) }{{\eta {T\eta }}}\left( {{D_{T,s,\gamma }}u}\right) ^{2}{dx%
} \\
\leq & c\int_{B\left( {{x_{0}},r}\right) }{{{\left\vert {{D_{T,s,\gamma }}u}%
\right\vert }^{2}}dx}\leq c \int_{B\left( {{x_{0}},r}\right) }{{\left ({\left\vert {{D_{T,s,\gamma }}u}%
\right\vert }^{p}+1\right)}dx}\\
\leq& c\int_{B\left( {{x_{0}},r}\right) }{\left( {\left(
{1+{{\left\vert {{\nabla _{H}}u}\right\vert }^{2}}}\right)^{\frac{{mp\left( {2q - 2-p} \right)}}{2[{m\left( {p - 2} \right) - 2p}]}} +{{\left\vert
u\right\vert }^{p}}}\right) dx}\leq \tilde{M}.
\end{align*}
The remaining proofs are consistent with Lemma \ref{Le31}. Therefore, Theorem \ref{Th11} with
$2 < p < 4$ and Remark \ref{Re32} also holds.
\end{remark}

\textbf{Proof of Theorem \ref{Th11} ($1 < p \le 2$)}
Let $\eta $ be a cut-off function between $B\left( {{x_0},\frac{r}{2}} \right)$ and $B\left( {{x_0},r} \right)$ with $\left| {{\nabla _H}\eta } \right| \le c$ and $\left| {T\eta } \right| \le c$. Similarly to the proof of Lemma \ref{Le31}, we can get \eqref{eq35}, only changing $\alpha  \in \left[ {0,\frac{1}{2}} \right)$ into $\alpha  = 0$, and then $\gamma  = \frac{1}{2}$.  Now we estimate ${I_1},\;{I_2},\;{I_3}$ and ${I_4}$.

By \eqref{eq13},
\begin{align}\label{eq314}
   {I_1}& \ge \int_{B\left( {{x_0},r} \right)} {a\left( x \right){\eta ^2}{{\left( {1 + {{\left| {{\nabla _H}u} \right|}^2} + {{\left| {{\nabla _H}u\left( {x{e^{sT}}} \right)} \right|}^2}} \right)}^{\frac{{p - 2}}{2}}}{{\left| {{D_{T,s,\gamma }}{\nabla _H}u} \right|}^2}dx} \nonumber \\
   & \ge c\int_{B\left( {{x_0},r} \right)} {{\eta ^2}{{\left( {1 + {{\left| {{\nabla _H}u} \right|}^2} + {{\left| {{\nabla _H}u\left( {x{e^{sT}}} \right)} \right|}^2}} \right)}^{\frac{{p - 2}}{2}}}{{\left| {{D_{T,s,\gamma }}{\nabla _H}u} \right|}^2}dx} .
\end{align}
Using \eqref{eq14}, H\"{o}lder's inequality, Young's inequality and Lemma \ref{Le23}, it yields
\begin{align}\label{eq315}
 {I_2}  & \le c\int_{B\left( {{x_0},r} \right)} {b\left( x \right){{\left( {1 + {{\left| {{\nabla _H}u} \right|}^2} + {{\left| {{\nabla _H}u\left( {x{e^{sT}}} \right)} \right|}^2}} \right)}^{\frac{{q - 2}}{2}}}\left| {{D_{T,s,\gamma }}{\nabla _H}u} \right|\left| {{D_{T,s,\gamma }}u} \right|dx} \nonumber \\
   &  = c{\left| s \right|^{ - \gamma }}\int_{B\left( {{x_0},r} \right)} {b\left( x \right){{\left( {1 + {{\left| {{\nabla _H}u} \right|}^2} + {{\left| {{\nabla _H}u\left( {x{e^{sT}}} \right)} \right|}^2}} \right)}^{\frac{{q - 2}}{2}}}\left| {{\Delta _{T,s}}{\nabla _H}u} \right|\left| {{D_{T,s,\gamma }}u} \right|dx}\nonumber \\
   & \le c{\left| s \right|^{ - \gamma }}\int_{B\left( {{x_0},r} \right)} {b\left( x \right){{\left( {1 + {{\left| {{\nabla _H}u} \right|}^2} + {{\left| {{\nabla _H}u\left( {x{e^{sT}}} \right)} \right|}^2}} \right)}^{\frac{{q - 1}}{2}}}\left| {{D_{T,s,\gamma }}u} \right|dx} \nonumber \\
   & \le c{\left| s \right|^{ - \gamma }}{\left\| b \right\|_{{L^m}\left( \Omega  \right)}}{\left( {\int_{B\left( {{x_0},r} \right)} {{{\left( {1 + {{\left| {{\nabla _H}u} \right|}^2} + {{\left| {{\nabla _H}u\left( {x{e^{sT}}} \right)} \right|}^2}} \right)}^{\frac{{m\left( {q - 1} \right)}}{{2\left( {m - 1} \right)}}}}{{\left| {{D_{T,s,\gamma }}u} \right|}^{\frac{m}{{m - 1}}}}dx} } \right)^{\frac{{m - 1}}{m}}}\nonumber \\
   & \le c{\left| s \right|^{ - \gamma }}{\left\| b \right\|_{{L^m}\left( \Omega  \right)}}\left( {\int_{B\left( {{x_0},r} \right)} {{{\left( {1 + {{\left| {{\nabla _H}u} \right|}^2} + {{\left| {{\nabla _H}u\left( {x{e^{sT}}} \right)} \right|}^2}} \right)}^{\frac{{m\left( {q - 1} \right)}}{{2\left( {m - 1} \right)}}}}{{\left| {{D_{T,s,\gamma }}u} \right|}^{\frac{m}{{m - 1}}}}dx}  + 1} \right)\nonumber \\
   & \le c{\left| s \right|^{ - \gamma }}{\left\| b \right\|_{{L^m}\left( \Omega  \right)}}\left( {\int_{B\left( {{x_0},r} \right)} {\left( {{{\left( {1 + {{\left| {{\nabla _H}u} \right|}^2} + {{\left| {{\nabla _H}u\left( {x{e^{sT}}} \right)} \right|}^2}} \right)}^{\frac{{mp\left( {q - 1} \right)}}{{2\left[ {m\left( {p - 1} \right) - p} \right]}}}} + {{\left| {{D_{T,s,\gamma }}u} \right|}^p}} \right)dx}  + 1} \right)\nonumber \\
  & \le c{\left| s \right|^{ - \gamma }}{\left\| b \right\|_{{L^m}\left( \Omega  \right)}}\left( {\int_{B\left( {{x_0},2r} \right)} {\left( {{{\left( {1 + {{\left| {{\nabla _H}u} \right|}^2}} \right)}^{\frac{{mp\left( {q - 1} \right)}}{{2\left[ {m\left( {p - 1} \right) - p} \right]}}}} + {{\left| {{\nabla _H}u} \right|}^p} + {{\left| u \right|}^p}} \right)dx}  + 1} \right)\nonumber \\
  &\le c{\left| s \right|^{ - \gamma }}{\left\| b \right\|_{{L^m}\left( \Omega  \right)}}\left( {\int_{B\left( {{x_0},2r} \right)} {\left( {{{\left( {1 + {{\left| {{\nabla _H}u} \right|}^2}} \right)}^{\frac{{mp\left( {q - 1} \right)}}{{2\left[ {m\left( {p - 1} \right) - p} \right]}}}} + {{\left| u \right|}^p}} \right)dx}  + 1} \right).
\end{align}
One applies \eqref{eq15}, H\"{o}lder's inequality and Young's inequality to have
\begin{align}\label{eq316}
   {I_3}& \le {\left| s \right|^{1 - \gamma }}\int_{B\left( {{x_0},r} \right)} {k\left( x \right){\eta ^2}{{\left( {1 + {{\left| {{\nabla _H}u\left( {x{e^{sT}}} \right)} \right|}^2}} \right)}^{\frac{{q - 1}}{2}}}\left| {{D_{T,s,\gamma }}{\nabla _H}u} \right|dx}  \nonumber\\
   & \le c{\left| s \right|^{1 - 2\gamma }}\int_{B\left( {{x_0},r} \right)} {k\left( x \right){\eta ^2}{{\left( {1 + {{\left| {{\nabla _H}u\left( {x{e^{sT}}} \right)} \right|}^2}} \right)}^{\frac{q}{2}}}dx}  \nonumber\\
  & \le c{\left| s \right|^{1 - 2\gamma }}{\left\| k \right\|_{{L^m}\left( \Omega  \right)}}{\left( {\int_{B\left( {{x_0},r} \right)} {{{\left( {1 + {{\left| {{\nabla _H}u\left( {x{e^{sT}}} \right)} \right|}^2}} \right)}^{\frac{{mq}}{{2\left( {m - 1} \right)}}}}dx} } \right)^{\frac{{m - 1}}{m}}} \nonumber\\
  & \le c{\left| s \right|^{1 - 2\gamma }}{\left\| k \right\|_{{L^m}\left( \Omega  \right)}}\left( {\int_{B\left( {{x_0},2r} \right)} {{{\left( {1 + {{\left| {{\nabla _H}u} \right|}^2}} \right)}^{\frac{{mq}}{{2\left( {m - 1} \right)}}}}dx}  + 1} \right).
\end{align}
Applying \eqref{eq15}, H\"{o}lder's inequality and Young's inequality, it gives
\begin{align}\label{eq317}
   {I_4}& \le {\left| s \right|^{1 - \gamma }}\int_{B\left( {{x_0},r} \right)} {k\left( x \right)\eta \left| {{\nabla _H}\eta } \right|{{\left( {1 + {{\left| {{\nabla _H}u\left( {x{e^{sT}}} \right)} \right|}^2}} \right)}^{\frac{{q - 1}}{2}}}\left| {{D_{T,s,\gamma }}u} \right|dx} \nonumber \\
   &  \le c{\left| s \right|^{1 - 2\gamma }}\int_{B\left( {{x_0},r} \right)} {k\left( x \right){{\left( {1 + {{\left| {{\nabla _H}u\left( {x{e^{sT}}} \right)} \right|}^2}} \right)}^{\frac{{q - 1}}{2}}}\left| {{\Delta _{T,s}}u} \right|dx} \nonumber \\
  & \le c{\left| s \right|^{1 - 2\gamma }}{\left\| k \right\|_{{L^m}\left( \Omega  \right)}}{\left( {\int_{B\left( {{x_0},r} \right)} {{{\left( {1 + {{\left| {{\nabla _H}u\left( {x{e^{sT}}} \right)} \right|}^2}} \right)}^{\frac{{m\left( {q - 1} \right)}}{{2\left( {m - 1} \right)}}}}{{\left| {{\Delta _{T,s}}u} \right|}^{\frac{m}{{m - 1}}}}dx} } \right)^{\frac{{m - 1}}{m}}} \nonumber \\
   &\le c{\left| s \right|^{1 - 2\gamma }}{\left\| k \right\|_{{L^m}\left( \Omega  \right)}}\left( {\int_{B\left( {{x_0},r} \right)} {{{\left( {1 + {{\left| {{\nabla _H}u\left( {x{e^{sT}}} \right)} \right|}^2}} \right)}^{\frac{{m\left( {q - 1} \right)}}{{2\left( {m - 1} \right)}}}}{{\left| {{\Delta _{T,s}}u} \right|}^{\frac{m}{{m - 1}}}}dx}  + 1} \right)\nonumber \\
   & \le c{\left| s \right|^{1 - 2\gamma }}{\left\| k \right\|_{{L^m}\left( \Omega  \right)}}\left( {\int_{B\left( {{x_0},2r} \right)} {\left( {{{\left( {1 + {{\left| {{\nabla _H}u} \right|}^2}} \right)}^{\frac{{mp\left( {q - 1} \right)}}{{2\left[ {m\left( {p - 1} \right) - p} \right]}}}} + {{\left| u \right|}^p}} \right)dx}  + 1} \right).
\end{align}
Substituting \eqref{eq314}-\eqref{eq317} into \eqref{eq35}, it follows
\begin{align}\label{eq318}
   & \int_{B\left( {{x_0},r} \right)} {{\eta ^2}{{\left( {1 + {{\left| {{\nabla _H}u} \right|}^2} + {{\left| {{\nabla _H}u\left( {x{e^{sT}}} \right)} \right|}^2}} \right)}^{\frac{{p - 2}}{2}}}{{\left| {{\Delta _{T,s}}{\nabla _H}u} \right|}^2}dx} \nonumber\\
   &  \le c{\left| s \right|^\gamma }\left( {{{\left\| b \right\|}_{{L^m}\left( \Omega  \right)}} + {{\left\| k \right\|}_{{L^m}\left( \Omega  \right)}}} \right)\left( {\int_{B\left( {{x_0},2r} \right)} {{{\left( {1 + {{\left| {{\nabla _H}u} \right|}^2}} \right)}^{\frac{{mp\left( {q - 1} \right)}}{{2\left[ {m\left( {p - 1} \right) - p} \right]}}}} + {{\left( {1 + {{\left| {{\nabla _H}u} \right|}^2}} \right)}^{\frac{{mq}}{{2\left( {m - 1} \right)}}}} + {{\left| u \right|}^p}dx}  + 1} \right)\nonumber\\
   & \le c{\left| s \right|^\gamma }\left( {{{\left\| b \right\|}_{{L^m}\left( \Omega  \right)}} + {{\left\| k \right\|}_{{L^m}\left( \Omega  \right)}}} \right)\left( {\int_{B\left( {{x_0},2r} \right)} {\left( {{{\left( {1 + {{\left| {{\nabla _H}u} \right|}^2}} \right)}^{\frac{{mp(q-1)}}{{2\left[ {m\left( {p - 1} \right) - p} \right]}}}} + {{\left| u \right|}^p}} \right)dx}  + 1} \right).
\end{align}
Noting
\begin{align*}
   & {\left( {1 + {{\left| {{\nabla _H}u} \right|}^2} + {{\left| {{\nabla _H}u\left( {x{e^{sT}}} \right)} \right|}^2}} \right)^{\frac{p}{2}}} \\
  = &  {\left( {1 + {{\left| {{\nabla _H}u} \right|}^2} + {{\left| {{\nabla _H}u\left( {x{e^{sT}}} \right)} \right|}^2}} \right)^{\frac{p}{2} - 1}}\left( {1 + {{\left| {{\nabla _H}u} \right|}^2} + {{\left| {{\nabla _H}u\left( {x{e^{sT}}} \right)} \right|}^2}} \right)\\
   \le &{\left( {1 + {{\left| {{\nabla _H}u} \right|}^2} + {{\left| {{\nabla _H}u\left( {x{e^{sT}}} \right)} \right|}^2}} \right)^{\frac{p}{2} - 1}}\left( {1 + {{\left| {{\nabla _H}u} \right|}^2} + {{\left( {\left| {{\nabla _H}u} \right| + \left| {{\Delta _{T,s}}{\nabla _H}u} \right|} \right)}^2}} \right)\\
   \le& {\left( {1 + {{\left| {{\nabla _H}u} \right|}^2} + {{\left| {{\nabla _H}u\left( {x{e^{sT}}} \right)} \right|}^2}} \right)^{\frac{p}{2} - 1}}\left( {1 + {{\left| {{\nabla _H}u} \right|}^2} + 2\left( {{{\left| {{\nabla _H}u} \right|}^2} + {{\left| {{\Delta _{T,s}}{\nabla _H}u} \right|}^2}} \right)} \right)\\
  \le &3{\left( {1 + {{\left| {{\nabla _H}u} \right|}^2} + {{\left| {{\nabla _H}u\left( {x{e^{sT}}} \right)} \right|}^2}} \right)^{\frac{p}{2} - 1}}\left( {1 + {{\left| {{\nabla _H}u} \right|}^2} + {{\left| {{\Delta _{T,s}}{\nabla _H}u} \right|}^2}} \right)\\
   \le& 3{\left( {1 + {{\left| {{\nabla _H}u} \right|}^2}} \right)^{\frac{p}{2}}} + 3{\left( {1 + {{\left| {{\nabla _H}u} \right|}^2} + {{\left| {{\nabla _H}u\left( {x{e^{sT}}} \right)} \right|}^2}} \right)^{\frac{p}{2} - 1}}{\left| {{\Delta _{T,s}}{\nabla _H}u} \right|^2},
\end{align*}
we obtain from \eqref{eq318} that
\begin{align}\label{eq319}
   & \int_{B\left( {{x_0},r} \right)} {{\eta ^2}{{\left( {1 + {{\left| {{\nabla _H}u} \right|}^2} + {{\left| {{\nabla _H}u\left( {x{e^{sT}}} \right)} \right|}^2}} \right)}^{\frac{p}{2}}}dx} \nonumber \\
   \le&  3\int_{B\left( {{x_0},r} \right)} {{\eta ^2}{{\left( {1 + {{\left| {{\nabla _H}u} \right|}^2}} \right)}^{\frac{p}{2}}}dx}  \nonumber \\
   &+ 3\int_{B\left( {{x_0},r} \right)} {{\eta ^2}{{\left( {1 + {{\left| {{\nabla _H}u} \right|}^2} + {{\left| {{\nabla _H}u\left( {x{e^{sT}}} \right)} \right|}^2}} \right)}^{\frac{{p - 2}}{2}}}{{\left| {{\Delta _{T,s}}{\nabla _H}u} \right|}^2}dx}\nonumber \\
  \le& c\left( {{{\left\| b \right\|}_{{L^m}\left( \Omega  \right)}} + {{\left\| k \right\|}_{{L^m}\left( \Omega  \right)}} + 1} \right)\left( {\int_{B\left( {{x_0},2r} \right)} {\left( {{{\left( {1 + {{\left| {{\nabla _H}u} \right|}^2}} \right)}^{\frac{{mp(q-1)}}{{2\left[ {m\left( {p - 1} \right) - p} \right]}}}} + {{\left| u \right|}^p}} \right)dx}  + 1} \right),
\end{align}
and then apply H\"{o}lder's inequality, \eqref{eq318} and \eqref{eq319} to deduce
\begin{align*}
   & \int_{B\left( {{x_0},r} \right)} {{\eta ^2}{{\left| {{\Delta _{T,s}}{\nabla _H}u} \right|}^p}dx}  \\
   =&  \int_{B\left( {{x_0},r} \right)} {{{\left( {{\eta ^2}{{\left( {1 + {{\left| {{\nabla _H}u} \right|}^2} + {{\left| {{\nabla _H}u\left( {x{e^{sT}}} \right)} \right|}^2}} \right)}^{\frac{{p - 2}}{2}}}{{\left| {{\Delta _{T,s}}{\nabla _H}u} \right|}^2}} \right)}^{\frac{p}{2}}}}  \cdot \\
   &{\left( {{\eta ^{\frac{4}{p}}}\left( {1 + {{\left| {{\nabla _H}u} \right|}^2} + {{\left| {{\nabla _H}u\left( {x{e^{sT}}} \right)} \right|}^2}} \right)} \right)^{\frac{{2 - p}}{2} \cdot \frac{p}{2}}}dx\\
  \le& {\left( {\int_{B\left( {{x_0},r} \right)} {{\eta ^2}{{\left( {1 + {{\left| {{\nabla _H}u} \right|}^2} + {{\left| {{\nabla _H}u\left( {x{e^{sT}}} \right)} \right|}^2}} \right)}^{\frac{{p - 2}}{2}}}{{\left| {{\Delta _{T,s}}{\nabla _H}u} \right|}^2}dx} } \right)^{\frac{p}{2}}} \cdot \\
  &{\left( {\int_{B\left( {{x_0},r} \right)} {{\eta ^2}{{\left( {1 + {{\left| {{\nabla _H}u} \right|}^2} + {{\left| {{\nabla _H}u\left( {x{e^{sT}}} \right)} \right|}^2}} \right)}^{\frac{p}{2}}}dx} } \right)^{\frac{{2 - p}}{2}}}\\
   \le& c{\left| s \right|^{\frac{{p\gamma }}{2}}}\left( {{{\left\| b \right\|}_{{L^m}\left( \Omega  \right)}} + {{\left\| k \right\|}_{{L^m}\left( \Omega  \right)}} + 1} \right)\left( {\int_{B\left( {{x_0},2r} \right)} {\left( {{{\left( {1 + {{\left| {{\nabla _H}u} \right|}^2}} \right)}^{\frac{{mp(q-1)}}{{2\left[ {m\left( {p - 1} \right) - p} \right]}}}} + {{\left| u \right|}^p}} \right)dx}  + 1} \right).
\end{align*}
Similarly to the proof of \eqref{eq312}, we see by using the above formula that
\begin{align}\label{eq320}
   & \int_{B\left( {{x_0},r} \right)} {{{\left| {{D_{T,s,\frac{\gamma }{2}}}{\nabla _H}\left( {{\eta ^2}u} \right)} \right|}^p}dx} \nonumber \\
   \le&  c\left( {{{\left\| b \right\|}_{{L^m}\left( \Omega  \right)}} + {{\left\| k \right\|}_{{L^m}\left( \Omega  \right)}} + 1} \right)\left( {\int_{B\left( {{x_0},2r} \right)} {\left( {{{\left( {1 + {{\left| {{\nabla _H}u} \right|}^2}} \right)}^{\frac{{mp(q-1)}}{{2\left[ {m\left( {p - 1} \right) - p} \right]}}}} + {{\left| u \right|}^p}} \right)dx}  + 1} \right),
\end{align}
so Lemma \ref{Le23} and \eqref{eq320} imply
\begin{align}\label{eq321}
   & \int_{B\left( {{x_0},r} \right)} {{{\left| {{D_{T, - s,\frac{1}{2}}}{D_{T,s,\frac{\gamma }{2}}}\left( {{\eta ^2}u} \right)} \right|}^p}dx} \nonumber \\
   \le&  c\left( {{{\left\| b \right\|}_{{L^m}\left( \Omega  \right)}} + {{\left\| k \right\|}_{{L^m}\left( \Omega  \right)}} + 1} \right)\left( {\int_{B\left( {{x_0},2r} \right)} {\left( {{{\left( {1 + {{\left| {{\nabla _H}u} \right|}^2}} \right)}^{\frac{{mp(q-1)}}{{2\left[ {m\left( {p - 1} \right) - p} \right]}}}} + {{\left| u \right|}^p}} \right)dx}  + 1} \right).
\end{align}
Then there exists $0 < h < 1$ such that
\begin{align}\label{eq322}
   & \mathop {\sup }\limits_{0 < \left| s \right| \le h} \frac{{{{\left\| {\Delta _{Z,s}^2\left( {{\eta ^2}u} \right)} \right\|}_{{L^p}\left( {B\left( {{x_0},r} \right)} \right)}}}}{{{{\left| s \right|}^{\frac{1}{2} + \frac{\gamma }{2}}}}}\nonumber \\
  \le & {\left[ {c\left( {{{\left\| b \right\|}_{{L^m}\left( \Omega  \right)}} + {{\left\| k \right\|}_{{L^m}\left( \Omega  \right)}} + 1} \right)\left( {\int_{B\left( {{x_0},2r} \right)} {\left( {{{\left( {1 + {{\left| {{\nabla _H}u} \right|}^2}} \right)}^{\frac{{mp(q-1)}}{{2\left[ {m\left( {p - 1} \right) - p} \right]}}}} + {{\left| u \right|}^p}} \right)dx}  + 1} \right)} \right]^{\frac{1}{p}}}.
\end{align}
Because of $\gamma  = \frac{1}{2}$, the orders of the corresponding fractional difference quotients are $\frac{1}{4}$ in \eqref{eq320} and $\frac{1}{2} + \frac{\gamma }{2} = \frac{3}{4}$ in \eqref{eq322}, respectively. Using Lemma \ref{Le24} and taking that $\eta $ is a cut-off function between $B\left( {{x_0},\frac{r}{{{2^l}}}} \right)$ and $B\left( {{x_0},\frac{r}{{{2^{l - 1}}}}} \right)$, then after the $l$ iterations, $l=1,2,\cdots,$ it gives
\begin{align}\label{eq323}
   & \int_{B\left( {{x_0},\frac{r}{{{2^{l - 1}}}}} \right)} {{{\left| {{D_{T,s,\frac{{{2^l} - 1}}{{{2^{l + 1}}}}}}{\nabla _H}\left( {{\eta ^2}u} \right)\left( x \right)} \right|}^p}dx} \nonumber \\
   \le&  c\left( {{{\left\| b \right\|}_{{L^m}\left( \Omega  \right)}} + {{\left\| k \right\|}_{{L^m}\left( \Omega  \right)}} + 1} \right)^l\left( {\int_{B\left( {{x_0},2r} \right)} {\left( {{{\left( {1 + {{\left| {{\nabla _H}u} \right|}^2}} \right)}^{\frac{{mp(q-1)}}{{2\left[ {m\left( {p - 1} \right) - p} \right]}}}} + {{\left| u \right|}^p}} \right)dx}  + 1} \right)
\end{align}
and
\begin{align}\label{eq324}
   & \mathop {\sup }\limits_{0 < \left| s \right| \le h} \frac{{{{\left\| {{\Delta _{Z,s}}\left( {{\eta ^2}u} \right)} \right\|}_{{L^p}\left( {B\left( {{x_0},\frac{r}{{{2^l}}}} \right)} \right)}}}}{{{{\left| s \right|}^{\frac{{{2^{l + 1}} - 1}}{{{2^{l + 1}}}}}}}}\nonumber \\
  \le &  {\left[ {c\left( {{{\left\| b \right\|}_{{L^m}\left( \Omega  \right)}} + {{\left\| k \right\|}_{{L^m}\left( \Omega  \right)}} + 1} \right)^l\left( {\int_{B\left( {{x_0},2r} \right)} {\left( {{{\left( {1 + {{\left| {{\nabla _H}u} \right|}^2}} \right)}^{\frac{{mp(q-1)}}{{2\left[ {m\left( {p - 1} \right) - p} \right]}}}} + {{\left| u \right|}^p}} \right)dx}  + 1} \right)} \right]^{\frac{1}{p}}}.
\end{align}

We consider now $l \in \mathbb{N}$ satisfying
\[\frac{1}{{{2^l} - 1}} < p - 1.\]
Then for
\[\iota  = \frac{{{2^l} - 1}}{{{2^{l + 1}}}}\;\; \hbox{and}\;\; \sigma  = \frac{{{2^{l + 1}} - 1}}{{{2^{l + 1}}}},\]
it follows
\[\iota \left( {p - 1} \right) + \sigma  > 1.\]
Let us take
\[\gamma  = \frac{{\iota \left( {p - 1} \right) + \sigma }}{2}\]
and return to \eqref{eq35} with a cut-off function $\eta $ between $B\left( {{x_0},\frac{r}{{{2^{l + 1}}}}} \right)$ and $B\left( {{x_0},\frac{r}{{{2^l}}}} \right)$. Then by ${\left| s \right|^{ - \iota \left( {p - 1} \right)}} \le {\left| s \right|^{ - \gamma }}$, H\"{o}lder's inequality, Young's inequality and \eqref{eq324}, one has
\begin{align*}
   {I_2}& \le \int_{B\left( {{x_0},\frac{r}{{{2^l}}}} \right)} {b\left( x \right){{\left( {1 + {{\left| {{\nabla _H}u} \right|}^2} + {{\left| {{\nabla _H}u\left( {x{e^{sT}}} \right)} \right|}^2}} \right)}^{\frac{{q - 2}}{2}}}\frac{{\left| {{\Delta _{T,s}}{\nabla _H}u} \right|}}{{{{\left| s \right|}^{\iota\left( {p - 1} \right)}}}}\frac{{\left| {{\Delta _{T,s}}u} \right|}}{{{{\left| s \right|}^\sigma}}}dx}  \\
   &  \le {\left| s \right|^{ - \gamma }}\int_{B\left( {{x_0},\frac{r}{{{2^l}}}} \right)} {b\left( x \right){{\left( {1 + {{\left| {{\nabla _H}u} \right|}^2} + {{\left| {{\nabla _H}u\left( {x{e^{sT}}} \right)} \right|}^2}} \right)}^{\frac{{q - 2}}{2}}}\left| {{\Delta _{T,s}}{\nabla _H}u} \right|\frac{{\left| {{\Delta _{T,s}}u} \right|}}{{{s^\sigma}}}dx} \\
    &\le c{\left| s \right|^{ - \gamma }}\int_{B\left( {{x_0},\frac{r}{{{2^l}}}} \right)} {b\left( x \right){{\left( {1 + {{\left| {{\nabla _H}u} \right|}^2} + {{\left| {{\nabla _H}u\left( {x{e^{sT}}} \right)} \right|}^2}} \right)}^{\frac{{q - 1}}{2}}}\frac{{\left| {{\Delta _{T,s}}u} \right|}}{{{s^\sigma}}}dx} \\
    & \le c{\left| s \right|^{ - \gamma }}{\left\| b \right\|_{{L^m}\left( \Omega  \right)}}{\left( {\int_{B\left( {{x_0},\frac{r}{{{2^l}}}} \right)} {{{\left( {1 + {{\left| {{\nabla _H}u} \right|}^2} + {{\left| {{\nabla _H}u\left( {x{e^{sT}}} \right)} \right|}^2}} \right)}^{\frac{{m\left( {q - 1} \right)}}{{2\left( {m - 1} \right)}}}}{{\left| {\frac{{\left| {{\Delta _{T,s}}u} \right|}}{{{s^\sigma}}}} \right|}^{\frac{m}{{m - 1}}}}dx} } \right)^{\frac{{m - 1}}{m}}}\\
    &\le c{\left| s \right|^{ - \gamma }}{\left\| b \right\|_{{L^m}\left( \Omega  \right)}}\left( {\int_{B\left( {{x_0},\frac{r}{{{2^l}}}} \right)} {{{\left( {1 + {{\left| {{\nabla _H}u} \right|}^2} + {{\left| {{\nabla _H}u\left( {x{e^{sT}}} \right)} \right|}^2}} \right)}^{\frac{{m\left( {q - 1} \right)}}{{2\left( {m - 1} \right)}}}}{{\left| {\frac{{\left| {{\Delta _{T,s}}u} \right|}}{{{s^\sigma}}}} \right|}^{\frac{m}{{m - 1}}}}dx}  + 1} \right)\\
   &\le c{\left| s \right|^{ - \gamma }}{\left\| b \right\|_{{L^m}\left( \Omega  \right)}}\left( {\int_{B\left( {{x_0},\frac{r}{{{2^l}}}} \right)} {\left( {{{\left( {1 + {{\left| {{\nabla _H}u} \right|}^2} + {{\left| {{\nabla _H}u\left( {x{e^{sT}}} \right)} \right|}^2}} \right)}^{\frac{{mp\left( {q - 1} \right)}}{{2\left[ {m\left( {p - 1} \right) - p} \right]}}}} + \frac{{{{\left| {{\Delta _{T,s}}u} \right|}^p}}}{{{s^{\sigma p}}}}} \right)dx}  + 1} \right)\\
   & \le c\kappa^{l+1} {\left| s \right|^{ - \gamma }}\left( {\int_{B\left( {{x_0},2r} \right)} {\left( {{{\left( {1 + {{\left| {{\nabla _H}u} \right|}^2}} \right)}^{\frac{{mp(q-1)}}{{2\left[ {m\left( {p - 1} \right) - p} \right]}}}} + {{\left| u \right|}^p}} \right)dx}  + 1} \right).
\end{align*}
By the similar discussion with \eqref{eq316}-\eqref{eq322}, it concludes
\[\mathop {\sup }\limits_{0 < \left| s \right| \le h} \frac{{{{\left\| {\Delta _{Z,s}^2\left( {{\eta ^2}u} \right)} \right\|}_{{L^p}\left( {B\left( {{x_0},\frac{r}{{{2^l}}}} \right)} \right)}}}}{{{{\left| s \right|}^{\frac{1}{2} + \frac{\gamma }{2}}}}} \le {\left[ {c\kappa^{l+1} \left( {\int_{B\left( {{x_0},2r} \right)} {\left( {{{\left( {1 + {{\left| {{\nabla _H}u} \right|}^2}} \right)}^{\frac{{mp(q-1)}}{{2\left[ {m\left( {p - 1} \right) - p} \right]}}}} + {{\left| u \right|}^p}} \right)dx}  + 1} \right)} \right]^{\frac{1}{p}}}.\]
Owing to $\frac{1}{2} + \frac{\gamma }{2} > 1$ and Remark \ref{Re25}, it means \eqref{eq19}. From \eqref{eq18} and \eqref{eq19}, it follows directly $Tu \in L_{loc}^p\left( \Omega  \right)$.

\begin{remark}\label{Re34}
From the above proof process, we see that Theorem \ref{Th11} with $1 < p \le 2$ also holds for the weak solution $u \in H{W^{1,{\tilde q}}}\left( \Omega  \right)$ to system \eqref{eq12}, just replacing ${\left| u \right|^p}$ in \eqref{eq19} with ${\left| u \right|^2}$. To be specific, when we deal with system \eqref{eq12}, an additional term $\int_{B\left( {{x_{0}},r}\right) }{{X_{i}}u\ \cdot {X_{n+i}}\varphi dx}-\int_{B\left( {{x_{0}},r}\right) \ }{{X_{n+i}}u\ \cdot {X_{i}}\varphi dx}$ will appear in \eqref{eq35}, denoted by ${I_5}$. To it, it follows by using the integration by parts and Lemma \ref{Le23} that
\begin{align*}
   {I_{5}}=& -\int_{B\left( {{x_{0}},r}\right) }{{D_{T,s,\gamma }X_{i}}u\cdot {%
X_{n+i}}\left( {{\eta ^{2}}{D_{T,s,\gamma }}u}\right) dx}+\int_{B\left( {{%
x_{0}},r}\right) }{{D_{T,s,\gamma }X_{n+i}}u\cdot {X_{i}}\left( {{\eta ^{2}}{%
D_{T,s,\gamma }}u}\right) dx} \\
=& -\int_{B\left( {{x_{0}},r}\right) }{{\eta ^{2}{X_{i}}}}\left( {{%
D_{T,s,\gamma }}u}\right) {\cdot {X_{n+i}}\left( {{D_{T,s,\gamma }}u}\right)
dx-2}\int_{B\left( {{x_{0}},r}\right) }{{\eta {X_{n+i}\eta \cdot X_{i}}}}%
\left( {{D_{T,s,\gamma }}u}\right) {\cdot {{D_{T,s,\gamma }}u}dx} \\
& +\int_{B\left( {{x_{0}},r}\right) }{{\eta ^{2}{X_{n+i}}}}\left( {{%
D_{T,s,\gamma }}u}\right) {\cdot {{X_{i}}}\left( {{D_{T,s,\gamma }}u}\right)
dx+2}\int_{B\left( {{x_{0}},r}\right) }{{\eta {{X_{i}}\eta \cdot {{X_{n+i}}}}%
}}\left( {{D_{T,s,\gamma }}u}\right) {\cdot {{D_{T,s,\gamma }}u}dx} \\
=&- {2}\int_{B\left( {{x_{0}},r}\right) }{{\eta {X_{n+i}\eta \cdot X_{i}}}}%
\left( {{D_{T,s,\gamma }}u}\right) {\cdot {{D_{T,s,\gamma }}u}dx+2}%
\int_{B\left( {{x_{0}},r}\right) }{{\eta {{X_{i}}\eta \cdot {{X_{n+i}}}}}}%
\left( {{D_{T,s,\gamma }}u}\right) {\cdot {{D_{T,s,\gamma }}u}dx} \\
=& -\int_{B\left( {{x_{0}},r}\right) }{{\eta {X_{n+i}\eta \cdot X_{i}}}}%
\left[ \left( {{D_{T,s,\gamma }}u}\right) ^{2}\right] {dx+}\int_{B\left( {{%
x_{0}},r}\right) }{{\eta {{X_{i}}\eta \cdot {{X_{n+i}}}}}}\left[ \left( {{%
D_{T,s,\gamma }}u}\right) ^{2}\right] {dx} \\
=& \int_{B\left( {{x_{0}},r}\right) }\left[ {{{X_{i}}}}\left( {{\eta {%
X_{n+i}\eta }}}\right) {{{{-{X_{n+i}}}}}}\left( {{\eta {{X_{i}}\eta }}}%
\right) \right] \left( {{D_{T,s,\gamma }}u}\right) ^{2}{dx=}\int_{B\left( {{%
x_{0}},r}\right) }{{\eta {T\eta }}}\left( {{D_{T,s,\gamma }}u}\right) ^{2}{dx%
} \\
   \le& c\int_{B\left( {{x_0},r} \right)} {{{\left| {{D_{T,s,\gamma }}u} \right|}^2}dx}  \le c\int_{B\left( {{x_0},r} \right)} {\left( {{{\left| {{\nabla _H}u} \right|}^2} + {{\left| u \right|}^2}} \right)dx} \\
\le& c\kappa \left( {\int_{B\left( {{x_0},r} \right)} {\left( {{{\left( {1 + \left| {{\nabla _H}{u^2}} \right|} \right)}^{\frac{{mp(q-1)}}{{2\left[ {m\left( {p - 1} \right) - p} \right]}}}} + {{\left| u \right|}^2}} \right)dx}  + 1} \right).
\end{align*}
The remaining proofs are consistent with the previous proof.
\end{remark}

\section{Proof of Theorem \ref{Th12}}

\textbf{Proof of Theorem \ref{Th12}.} Let ${i_0} \in \left\{ {1, \cdots ,n} \right\},\;s > 0$, and take the test function in \eqref{eq212}
\[\varphi  = {D_{{X_{{i_0}}}, - s,1}}{D_{{X_{{i_0}}},s,1}}\left( {{\eta ^4}u} \right)\left( x \right),\]
where $\eta $ is a cut-off function between $B\left( {{x_0},\frac{r}{{{2^{l + 2}}}}} \right)$ and $B\left( {{x_0},\frac{r}{{{2^{l + 1}}}}} \right)$ with $\left| {{\nabla _H}\eta } \right| \le c$ and $\left| {T\eta } \right| \le c$.

For $i \ne n + {i_0}$, it holds from the commutativity of ${X_i}$, ${D_{{X_{{i_0}}}, - s,1}}$ and ${D_{{X_{{i_0}}},s,1}}$ that
\[{X_i}\left( {{D_{{X_{{i_0}}}, - s,1}}{D_{{X_{{i_0}}},s,1}}\left( {{\eta ^4}u} \right)\left( x \right)} \right) = {D_{{X_{{i_0}}}, - s,1}}{D_{{X_{{i_0}}},s,1}}\left( {{X_i}\left( {{\eta ^4}u} \right)\left( x \right)} \right).\]
For $i = n + {i_0}$, it infers by the non-commutativity of ${X_i}$, ${D_{{X_{{i_0}}}, - s,1}}$ and ${D_{{X_{{i_0}}},s,1}}$ and Lemma \ref{Le21} that
\begin{align}\label{eq41}
  &{X_{n + {i_0}}}\left( {{D_{{X_{{i_0}}}, - s,1}}{D_{{X_{{i_0}}},s,1}}\left( {{\eta ^4}u} \right)\left( x \right)} \right)\nonumber\\
 =& {D_{{X_{{i_0}}}, - s,1}}{X_{n + {i_0}}}\left( {{D_{{X_{{i_0}}},s,1}}\left( {{\eta ^4}u} \right)\left( x \right)} \right) - \left[ {{X_{n + {i_0}}},{X_{{i_0}}}} \right]{D_{{X_{{i_0}}},s,1}}\left( {\left( {{\eta ^4}u} \right)\left( {x{e^{ - s{X_{{i_0}}}}}} \right)} \right)\nonumber\\
= &{D_{{X_{{i_0}}}, - s,1}}{D_{{X_{{i_0}}},s,1}}\left( {{X_{n + {i_0}}}\left( {{\eta ^4}u} \right)\left( x \right)} \right)\nonumber\\
& - \left[ {{D_{{X_{{i_0}}}, - s,1}}\left( {T\left( {{\eta ^4}u} \right)\left( {x{e^{s{X_{{i_0}}}}}} \right)} \right) - T{D_{{X_{{i_0}}},s,1}}\left( {\left( {{\eta ^4}u} \right)\left( {x{e^{ - s{X_{{i_0}}}}}} \right)} \right)} \right]\nonumber\\
 = &{D_{{X_{{i_0}}}, - s,1}}{D_{{X_{{i_0}}},s,1}}\left( {{X_{n + {i_0}}}\left( {{\eta ^4}u} \right)\left( x \right)} \right)\nonumber\\
 &- \left[ {{D_{{X_{{i_0}}},s,1}}\left( {T\left( {{\eta ^4}u} \right)\left( x \right)} \right) + {D_{{X_{{i_0}}}, - s,1}}\left( {T\left( {{\eta ^4}u} \right)\left( x \right)} \right)} \right],
\end{align}
so \eqref{eq212} becomes
\begin{align*}
   & \int_\Omega  {\langle A\left( {x,{\nabla _H}u} \right),{D_{{X_{{i_0}}}, - s,1}}{D_{{X_{{i_0}}},s,1}}\left( {{\nabla _H}\left( {{\eta ^4}u} \right)\left( x \right)} \right)\rangle dx}  \\
  = & \int_\Omega  {{A_{n + {i_0}}}\left( {x,{\nabla _H}u} \right)\left[ {{D_{{X_{{i_0}}},s,1}}\left( {T\left( {{\eta ^4}u} \right)\left( x \right)} \right) + {D_{{X_{{i_0}}}, - s,1}}\left( {T\left( {{\eta ^4}u} \right)\left( x \right)} \right)} \right]dx} ,
\end{align*}
The property of difference quotients \eqref{eq116} yields
\begin{align}\label{eq42}
   & \int_\Omega  {\langle{D_{{X_{{i_0}}},s,1}}\left( {A\left( {x,{\nabla _H}u} \right)} \right),{D_{{X_{{i_0}}},s,1}}\left( {{\nabla _H}\left( {{\eta ^4}u} \right)} \right)\left( x \right)\rangle dx} \nonumber \\
  = &  - \int_\Omega  {{A_{n + {i_0}}}\left( {x,{\nabla _H}u} \right)\left[ {{D_{{X_{{i_0}}},s,1}}T\left( {{\eta ^4}u} \right)\left( x \right) + {D_{{X_{{i_0}}}, - s,1}}T\left( {{\eta ^4}u} \right)\left( x \right)} \right]dx} .
\end{align}
Referring to the equality below (5.2) in \cite{DA04}, we know
\begin{align*}
  {D_{{X_{{i_0}}},s,1}}\left( {{\nabla _H}\left( {{\eta ^4}u} \right)} \right)\left( x \right) = & {D_{{X_{{i_0}}},s,1}}\left( {4{\eta ^3}{\nabla _H}\eta  \otimes u + {\eta ^4}{\nabla _H}u} \right)\left( x \right) \\
  = & 4{D_{{X_{{i_0}}},s,1}}\eta \left( x \right) \cdot \eta {\left( {x{e^{s{X_{{i_0}}}}}} \right)^2}{\nabla _H}\eta \left( {x{e^{s{X_{{i_0}}}}}} \right) \otimes u\left( {x{e^{s{X_{{i_0}}}}}} \right)\\
  & + 4\eta \left( x \right){D_{{X_{{i_0}}},s,1}}\eta \left( x \right)\eta \left( {x{e^{s{X_{{i_0}}}}}} \right){\nabla _H}\eta \left( {x{e^{s{X_{{i_0}}}}}} \right) \otimes u\left( {x{e^{s{X_{{i_0}}}}}} \right)\\
  & + 4\eta {\left( x \right)^2}{D_{{X_{{i_0}}},s,1}}\eta \left( x \right){\nabla _H}\eta \left( {x{e^{s{X_{{i_0}}}}}} \right) \otimes u\left( {x{e^{s{X_{{i_0}}}}}} \right)\\
  & + 4\eta {\left( x \right)^3}{D_{{X_{{i_0}}},s,1}}{\nabla _H}\eta \left( x \right) \otimes u\left( {x{e^{s{X_{{i_0}}}}}} \right)\\
  & + 4\eta {\left( x \right)^3}{\nabla _H}\eta \left( x \right) \otimes {D_{{X_{{i_0}}},s,1}}u\left( x \right)\\
  &+ {D_{{X_{{i_0}}},s,1}}\eta \left( x \right)\eta {\left( {x{e^{s{X_{{i_0}}}}}} \right)^3}{\nabla _H}u\left( {x{e^{s{X_{{i_0}}}}}} \right)\\
  & + \eta \left( x \right){D_{{X_{{i_0}}},s,1}}\eta \left( x \right)\eta {\left( {x{e^{s{X_{{i_0}}}}}} \right)^2}{\nabla _H}u\left( {x{e^{s{X_{{i_0}}}}}} \right)\\
  & + \eta {\left( x \right)^2}{D_{{X_{{i_0}}},s,1}}\eta \left( x \right)\eta \left( {x{e^{s{X_{{i_0}}}}}} \right){\nabla _H}u\left( {x{e^{s{X_{{i_0}}}}}} \right)\\
  & + \eta {\left( x \right)^3}{D_{{X_{{i_0}}},s,1}}\eta \left( x \right){\nabla _H}u\left( {x{e^{s{X_{{i_0}}}}}} \right)\\
  & + \eta {\left( x \right)^4}{D_{{X_{{i_0}}},s,1}}{\nabla _H}u\left( x \right),
\end{align*}
and infer by substituting the above relationship into \eqref{eq42} that
\begin{align}\label{eq43}
   {J_0}: =& \int_\Omega  {{\eta ^4}\langle{D_{{X_{{i_0}}},s,1}}\left( {A\left( {x,{\nabla _H}u} \right)} \right),{D_{{X_{{i_0}}},s,1}}{\nabla _H}u \rangle dx} \nonumber \\
= & \int_\Omega  {{D_{{X_{{i_0}}}, - s,1}}\left( {{A_{n + {i_0}}}\left( {x,{\nabla _H}u} \right)} \right)T\left( {{\eta ^4}u} \right)dx}  + \int_\Omega  {{D_{{X_{{i_0}}},s,1}}\left( {{A_{n + {i_0}}}\left( {x,{\nabla _H}u} \right)} \right)T\left( {{\eta ^4}u} \right)dx} \nonumber \\
& - 4\int_\Omega  {\langle{D_{{X_{{i_0}}},s,1}}\left( {A\left( {x,{\nabla _H}u} \right)} \right),{D_{{X_{{i_0}}},s,1}}\eta \left( x \right)  \eta {{\left( {x{e^{s{X_{{i_0}}}}}} \right)}^2}{\nabla _H}\eta \left( {x{e^{s{X_{{i_0}}}}}} \right) \otimes u\left( {x{e^{s{X_{{i_0}}}}}} \right)\rangle dx} \nonumber \\
& - 4\int_\Omega  {\langle{D_{{X_{{i_0}}},s,1}}\left( {A\left( {x,{\nabla _H}u} \right)} \right),\eta \left( x \right){D_{{X_{{i_0}}},s,1}}\eta \left( x \right)\eta \left( {x{e^{s{X_{{i_0}}}}}} \right){\nabla _H}\eta \left( {x{e^{s{X_{{i_0}}}}}} \right) \otimes u\left( {x{e^{s{X_{{i_0}}}}}} \right)\rangle dx} \nonumber \\
&- 4\int_\Omega  {\langle{D_{{X_{{i_0}}},s,1}}\left( {A\left( {x,{\nabla _H}u} \right)} \right),\eta {{\left( x \right)}^2}{D_{{X_{{i_0}}},s,1}}\eta \left( x \right){\nabla _H}\eta \left( {x{e^{s{X_{{i_0}}}}}} \right) \otimes u\left( {x{e^{s{X_{{i_0}}}}}} \right)\rangle dx} \nonumber \\
& - 4\int_\Omega  {\langle{D_{{X_{{i_0}}},s,1}}\left( {A\left( {x,{\nabla _H}u} \right)} \right),\eta {{\left( x \right)}^3}{D_{{X_{{i_0}}},s,1}}{\nabla _H}\eta \left( x \right) \otimes u\left( {x{e^{s{X_{{i_0}}}}}} \right)\rangle dx} \nonumber \\
& - 4\int_\Omega  {\langle{D_{{X_{{i_0}}},s,1}}\left( {A\left( {x,{\nabla _H}u} \right)} \right),\eta {{\left( x \right)}^3}{\nabla _H}\eta \left( x \right) \otimes {D_{{X_{{i_0}}},s,1}}u\left( x \right)\rangle dx} \nonumber \\
& - \int_\Omega  {\langle{D_{{X_{{i_0}}},s,1}}\left( {A\left( {x,{\nabla _H}u} \right)} \right),{D_{{X_{{i_0}}},s,1}}\eta \left( x \right)\eta {{\left( {x{e^{s{X_{{i_0}}}}}} \right)}^3}{\nabla _H}u\left( {x{e^{s{X_{{i_0}}}}}} \right)\rangle dx} \nonumber \\
& - \int_\Omega  {\langle{D_{{X_{{i_0}}},s,1}}\left( {A\left( {x,{\nabla _H}u} \right)} \right),\eta \left( x \right){D_{{X_{{i_0}}},s,1}}\eta \left( x \right)\eta {{\left( {x{e^{s{X_{{i_0}}}}}} \right)}^2}{\nabla _H}u\left( {x{e^{s{X_{{i_0}}}}}} \right)\rangle dx} \nonumber \\
& - \int_\Omega  {\langle{D_{{X_{{i_0}}},s,1}}\left( {A\left( {x,{\nabla _H}u} \right)} \right),\eta {{\left( x \right)}^2}{D_{{X_{{i_0}}},s,1}}\eta \left( x \right)\eta \left( {x{e^{s{X_{{i_0}}}}}} \right){\nabla _H}u\left( {x{e^{s{X_{{i_0}}}}}} \right)\rangle dx} \nonumber \\
& - \int_\Omega  {\langle{D_{{X_{{i_0}}},s,1}}\left( {A\left( {x,{\nabla _H}u} \right)} \right),\eta {{\left( x \right)}^3}{D_{{X_{{i_0}}},s,1}}\eta \left( x \right){\nabla _H}u\left( {x{e^{s{X_{{i_0}}}}}} \right)\rangle dx} \nonumber \\
: =& \sum\limits_{i = 1}^{10} {{J_i}} .
\end{align}
Next, let us estimate ${J_i}\left( {i = 0,1, \cdots ,10} \right)$. For the sake of brevity, denote
\[{B_1}: = \frac{{A\left( {x,{\nabla _H}u\left( {x{e^{s{X_{{i_0}}}}}} \right)} \right) - A\left( {x,{\nabla _H}u\left( x \right)} \right)}}{s}\]
and
\[{B_2}: = \frac{{A\left( {x{e^{s{X_{{i_0}}}}},{\nabla _H}u\left( {x{e^{s{X_{{i_0}}}}}} \right)} \right) - A\left( {x,{\nabla _H}u\left( {x{e^{s{X_{{i_0}}}}}} \right)} \right)}}{s},\]
then
\[{D_{{X_{{i_0}}},s,1}}\left( {A\left( {x,{\nabla _H}u\left( x \right)} \right)} \right) = {B_1} + {B_2}.\]
Thus we have
\[{J_0} = \int_\Omega  {{\eta ^4}\langle{B_1} , {D_{{X_{{i_0}}},s,1}}{\nabla _H}u\left( x \right)\rangle dx}  + \int_\Omega  {{\eta ^4}\langle{B_2} , {D_{{X_{{i_0}}},s,1}}{\nabla _H}u\left( x \right)\rangle dx}  = :{J_{01}} + {J_{02}}.\]
Similarly, we also have
\[{J_i} = :{J_{i1}} + {J_{i2}},\;(i = 1,2, \cdots ,10).\]

It yields from \eqref{eq13} that
\begin{equation}\label{eq44}
{J_{01}} \ge \int_\Omega  {a\left( x \right){\eta ^4}{{\left( {1 + {{\left| {{\nabla _H}u} \right|}^2} + {{\left| {{\nabla _H}u\left( {x{e^{s{X_{{i_0}}}}}} \right)} \right|}^2}} \right)}^{\frac{{p - 2}}{2}}}{{\left| {{D_{{X_{{i_0}}},s,1}}{\nabla _H}u} \right|}^2}dx}
\end{equation}
and derives by using \eqref{eq15}, Young's inequality and H\"{o}lder's inequality that
\begin{align}\label{eq45}
  {J_{02}} \le & \int_\Omega  {{\eta ^4}k\left( x \right){{\left( {1 + {{\left| {{\nabla _H}u\left( {x{e^{s{X_{{i_0}}}}}} \right)} \right|}^2}} \right)}^{\frac{{q - 1}}{2}}}\left| {{D_{{X_{{i_0}}},s,1}}{\nabla _H}u} \right|dx} \nonumber \\
  \le & \varepsilon \int_\Omega  {{\eta ^4}a\left( x \right){{\left( {1 + {{\left| {{\nabla _H}u} \right|}^2} + {{\left| {{\nabla _H}u\left( {x{e^{s{X_{{i_0}}}}}} \right)} \right|}^2}} \right)}^{\frac{{p - 2}}{2}}}{{\left| {{D_{{X_{{i_0}}},s,1}}{\nabla _H}u} \right|}^2}dx}   \nonumber \\
 & + {c_\varepsilon }\int_\Omega  {{\eta ^4}\frac{{k{{\left( x \right)}^2}}}{{a\left( x \right)}}{{\left( {1 + {{\left| {{\nabla _H}u\left( {x{e^{s{X_{{i_0}}}}}} \right)} \right|}^2}} \right)}^{\frac{{2q - p}}{2}}}dx}  \nonumber \\
\le & \varepsilon \int_\Omega  {{\eta ^4}a\left( x \right){{\left( {1 + {{\left| {{\nabla _H}u} \right|}^2} + {{\left| {{\nabla _H}u\left( {x{e^{s{X_{{i_0}}}}}} \right)} \right|}^2}} \right)}^{\frac{{p - 2}}{2}}}{{\left| {{D_{{X_{{i_0}}},s,1}}{\nabla _H}u} \right|}^2}dx}   \nonumber \\
 & + {c_\varepsilon }{\left\| {{a^{ - 1}}} \right\|_{{L^\infty }\left( \Omega  \right)}}\left\| k \right\|_{{L^m}\left( \Omega  \right)}^2{\left( {\int_{B\left( {{x_0},\frac{r}{{{2^l}}}} \right)} {{{\left( {1 + {{\left| {{\nabla _H}u} \right|}^2}} \right)}^{\frac{{m\left( {2q - p} \right)}}{{2\left( {m - 2} \right)}}}}dx} } \right)^{\frac{{m - 2}}{m}}} \nonumber \\
 \le& \varepsilon \int_\Omega  {{\eta ^4}a\left( x \right){{\left( {1 + {{\left| {{\nabla _H}u} \right|}^2} + {{\left| {{\nabla _H}u\left( {x{e^{s{X_{{i_0}}}}}} \right)} \right|}^2}} \right)}^{\frac{{p - 2}}{2}}}{{\left| {{D_{{X_{{i_0}}},s,1}}{\nabla _H}u} \right|}^2}dx} \nonumber \\
 & + {c_\varepsilon }{\left\| {{a^{ - 1}}} \right\|_{{L^\infty }\left( \Omega  \right)}}\left\| k \right\|_{{L^m}\left( \Omega  \right)}^2\left( {\int_{B\left( {{x_0},\frac{r}{{{2^l}}}} \right)} {{{\left( {1 + {{\left| {{\nabla _H}u} \right|}^2}} \right)}^{\frac{{m\left( {2q - p} \right)}}{{2\left( {m - 2} \right)}}}}dx}  + 1} \right).
\end{align}

It deduces from \eqref{eq14} that
\begin{align*}
  {J_{11}} \le & \int_\Omega  {b\left( x \right){{\left( {1 + {{\left| {{\nabla _H}u\left( x \right)} \right|}^2} + {{\left| {{\nabla _H}u\left( {x{e^{ - s{X_{{i_0}}}}}} \right)} \right|}^2}} \right)}^{\frac{{q - 2}}{2}}}\left| {{D_{{X_{{i_0}}}, - s,1}}{\nabla _H}u\left( x \right)} \right|\left( {{\eta ^4}\left| {Tu} \right| + 4{\eta ^3}\left| {T\eta } \right|\left| u \right|} \right)dx} \\
   &  + \int_\Omega  {b\left( x \right){{\left( {1 + {{\left| {{\nabla _H}u\left( x \right)} \right|}^2} + {{\left| {{\nabla _H}u\left( {x{e^{s{X_{{i_0}}}}}} \right)} \right|}^2}} \right)}^{\frac{{q - 2}}{2}}}\left| {{D_{{X_{{i_0}}},s,1}}{\nabla _H}u\left( x \right)} \right|\left( {{\eta ^4}\left| {Tu} \right| + 4{\eta ^3}\left| {T\eta } \right|\left| u \right|} \right)dx} \\
   = :&{J_{111}} + {J_{112}},
\end{align*}
and then applies Young's inequality and H\"{o}lder's inequality to get
\begin{align*}
   {J_{111}} \le & \varepsilon \int_\Omega  {{\eta ^4}a\left( x \right){{\left( {1 + {{\left| {{\nabla _H}u} \right|}^2} + {{\left| {{\nabla _H}u\left( {x{e^{ - s{X_{{i_0}}}}}} \right)} \right|}^2}} \right)}^{\frac{{p - 2}}{2}}}{{\left| {{D_{{X_{{i_0}}}, - s,1}}{\nabla _H}u} \right|}^2}dx}  \\
   &  + {c_\varepsilon }\int_\Omega  {{\eta ^4}\frac{{b{{\left( x \right)}^2}}}{{a\left( x \right)}}{{\left( {1 + {{\left| {{\nabla _H}u} \right|}^2} + {{\left| {{\nabla _H}u\left( {x{e^{ - s{X_{{i_0}}}}}} \right)} \right|}^2}} \right)}^{\frac{{2q - 2 - p}}{2}}}{{\left| {Tu} \right|}^2}dx} \\
   &+ {c_\varepsilon }\int_\Omega  {{\eta ^2}{{\left| {T\eta } \right|}^2}\frac{{b{{\left( x \right)}^2}}}{{a\left( x \right)}}{{\left( {1 + {{\left| {{\nabla _H}u} \right|}^2} + {{\left| {{\nabla _H}u\left( {x{e^{ - s{X_{{i_0}}}}}} \right)} \right|}^2}} \right)}^{\frac{{2q - 2 - p}}{2}}}{{\left| u \right|}^2}dx} \\
    \le &\varepsilon \int_\Omega  {{\eta ^4}a\left( x \right){{\left( {1 + {{\left| {{\nabla _H}u} \right|}^2} + {{\left| {{\nabla _H}u\left( {x{e^{ - s{X_{{i_0}}}}}} \right)} \right|}^2}} \right)}^{\frac{{p - 2}}{2}}}{{\left| {{D_{{X_{{i_0}}}, - s,1}}{\nabla _H}u} \right|}^2}dx} \\
   & + {c_\varepsilon }{\left\| {{a^{ - 1}}} \right\|_{{L^\infty }\left( \Omega  \right)}}\left\| b \right\|_{{L^m}\left( \Omega  \right)}^2\cdot\\
   &\;\;\;\;{\left( {\int_{B\left( {{x_0},\frac{r}{{{2^{l + 1}}}}} \right)} {{{\left( {{{\left( {1 + {{\left| {{\nabla _H}u} \right|}^2} + {{\left| {{\nabla _H}u\left( {x{e^{ - s{X_{{i_0}}}}}} \right)} \right|}^2}} \right)}^{\frac{{2q - 2 - p}}{2}}}{{\left| {Tu} \right|}^2}} \right)}^{\frac{m}{{m - 2}}}}dx} } \right)^{\frac{{m - 2}}{m}}}\\
   & + {c_\varepsilon }{\left\| {{a^{ - 1}}} \right\|_{{L^\infty }\left( \Omega  \right)}}\left\| b \right\|_{{L^m}\left( \Omega  \right)}^2\cdot\\
   &\;\;\;\;{\left( {\int_{B\left( {{x_0},\frac{r}{{{2^{l + 1}}}}} \right)} {{{\left( {{{\left( {1 + {{\left| {{\nabla _H}u} \right|}^2} + {{\left| {{\nabla _H}u\left( {x{e^{ - s{X_{{i_0}}}}}} \right)} \right|}^2}} \right)}^{\frac{{2q - 2 - p}}{2}}}{{\left| u \right|}^2}} \right)}^{\frac{m}{{m - 2}}}}dx} } \right)^{\frac{{m - 2}}{m}}}\\
   \le &\varepsilon \int_\Omega  {{\eta ^4}a\left( x \right){{\left( {1 + {{\left| {{\nabla _H}u} \right|}^2} + {{\left| {{\nabla _H}u\left( {x{e^{ - s{X_{{i_0}}}}}} \right)} \right|}^2}} \right)}^{\frac{{p - 2}}{2}}}{{\left| {{D_{{X_{{i_0}}}, - s,1}}{\nabla _H}u} \right|}^2}dx} \\
   & + {c_\varepsilon }{\left\| {{a^{ - 1}}} \right\|_{{L^\infty }\left( \Omega  \right)}}\left\| b \right\|_{{L^m}\left( \Omega  \right)}^2\cdot\\
   &\;\;\;\;\left( {\int_{B\left( {{x_0},\frac{r}{{{2^{l + 1}}}}} \right)} {{{\left( {{{\left( {1 + {{\left| {{\nabla _H}u} \right|}^2} + {{\left| {{\nabla _H}u\left( {x{e^{ - s{X_{{i_0}}}}}} \right)} \right|}^2}} \right)}^{\frac{{2q - 2 - p}}{2}}}{{\left| {Tu} \right|}^2}} \right)}^{\frac{m}{{m - 2}}}}dx}  + 1} \right)\\
   & + {c_\varepsilon }{\left\| {{a^{ - 1}}} \right\|_{{L^\infty }\left( \Omega  \right)}}\left\| b \right\|_{{L^m}\left( \Omega  \right)}^2\cdot\\
   &\;\;\;\;\left( {\int_{B\left( {{x_0},\frac{r}{{{2^{l + 1}}}}} \right)} {{{\left( {{{\left( {1 + {{\left| {{\nabla _H}u} \right|}^2} + {{\left| {{\nabla _H}u\left( {x{e^{ - s{X_{{i_0}}}}}} \right)} \right|}^2}} \right)}^{\frac{{2q - 2 - p}}{2}}}{{\left| u \right|}^2}} \right)}^{\frac{m}{{m - 2}}}}dx}  + 1} \right)\\
   \le &\varepsilon \int_\Omega  {{\eta ^4}a\left( x \right){{\left( {1 + {{\left| {{\nabla _H}u} \right|}^2} + {{\left| {{\nabla _H}u\left( {x{e^{ - s{X_{{i_0}}}}}} \right)} \right|}^2}} \right)}^{\frac{{p - 2}}{2}}}{{\left| {{D_{{X_{{i_0}}}, - s,1}}{\nabla _H}u} \right|}^2}dx} \\
   & + {c_\varepsilon }{\left\| {{a^{ - 1}}} \right\|_{{L^\infty }\left( \Omega  \right)}}\left\| b \right\|_{{L^m}\left( \Omega  \right)}^2\left( {\int_{B\left( {{x_0},\frac{r}{{{2^l}}}} \right)} {\left( {{{\left( {1 + {{\left| {{\nabla _H}u} \right|}^2}} \right)}^{\frac{{pm\left( {2q - 2 - p} \right)}}{{2\left[ {m\left( {p - 2} \right) - 2p} \right]}}}} + {{\left| u \right|}^p}} \right)dx}  + 1} \right)\\
   & + {c_\varepsilon }{\left\| {{a^{ - 1}}} \right\|_{{L^\infty }\left( \Omega  \right)}}\left\| b \right\|_{{L^m}\left( \Omega  \right)}^2\int_{B\left( {{x_0},\frac{r}{{{2^{l + 1}}}}} \right)} {{{\left| {Tu} \right|}^p}dx} ;
\end{align*}
similarly,
\begin{align*}
   {J_{112}} \le& \varepsilon \int_\Omega  {{\eta ^4}a\left( x \right){{\left( {1 + {{\left| {{\nabla _H}u} \right|}^2} + {{\left| {{\nabla _H}u\left( {x{e^{s{X_{{i_0}}}}}} \right)} \right|}^2}} \right)}^{\frac{{p - 2}}{2}}}{{\left| {{D_{{X_{{i_0}}},s,1}}{\nabla _H}u} \right|}^2}dx}  \\
   &  + {c_\varepsilon }{\left\| {{a^{ - 1}}} \right\|_{{L^\infty }\left( \Omega  \right)}}\left\| b \right\|_{{L^m}\left( \Omega  \right)}^2\left( {\int_{B\left( {{x_0},\frac{r}{{{2^l}}}} \right)} {\left( {{{\left( {1 + {{\left| {{\nabla _H}u} \right|}^2}} \right)}^{\frac{{pm\left( {2q - 2 - p} \right)}}{{2\left[ {m\left( {p - 2} \right) - 2p} \right]}}}} + {{\left| u \right|}^p}} \right)dx}  + 1} \right)\\
   & + {c_\varepsilon }{\left\| {{a^{ - 1}}} \right\|_{{L^\infty }\left( \Omega  \right)}}\left\| b \right\|_{{L^m}\left( \Omega  \right)}^2\int_{B\left( {{x_0},\frac{r}{{{2^{l + 1}}}}} \right)} {{{\left| {Tu} \right|}^p}dx} ,\\
\end{align*}
so
\begin{align}\label{eq46}
   {J_{11}} \le& \varepsilon \int_\Omega  {{\eta ^4}a\left( x \right){{\left( {1 + {{\left| {{\nabla _H}u} \right|}^2} + {{\left| {{\nabla _H}u\left( {x{e^{ - s{X_{{i_0}}}}}} \right)} \right|}^2}} \right)}^{\frac{{p - 2}}{2}}}{{\left| {{D_{{X_{{i_0}}}, - s,1}}{\nabla _H}u} \right|}^2}dx} \nonumber \\
   &  + \varepsilon \int_\Omega  {{\eta ^4}a\left( x \right){{\left( {1 + {{\left| {{\nabla _H}u} \right|}^2} + {{\left| {{\nabla _H}u\left( {x{e^{s{X_{{i_0}}}}}} \right)} \right|}^2}} \right)}^{\frac{{p - 2}}{2}}}{{\left| {{D_{{X_{{i_0}}}, - s,1}}{\nabla _H}u} \right|}^2}dx} \nonumber \\
   & + {c_\varepsilon }{\left\| {{a^{ - 1}}} \right\|_{{L^\infty }\left( \Omega  \right)}}\left\| b \right\|_{{L^m}\left( \Omega  \right)}^2\left( {\int_{B\left( {{x_0},\frac{r}{{{2^l}}}} \right)} {\left( {{{\left( {1 + {{\left| {{\nabla _H}u} \right|}^2}} \right)}^{\frac{{pm\left( {2q - 2 - p} \right)}}{{2\left[ {m\left( {p - 2} \right) - 2p} \right]}}}} + {{\left| u \right|}^p}} \right)dx}  + 1} \right)\nonumber \\
  & + {c_\varepsilon }{\left\| {{a^{ - 1}}} \right\|_{{L^\infty }\left( \Omega  \right)}}\left\| b \right\|_{{L^m}\left( \Omega  \right)}^2\int_{B\left( {{x_0},\frac{r}{{{2^{l + 1}}}}} \right)} {{{\left| {Tu} \right|}^p}dx} .
\end{align}
Using \eqref{eq15},
\begin{align*}
   {J_{12}} \le & \int_\Omega  {k\left( x \right){{\left( {1 + {{\left| {{\nabla _H}u\left( {x{e^{ - s{X_{{i_0}}}}}} \right)} \right|}^2}} \right)}^{\frac{{q - 1}}{2}}}\left( {{\eta ^4}\left| {Tu} \right| + 4{\eta ^3}\left| {T\eta } \right|\left| u \right|} \right)dx}  \\
   &  + \int_\Omega  {k\left( x \right){{\left( {1 + {{\left| {{\nabla _H}u\left( {x{e^{s{X_{{i_0}}}}}} \right)} \right|}^2}} \right)}^{\frac{{q - 1}}{2}}}\left( {{\eta ^4}\left| {Tu} \right| + 4{\eta ^3}\left| {T\eta } \right|\left| u \right|} \right)dx} \\
   = :&{J_{121}} + {J_{122}},
\end{align*}
then we apply H\"{o}lder's inequality and Young's inequality to see
\begin{align*}
   {J_{121}} \le& c{\left\| k \right\|_{{L^m}\left( \Omega  \right)}}{\left( {\int_{B\left( {{x_0},\frac{r}{{{2^{l + 1}}}}} \right)} {{{\left( {{{\left( {1 + {{\left| {{\nabla _H}u\left( {x{e^{ - s{X_{{i_0}}}}}} \right)} \right|}^2}} \right)}^{\frac{{q - 1}}{2}}}\left| {Tu} \right|} \right)}^{\frac{m}{{m - 1}}}}dx} } \right)^{\frac{{m - 1}}{m}}} \\
   &  + c{\left\| k \right\|_{{L^m}\left( \Omega  \right)}}{\left( {\int_{B\left( {{x_0},\frac{r}{{{2^{l + 1}}}}} \right)} {{{\left( {{{\left( {1 + {{\left| {{\nabla _H}u\left( {x{e^{ - s{X_{{i_0}}}}}} \right)} \right|}^2}} \right)}^{\frac{{q - 1}}{2}}}\left| u \right|} \right)}^{\frac{m}{{m - 1}}}}dx} } \right)^{\frac{{m - 1}}{m}}}\\
    \le& c{\left\| k \right\|_{{L^m}\left( \Omega  \right)}}\left( {\int_{B\left( {{x_0},\frac{r}{{{2^{l + 1}}}}} \right)} {{{\left( {{{\left( {1 + {{\left| {{\nabla _H}u\left( {x{e^{ - s{X_{{i_0}}}}}} \right)} \right|}^2}} \right)}^{\frac{{q - 1}}{2}}}\left| {Tu} \right|} \right)}^{\frac{m}{{m - 1}}}}dx}  + 1} \right)\\
    & + c{\left\| k \right\|_{{L^m}\left( \Omega  \right)}}\left( {\int_{B\left( {{x_0},\frac{r}{{{2^{l + 1}}}}} \right)} {{{\left( {{{\left( {1 + {{\left| {{\nabla _H}u\left( {x{e^{ - s{X_{{i_0}}}}}} \right)} \right|}^2}} \right)}^{\frac{{q - 1}}{2}}}\left| u \right|} \right)}^{\frac{m}{{m - 1}}}}dx}  + 1} \right)\\
    \le &c{\left\| k \right\|_{{L^m}\left( \Omega  \right)}}\left( {\int_{B\left( {{x_0},\frac{r}{{{2^l}}}} \right)} {\left( {{{\left( {1 + {{\left| {{\nabla _H}u} \right|}^2}} \right)}^{\frac{{pm\left( {q - 1} \right)}}{{2\left[ {m\left( {p - 1} \right) - p} \right]}}}} + {{\left| u \right|}^p}} \right)dx}  + 1} \right)\\
    & + c{\left\| k \right\|_{{L^m}\left( \Omega  \right)}}\int_{B\left( {{x_0},\frac{r}{{{2^{l + 1}}}}} \right)} {{{\left| {Tu} \right|}^p}dx} ;
\end{align*}
similarly,
\begin{align*}
   {J_{122}} \le& c{\left\| k \right\|_{{L^m}\left( \Omega  \right)}}\left( {\int_{B\left( {{x_0},\frac{r}{{{2^l}}}} \right)} {\left( {{{\left( {1 + {{\left| {{\nabla _H}u} \right|}^2}} \right)}^{\frac{{pm\left( {q - 1} \right)}}{{2\left[ {m\left( {p - 1} \right) - p} \right]}}}} + {{\left| u \right|}^p}} \right)dx}  + 1} \right) \\
   &  + c{\left\| k \right\|_{{L^m}\left( \Omega  \right)}}\int_{B\left( {{x_0},\frac{r}{{{2^{l + 1}}}}} \right)} {{{\left| {Tu} \right|}^p}dx} ,
\end{align*}
so
\begin{align}\label{eq47}
   {J_{12}} \le& c{\left\| k \right\|_{{L^m}\left( \Omega  \right)}}\left( {\int_{B\left( {{x_0},\frac{r}{{{2^l}}}} \right)} {\left( {{{\left( {1 + {{\left| {{\nabla _H}u} \right|}^2}} \right)}^{\frac{{pm\left( {q - 1} \right)}}{{2\left[ {m\left( {p - 1} \right) - p} \right]}}}} + {{\left| u \right|}^p}} \right)dx}  + 1} \right)\nonumber \\
   &  + c{\left\| k \right\|_{{L^m}\left( \Omega  \right)}}\int_{B\left( {{x_0},\frac{r}{{{2^{l + 1}}}}} \right)} {{{\left| {Tu} \right|}^p}dx} .
\end{align}

By \eqref{eq14}, we have
\begin{align*}
{J_{21}} \le& \int_\Omega  {b\left( x \right){{\left( {1 + {{\left| {{\nabla _H}u} \right|}^2} + {{\left| {{\nabla _H}u\left( {x{e^{s{X_{{i_0}}}}}} \right)} \right|}^2}} \right)}^{\frac{{q - 2}}{2}}}\left| {{D_{{X_{{i_0}}},s,1}}{\nabla _H}u} \right| \cdot }  \\
&\;\;\;\;\left| {{D_{{X_{{i_0}}},s,1}}\eta } \right|{\eta ^2}\left| {{\nabla _H}\eta \left( {x{e^{s{X_{{i_0}}}}}} \right)} \right|\left| {u\left( {x{e^{s{X_{{i_0}}}}}} \right)} \right|dx\\
   &+ \int_\Omega  {b\left( x \right){{\left( {1 + {{\left| {{\nabla _H}u} \right|}^2} + {{\left| {{\nabla _H}u\left( {x{e^{s{X_{{i_0}}}}}} \right)} \right|}^2}} \right)}^{\frac{{q - 2}}{2}}}\left| {{D_{{X_{{i_0}}},s,1}}{\nabla _H}u} \right| \cdot } \\
&\;\;\;\;\left| {{D_{{X_{{i_0}}},s,1}}\eta } \right|s\left| {\frac{{\eta {{\left( {x{e^{s{X_{{i_0}}}}}} \right)}^2} - {\eta ^2}}}{s}} \right|\left| {{\nabla _H}\eta \left( {x{e^{s{X_{{i_0}}}}}} \right)} \right|\left| {u\left( {x{e^{s{X_{{i_0}}}}}} \right)} \right|dx\\
 =& \int_\Omega  {\left[ {{\eta ^2}a{{\left( x \right)}^{\frac{1}{2}}}{{\left( {1 + {{\left| {{\nabla _H}u} \right|}^2} + {{\left| {{\nabla _H}u\left( {x{e^{s{X_{{i_0}}}}}} \right)} \right|}^2}} \right)}^{\frac{{p - 2}}{4}}}\left| {{D_{{X_{{i_0}}},s,1}}{\nabla _H}u} \right|} \right] \cdot } \\
 &\;\;\;\;\left[ {\left| {{D_{{X_{{i_0}}},s,1}}\eta } \right|\left| {{\nabla _H}\eta \left( {x{e^{s{X_{{i_0}}}}}} \right)} \right|a{{\left( x \right)}^{ - \frac{1}{2}}}b\left( x \right){{\left( {1 + {{\left| {{\nabla _H}u} \right|}^2} + {{\left| {{\nabla _H}u\left( {x{e^{s{X_{{i_0}}}}}} \right)} \right|}^2}} \right)}^{\frac{{2q - 2 - p}}{4}}}\left| {u\left( {x{e^{s{X_{{i_0}}}}}} \right)} \right|} \right]dx\\
 &+ \int_\Omega  {\left| {{D_{{X_{{i_0}}},s,1}}\eta } \right|\left| {\frac{{\eta {{\left( {x{e^{s{X_{{i_0}}}}}} \right)}^2} - {\eta ^2}}}{s}} \right|\left| {{\nabla _H}\eta \left( {x{e^{s{X_{{i_0}}}}}} \right)} \right|b\left( x \right) \cdot } \\
 &\;\;\;\;{\left( {1 + {{\left| {{\nabla _H}u} \right|}^2} + {{\left| {{\nabla _H}u\left( {x{e^{s{X_{{i_0}}}}}} \right)} \right|}^2}} \right)^{\frac{{q - 2}}{2}}}s\left| {{D_{{X_{{i_0}}},s,1}}{\nabla _H}u} \right|\left| {u\left( {x{e^{s{X_{{i_0}}}}}} \right)} \right|dx,
\end{align*}
and so by
\[s\left| {{D_{{X_{{i_0}}},s,1}}{\nabla _H}u} \right| = \left| {{\Delta _{{X_{{i_0}}},s}}{\nabla _H}u} \right| \le {\left( {1 + {{\left| {{\nabla _H}u} \right|}^2} + {{\left| {{\nabla _H}u\left( {x{e^{s{X_{{i_0}}}}}} \right)} \right|}^2}} \right)^{\frac{1}{2}}}\]
and
\[\left| {\frac{{\eta {{\left( {x{e^{s{X_{{i_0}}}}}} \right)}^2} - {\eta ^2}}}{s}} \right| \le \left| {{D_{{X_{{i_0}}},s,1}}\eta } \right|,\]
we apply Young's inequality to see
\begin{align*}
   {J_{21}} \le& \varepsilon \int_\Omega  {{\eta ^4}a\left( x \right){{\left( {1 + {{\left| {{\nabla _H}u} \right|}^2} + {{\left| {{\nabla _H}u\left( {x{e^{s{X_{{i_0}}}}}} \right)} \right|}^2}} \right)}^{\frac{{p - 2}}{2}}}{{\left| {{D_{{X_{{i_0}}},s,1}}{\nabla _H}u} \right|}^2}dx}  \\
   &  + {c_\varepsilon }\int_\Omega  {{{\left| {{D_{{X_{{i_0}}},s,1}}\eta } \right|}^2}{{\left| {{\nabla _H}\eta \left( {x{e^{s{X_{{i_0}}}}}} \right)} \right|}^2}\frac{{b{{\left( x \right)}^2}}}{{a\left( x \right)}}{{\left( {1 + {{\left| {{\nabla _H}u} \right|}^2} + {{\left| {{\nabla _H}u\left( {x{e^{s{X_{{i_0}}}}}} \right)} \right|}^2}} \right)}^{\frac{{2q - 2 - p}}{2}}}{{\left| {u\left( {x{e^{s{X_{{i_0}}}}}} \right)} \right|}^2}dx}\\
     &+ c\int_\Omega  {{{\left| {{D_{{X_{{i_0}}},s,1}}\eta } \right|}^2}\left| {{\nabla _H}\eta \left( {x{e^{s{X_{{i_0}}}}}} \right)} \right|b\left( x \right){{\left( {1 + {{\left| {{\nabla _H}u} \right|}^2} + {{\left| {{\nabla _H}u\left( {x{e^{s{X_{{i_0}}}}}} \right)} \right|}^2}} \right)}^{\frac{{q - 1}}{2}}}\left| {u\left( {x{e^{s{X_{{i_0}}}}}} \right)} \right|dx} .
\end{align*}
Then, it derives from ${\left| {{D_{{X_{{i_0}}},s,1}}\eta } \right|^2}\left| {{\nabla _H}\eta \left( {x{e^{s{X_{{i_0}}}}}} \right)} \right| \le c$, H\"{o}lder's inequality and Young's inequality that
\begin{align}\label{eq48}
   {J_{21}} \le& \varepsilon \int_\Omega  {{\eta ^4}a\left( x \right){{\left( {1 + {{\left| {{\nabla _H}u} \right|}^2} + {{\left| {{\nabla _H}u\left( {x{e^{s{X_{{i_0}}}}}} \right)} \right|}^2}} \right)}^{\frac{{p - 2}}{2}}}{{\left| {{D_{{X_{{i_0}}},s,1}}{\nabla _H}u} \right|}^2}dx} \nonumber \\
 & + {c_\varepsilon }\int_\Omega  {\frac{{b{{\left( x \right)}^2}}}{{a\left( x \right)}}{{\left( {1 + {{\left| {{\nabla _H}u} \right|}^2} + {{\left| {{\nabla _H}u\left( {x{e^{s{X_{{i_0}}}}}} \right)} \right|}^2}} \right)}^{\frac{{2q - 2 - p}}{2}}}{{\left| {u\left( {x{e^{s{X_{{i_0}}}}}} \right)} \right|}^2}dx} \nonumber \\
& + c\int_\Omega  {b\left( x \right){{\left( {1 + {{\left| {{\nabla _H}u} \right|}^2} + {{\left| {{\nabla _H}u\left( {x{e^{s{X_{{i_0}}}}}} \right)} \right|}^2}} \right)}^{\frac{{q - 1}}{2}}}\left| {u\left( {x{e^{s{X_{{i_0}}}}}} \right)} \right|dx}  \nonumber \\
   \le& \varepsilon \int_\Omega  {{\eta ^4}a\left( x \right){{\left( {1 + {{\left| {{\nabla _H}u} \right|}^2} + {{\left| {{\nabla _H}u\left( {x{e^{s{X_{{i_0}}}}}} \right)} \right|}^2}} \right)}^{\frac{{p - 2}}{2}}}{{\left| {{D_{{X_{{i_0}}},s,1}}{\nabla _H}u} \right|}^2}dx} \nonumber \\
   &   + {c_\varepsilon }{\left\| {{a^{ - 1}}} \right\|_{{L^\infty }\left( \Omega  \right)}}\left\| b \right\|_{{L^m}\left( \Omega  \right)}^2\cdot\nonumber \\
   &\;\;\;\;{\left( {\int_{B\left( {{x_0},\frac{r}{{{2^{l + 1}}}}} \right)} {{{\left( {1 + {{\left| {{\nabla _H}u} \right|}^2} + {{\left| {{\nabla _H}u\left( {x{e^{s{X_{{i_0}}}}}} \right)} \right|}^2}} \right)}^{\frac{{m\left( {2q - 2 - p} \right)}}{{2\left( {m - 2} \right)}}}}{{\left| {u\left( {x{e^{s{X_{{i_0}}}}}} \right)} \right|}^{\frac{{2m}}{{m - 2}}}}dx} } \right)^{\frac{{m - 2}}{m}}}\nonumber \\
   &  + c{\left\| b \right\|_{{L^m}\left( \Omega  \right)}}{\left( {\int_{B\left( {{x_0},\frac{r}{{{2^{l + 1}}}}} \right)} {{{\left( {{{\left( {1 + {{\left| {{\nabla _H}u} \right|}^2} + {{\left| {{\nabla _H}u\left( {x{e^{s{X_{{i_0}}}}}} \right)} \right|}^2}} \right)}^{\frac{{q - 1}}{2}}}\left| {u\left( {x{e^{s{X_{{i_0}}}}}} \right)} \right|} \right)}^{\frac{m}{{m - 1}}}}dx} } \right)^{\frac{{m - 1}}{m}}}\nonumber \\
  \le & \varepsilon \int_\Omega  {{\eta ^4}a\left( x \right){{\left( {1 + {{\left| {{\nabla _H}u} \right|}^2} + {{\left| {{\nabla _H}u\left( {x{e^{s{X_{{i_0}}}}}} \right)} \right|}^2}} \right)}^{\frac{{p - 2}}{2}}}{{\left| {{D_{{X_{{i_0}}},s,1}}{\nabla _H}u} \right|}^2}dx}  \nonumber \\
  & + {c_\varepsilon }{\left\| {{a^{ - 1}}} \right\|_{{L^\infty }\left( \Omega  \right)}}\left\| b \right\|_{{L^m}\left( \Omega  \right)}^2\left( {\int_{B\left( {{x_0},\frac{r}{{{2^l}}}} \right)} {\left( {{{\left( {1 + {{\left| {{\nabla _H}u} \right|}^2}} \right)}^{\frac{{mp\left( {2q - 2 - p} \right)}}{{2\left[ {m\left( {p - 2} \right) - 2p} \right]}}}} + {{\left| u \right|}^p}} \right)dx}  + 1} \right) \nonumber \\
   & + c{\left\| b \right\|_{{L^m}\left( \Omega  \right)}}\left( {\int_{B\left( {{x_0},\frac{r}{{{2^l}}}} \right)} {\left( {{{\left( {1 + {{\left| {{\nabla _H}u} \right|}^2}} \right)}^{\frac{{pm\left( {q - 1} \right)}}{{2\left[ {m\left( {p - 1} \right) - p} \right]}}}} + {{\left| u \right|}^p}} \right)dx}  + 1} \right).
\end{align}
We use \eqref{eq15}, H\"{o}lder's inequality and Young's inequality to obtain
\begin{align}\label{eq49}
   {J_{22}}& \le \int_\Omega  {\left| {{D_{{X_{{i_0}}},s,1}}\eta } \right|\eta {{\left( {x{e^{s{X_{{i_0}}}}}} \right)}^2}\left| {{\nabla _H}\eta \left( {x{e^{s{X_{{i_0}}}}}} \right)} \right|k\left( x \right){{\left( {1 + {{\left| {{\nabla _H}u\left( {x{e^{s{X_{{i_0}}}}}} \right)} \right|}^2}} \right)}^{\frac{{q - 1}}{2}}}\left| {u\left( {x{e^{s{X_{{i_0}}}}}} \right)} \right|dx}  \nonumber \\
  &\le c \int_\Omega  {k\left( x \right){{\left( {1 + {{\left| {{\nabla _H}u\left( {x{e^{s{X_{{i_0}}}}}} \right)} \right|}^2}} \right)}^{\frac{{q - 1}}{2}}}\left| {u\left( {x{e^{s{X_{{i_0}}}}}} \right)} \right|dx} \nonumber \\
   & \le c{\left\| k \right\|_{{L^m}\left( \Omega  \right)}}{\left( {\int_{B\left( {{x_0},\frac{r}{{{2^{l + 1}}}}} \right)} {{{\left( {{{\left( {1 + {{\left| {{\nabla _H}u\left( {x{e^{s{X_{{i_0}}}}}} \right)} \right|}^2}} \right)}^{\frac{{q - 1}}{2}}}\left| {u\left( {x{e^{s{X_{{i_0}}}}}} \right)} \right|} \right)}^{\frac{m}{{m - 1}}}}dx} } \right)^{\frac{{m - 1}}{m}}}\nonumber \\
 & \le c{\left\| k \right\|_{{L^m}\left( \Omega  \right)}}\left( {\int_{B\left( {{x_0},\frac{r}{{{2^{l + 1}}}}} \right)} {{{\left( {{{\left( {1 + {{\left| {{\nabla _H}u\left( {x{e^{s{X_{{i_0}}}}}} \right)} \right|}^2}} \right)}^{\frac{{q - 1}}{2}}}\left| {u\left( {x{e^{s{X_{{i_0}}}}}} \right)} \right|} \right)}^{\frac{m}{{m - 1}}}}dx}  + 1} \right) \nonumber \\
&\le c{\left\| k \right\|_{{L^m}\left( \Omega  \right)}}\left( {\int_{B\left( {{x_0},\frac{r}{{{2^l}}}} \right)} {\left( {{{\left( {1 + {{\left| {{\nabla _H}u} \right|}^2}} \right)}^{\frac{{pm\left( {q - 1} \right)}}{{2\left[ {m\left( {p - 1} \right) - p} \right]}}}} + {{\left| u \right|}^p}} \right)dx}  + 1} \right).
\end{align}

The estimations of ${J_{31}}$ and ${J_{32}}$ are similar to those of ${J_{21}}$ and ${J_{22}}$, respectively, and we have
\begin{align}\label{eq410}
   {J_{31}} \le& \varepsilon \int_\Omega  {{\eta ^4}a\left( x \right){{\left( {1 + {{\left| {{\nabla _H}u} \right|}^2} + {{\left| {{\nabla _H}u\left( {x{e^{s{X_{{i_0}}}}}} \right)} \right|}^2}} \right)}^{\frac{{p - 2}}{2}}}{{\left| {{D_{{X_{{i_0}}},s,1}}{\nabla _H}u} \right|}^2}dx} \nonumber \\
   &  + {c_\varepsilon }{\left\| {{a^{ - 1}}} \right\|_{{L^\infty }\left( \Omega  \right)}}\left\| b \right\|_{{L^m}\left( \Omega  \right)}^2\left( {\int_{B\left( {{x_0},\frac{r}{{{2^l}}}} \right)} {\left( {{{\left( {1 + {{\left| {{\nabla _H}u} \right|}^2}} \right)}^{\frac{{mp\left( {2q - 2 - p} \right)}}{{2\left[ {m\left( {p - 2} \right) - 2p} \right]}}}} + {{\left| u \right|}^p}} \right)dx}  + 1} \right)\nonumber \\
   & + c{\left\| b \right\|_{{L^m}\left( \Omega  \right)}}\left( {\int_{B\left( {{x_0},\frac{r}{{{2^l}}}} \right)} {\left( {{{\left( {1 + {{\left| {{\nabla _H}u} \right|}^2}} \right)}^{\frac{{pm\left( {q - 1} \right)}}{{2\left[ {m\left( {p - 1} \right) - p} \right]}}}} + {{\left| u \right|}^p}} \right)dx}  + 1} \right)
\end{align}
and
\begin{equation}\label{eq411}
{J_{32}} \le c{\left\| k \right\|_{{L^m}\left( \Omega  \right)}}\left( {\int_{B\left( {{x_0},\frac{r}{{{2^l}}}} \right)} {\left( {{{\left( {1 + {{\left| {{\nabla _H}u} \right|}^2}} \right)}^{\frac{{pm\left( {q - 1} \right)}}{{2\left[ {m\left( {p - 1} \right) - p} \right]}}}} + {{\left| u \right|}^p}} \right)dx}  + 1} \right).
\end{equation}

Applying \eqref{eq14}, Young's inequality and H\"{o}lder's inequality, it gives
\begin{align}\label{eq412}
   {J_{41}} \le& \int_\Omega  {b\left( x \right){{\left( {1 + {{\left| {{\nabla _H}u\left( x \right)} \right|}^2} + {{\left| {{\nabla _H}u\left( {x{e^{s{X_{{i_0}}}}}} \right)} \right|}^2}} \right)}^{\frac{{q - 2}}{2}}}\left| {{D_{{X_{{i_0}}},s,1}}{\nabla _H}u\left( x \right)} \right| \cdot } \nonumber \\
   &\;\;\;\;{\eta ^2}\left| {{D_{{X_{{i_0}}},s,1}}\eta } \right|\left| {{\nabla _H}\eta \left( {x{e^{s{X_{{i_0}}}}}} \right)} \right|\left| {u\left( {x{e^{s{X_{{i_0}}}}}} \right)} \right|dx\nonumber \\
   =& \int_\Omega  {\left[ {{\eta ^2}a{{\left( x \right)}^{\frac{1}{2}}}{{\left( {1 + {{\left| {{\nabla _H}u\left( x \right)} \right|}^2} + {{\left| {{\nabla _H}u\left( {x{e^{s{X_{{i_0}}}}}} \right)} \right|}^2}} \right)}^{\frac{{p - 2}}{4}}}\left| {{D_{{X_{{i_0}}},s,1}}{\nabla _H}u\left( x \right)} \right|} \right] \cdot } \nonumber \\
   &\;\;\;\;\left[ {\left| {{D_{{X_{{i_0}}},s,1}}\eta } \right|\left| {{\nabla _H}\eta \left( {x{e^{s{X_{{i_0}}}}}} \right)} \right|a{{\left( x \right)}^{ - \frac{1}{2}}}b\left( x \right){{\left( {1 + {{\left| {{\nabla _H}u} \right|}^2} + {{\left| {{\nabla _H}u\left( {x{e^{s{X_{{i_0}}}}}} \right)} \right|}^2}} \right)}^{\frac{{2q - 2 - p}}{4}}}\left| {u\left( {x{e^{s{X_{{i_0}}}}}} \right)} \right|} \right]dx \nonumber \\
   \le& \varepsilon \int_\Omega  {{\eta ^4}a\left( x \right){{\left( {1 + {{\left| {{\nabla _H}u} \right|}^2} + {{\left| {{\nabla _H}u\left( {x{e^{s{X_{{i_0}}}}}} \right)} \right|}^2}} \right)}^{\frac{{p - 2}}{2}}}{{\left| {{D_{{X_{{i_0}}},s,1}}{\nabla _H}u} \right|}^2}dx}\nonumber \\
   & + {c_\varepsilon }\int_\Omega  {{{\left| {{D_{{X_{{i_0}}},s,1}}\eta } \right|}^2}{{\left| {{\nabla _H}\eta \left( {x{e^{s{X_{{i_0}}}}}} \right)} \right|}^2}\frac{{b{{\left( x \right)}^2}}}{{a\left( x \right)}}{{\left( {1 + {{\left| {{\nabla _H}u} \right|}^2} + {{\left| {{\nabla _H}u\left( {x{e^{s{X_{{i_0}}}}}} \right)} \right|}^2}} \right)}^{\frac{{2q - 2 - p}}{2}}}{{\left| {u\left( {x{e^{s{X_{{i_0}}}}}} \right)} \right|}^2}dx} \nonumber \\
   \le& \varepsilon \int_\Omega  {{\eta ^4}a\left( x \right){{\left( {1 + {{\left| {{\nabla _H}u} \right|}^2} + {{\left| {{\nabla _H}u\left( {x{e^{s{X_{{i_0}}}}}} \right)} \right|}^2}} \right)}^{\frac{{p - 2}}{2}}}{{\left| {{D_{{X_{{i_0}}},s,1}}{\nabla _H}u} \right|}^2}dx}\nonumber \\
   & + {c_\varepsilon }\int_\Omega  {\frac{{b{{\left( x \right)}^2}}}{{a\left( x \right)}}{{\left( {1 + {{\left| {{\nabla _H}u} \right|}^2} + {{\left| {{\nabla _H}u\left( {x{e^{s{X_{{i_0}}}}}} \right)} \right|}^2}} \right)}^{\frac{{2q - 2 - p}}{2}}}{{\left| {u\left( {x{e^{s{X_{{i_0}}}}}} \right)} \right|}^2}dx} \nonumber \\
   \le &\varepsilon \int_\Omega  {{\eta ^4}a\left( x \right){{\left( {1 + {{\left| {{\nabla _H}u} \right|}^2} + {{\left| {{\nabla _H}u\left( {x{e^{s{X_{{i_0}}}}}} \right)} \right|}^2}} \right)}^{\frac{{p - 2}}{2}}}{{\left| {{D_{{X_{{i_0}}},s,1}}{\nabla _H}u} \right|}^2}dx}  \nonumber \\
   &  + {c_\varepsilon }{\left\| {{a^{ - 1}}} \right\|_{{L^\infty }\left( \Omega  \right)}}\left\| b \right\|_{{L^m}\left( \Omega  \right)}^2\cdot\nonumber \\
   &\;\;\;\;{\left( {\int_{B\left( {{x_0},\frac{r}{{{2^{l + 1}}}}} \right)} {{{\left( {1 + {{\left| {{\nabla _H}u} \right|}^2} + {{\left| {{\nabla _H}u\left( {x{e^{s{X_{{i_0}}}}}} \right)} \right|}^2}} \right)}^{\frac{{m\left( {2q - 2 - p} \right)}}{{2\left( {m - 2} \right)}}}}{{\left| {u\left( {x{e^{s{X_{{i_0}}}}}} \right)} \right|}^{\frac{{2m}}{{m - 2}}}}dx} } \right)^{\frac{{m - 2}}{m}}}\nonumber \\
    \le& \varepsilon \int_\Omega  {{\eta ^4}a\left( x \right){{\left( {1 + {{\left| {{\nabla _H}u} \right|}^2} + {{\left| {{\nabla _H}u\left( {x{e^{s{X_{{i_0}}}}}} \right)} \right|}^2}} \right)}^{\frac{{p - 2}}{2}}}{{\left| {{D_{{X_{{i_0}}},s,1}}{\nabla _H}u} \right|}^2}dx} \nonumber \\
   &  + {c_\varepsilon }{\left\| {{a^{ - 1}}} \right\|_{{L^\infty }\left( \Omega  \right)}}\left\| b \right\|_{{L^m}\left( \Omega  \right)}^2\cdot\nonumber \\
   &\;\;\;\;\left( {\int_{B\left( {{x_0},\frac{r}{{{2^{l + 1}}}}} \right)} {{{\left( {1 + {{\left| {{\nabla _H}u} \right|}^2} + {{\left| {{\nabla _H}u\left( {x{e^{s{X_{{i_0}}}}}} \right)} \right|}^2}} \right)}^{\frac{{m\left( {2q - 2 - p} \right)}}{{2\left( {m - 2} \right)}}}}{{\left| {u\left( {x{e^{s{X_{{i_0}}}}}} \right)} \right|}^{\frac{{2m}}{{m - 2}}}}dx}  + 1} \right)\nonumber \\
    \le &\varepsilon \int_\Omega  {{\eta ^4}a\left( x \right){{\left( {1 + {{\left| {{\nabla _H}u} \right|}^2} + {{\left| {{\nabla _H}u\left( {x{e^{s{X_{{i_0}}}}}} \right)} \right|}^2}} \right)}^{\frac{{p - 2}}{2}}}{{\left| {{D_{{X_{{i_0}}},s,1}}{\nabla _H}u} \right|}^2}dx} \nonumber \\
    & + {c_\varepsilon }{\left\| {{a^{ - 1}}} \right\|_{{L^\infty }\left( \Omega  \right)}}\left\| b \right\|_{{L^m}\left( \Omega  \right)}^2\left( {\int_{B\left( {{x_0},\frac{r}{{{2^l}}}} \right)} {\left( {{{\left( {1 + {{\left| {{\nabla _H}u} \right|}^2}} \right)}^{\frac{{pm\left( {2q - 2 - p} \right)}}{{2\left[ {m\left( {p - 2} \right) - 2p} \right]}}}} + {{\left| u \right|}^p}} \right)dx}  + 1} \right).
\end{align}
Similarly to the estimation of ${J_{22}}$, we use \eqref{eq15}, H\"{o}lder's inequality and Young's inequality to gain
\begin{align}\label{eq413}
   {J_{42}}& \le \int_\Omega  {\left| {{D_{{X_{{i_0}}},s,1}}\eta } \right|{\eta ^2}\left| {{\nabla _H}\eta \left( {x{e^{s{X_{{i_0}}}}}} \right)} \right|k\left( x \right){{\left( {1 + {{\left| {{\nabla _H}u\left( {x{e^{s{X_{{i_0}}}}}} \right)} \right|}^2}} \right)}^{\frac{{q - 1}}{2}}}\left| {u\left( {x{e^{s{X_{{i_0}}}}}} \right)} \right|dx} \nonumber \\
   &  \le c{\left\| k \right\|_{{L^m}\left( \Omega  \right)}}\left( {\int_{B\left( {{x_0},\frac{r}{{{2^l}}}} \right)} {\left( {{{\left( {1 + {{\left| {{\nabla _H}u} \right|}^2}} \right)}^{\frac{{pm\left( {q - 1} \right)}}{{2\left[ {m\left( {p - 1} \right) - p} \right]}}}} + {{\left| u \right|}^p}} \right)dx}  + 1} \right).
\end{align}

Similarly to the estimation of ${J_{41}}$, it shows by \eqref{eq14}, Young's inequality and H\"{o}lder's inequality that
\begin{align}\label{eq414}
  {J_{51}} \le & \int_\Omega  {b\left( x \right){{\left( {1 + {{\left| {{\nabla _H}u} \right|}^2} + {{\left| {{\nabla _H}u\left( {x{e^{s{X_{{i_0}}}}}} \right)} \right|}^2}} \right)}^{\frac{{q - 2}}{2}}}\left| {{D_{{X_{{i_0}}},s,1}}{\nabla _H}u} \right|{\eta ^3}\left| {{D_{{X_{{i_0}}},s,1}}{\nabla _H}\eta } \right|\left| {u\left( {x{e^{s{X_{{i_0}}}}}} \right)} \right|dx} \nonumber \\
  \le &  \varepsilon \int_\Omega  {{\eta ^4}a\left( x \right){{\left( {1 + {{\left| {{\nabla _H}u} \right|}^2} + {{\left| {{\nabla _H}u\left( {x{e^{s{X_{{i_0}}}}}} \right)} \right|}^2}} \right)}^{\frac{{p - 2}}{2}}}{{\left| {{D_{{X_{{i_0}}},s,1}}{\nabla _H}u} \right|}^2}dx}\nonumber \\
  & + {c_\varepsilon }\int_\Omega  {{\eta ^2}{{\left| {{D_{{X_{{i_0}}},s,1}}{\nabla _H}\eta } \right|}^2}\frac{{b{{\left( x \right)}^2}}}{{a\left( x \right)}}{{\left( {1 + {{\left| {{\nabla _H}u} \right|}^2} + {{\left| {{\nabla _H}u\left( {x{e^{s{X_{{i_0}}}}}} \right)} \right|}^2}} \right)}^{\frac{{2q - 2 - p}}{2}}}{{\left| {u\left( {x{e^{s{X_{{i_0}}}}}} \right)} \right|}^2}dx} \nonumber \\
   \le &\varepsilon \int_\Omega  {{\eta ^4}a\left( x \right){{\left( {1 + {{\left| {{\nabla _H}u} \right|}^2} + {{\left| {{\nabla _H}u\left( {x{e^{s{X_{{i_0}}}}}} \right)} \right|}^2}} \right)}^{\frac{{p - 2}}{2}}}{{\left| {{D_{{X_{{i_0}}},s,1}}{\nabla _H}u} \right|}^2}dx} \nonumber \\
   & + {c_\varepsilon }{\left\| {{a^{ - 1}}} \right\|_{{L^\infty }\left( \Omega  \right)}}\left\| b \right\|_{{L^m}\left( \Omega  \right)}^2\left( {\int_{B\left( {{x_0},\frac{r}{{{2^l}}}} \right)} {\left( {{{\left( {1 + {{\left| {{\nabla _H}u} \right|}^2}} \right)}^{\frac{{pm\left( {2q - 2 - p} \right)}}{{2\left[ {m\left( {p - 2} \right) - 2p} \right]}}}} + {{\left| u \right|}^p}} \right)dx}  + 1} \right).
\end{align}
Similarly to the estimation of ${J_{22}}$, we use \eqref{eq15}, H\"{o}lder's inequality and Young's inequality to derive
\begin{align}\label{eq415}
   {J_{52}}& \le \int_\Omega  {{\eta ^3}\left| {{D_{{X_{{i_0}}},s,1}}{\nabla _H}\eta } \right|k\left( x \right){{\left( {1 + {{\left| {{\nabla _H}u\left( {x{e^{s{X_{{i_0}}}}}} \right)} \right|}^2}} \right)}^{\frac{{q - 1}}{2}}}\left| {u\left( {x{e^{s{X_{{i_0}}}}}} \right)} \right|dx} \nonumber \\
   &  \le c{\left\| k \right\|_{{L^m}\left( \Omega  \right)}}\left( {\int_{B\left( {{x_0},\frac{r}{{{2^l}}}} \right)} {\left( {{{\left( {1 + {{\left| {{\nabla _H}u} \right|}^2}} \right)}^{\frac{{pm\left( {q - 1} \right)}}{{2\left[ {m\left( {p - 1} \right) - p} \right]}}}} + {{\left| u \right|}^p}} \right)dx}  + 1} \right).
\end{align}

As in the estimation of  ${J_{41}}$, it gets by using \eqref{eq14}, Young's inequality and H\"{o}lder's inequality that
\begin{align}\label{eq416}
  {J_{61}} \le & \int_\Omega  {b\left( x \right){{\left( {1 + {{\left| {{\nabla _H}u} \right|}^2} + {{\left| {{\nabla _H}u\left( {x{e^{s{X_{{i_0}}}}}} \right)} \right|}^2}} \right)}^{\frac{{q - 2}}{2}}}\left| {{D_{{X_{{i_0}}},s,1}}{\nabla _H}u} \right|{\eta ^3}\left| {{\nabla _H}\eta } \right|\left| {{D_{{X_{{i_0}}},s,1}}u} \right|dx} \nonumber \\
 \le  &  \varepsilon \int_\Omega  {{\eta ^4}a\left( x \right){{\left( {1 + {{\left| {{\nabla _H}u} \right|}^2} + {{\left| {{\nabla _H}u\left( {x{e^{s{X_{{i_0}}}}}} \right)} \right|}^2}} \right)}^{\frac{{p - 2}}{2}}}{{\left| {{D_{{X_{{i_0}}},s,1}}{\nabla _H}u} \right|}^2}dx}\nonumber \\
 & + {c_\varepsilon }\int_\Omega  {{\eta ^2}{{\left| {{\nabla _H}\eta } \right|}^2}\frac{{b{{\left( x \right)}^2}}}{{a\left( x \right)}}{{\left( {1 + {{\left| {{\nabla _H}u} \right|}^2} + {{\left| {{\nabla _H}u\left( {x{e^{s{X_{{i_0}}}}}} \right)} \right|}^2}} \right)}^{\frac{{2q - 2 - p}}{2}}}{{\left| {{D_{{X_{{i_0}}},s,1}}u} \right|}^2}dx} \nonumber \\
  \le &\varepsilon \int_\Omega  {{\eta ^4}a\left( x \right){{\left( {1 + {{\left| {{\nabla _H}u} \right|}^2} + {{\left| {{\nabla _H}u\left( {x{e^{s{X_{{i_0}}}}}} \right)} \right|}^2}} \right)}^{\frac{{p - 2}}{2}}}{{\left| {{D_{{X_{{i_0}}},s,1}}{\nabla _H}u} \right|}^2}dx} \nonumber \\
& + {c_\varepsilon }{\left\| {{a^{ - 1}}} \right\|_{{L^\infty }\left( \Omega  \right)}}\left\| b \right\|_{{L^m}\left( \Omega  \right)}^2\left( {\int_{B\left( {{x_0},\frac{r}{{{2^l}}}} \right)} {{{\left( {1 + {{\left| {{\nabla _H}u} \right|}^2}} \right)}^{\frac{{pm\left( {2q - 2 - p} \right)}}{{2\left[ {m\left( {p - 2} \right) - 2p} \right]}}}}dx}  + 1} \right)\nonumber \\
&+ {c_\varepsilon }{\left\| {{a^{ - 1}}} \right\|_{{L^\infty }\left( \Omega  \right)}}\left\| b \right\|_{{L^m}\left( \Omega  \right)}^2\int_{B\left( {{x_0},\frac{r}{{{2^{l + 1}}}}} \right)} {{{\left| {{D_{{X_{{i_0}}},s,1}}u} \right|}^p}dx} .
\end{align}
Similarly to the estimation of ${J_{22}}$, it yields from \eqref{eq15}, H\"{o}lder's inequality and Young's inequality that
\begin{align}\label{eq417}
  {J_{62}} \le & \int_\Omega  {{\eta ^3}\left| {{\nabla _H}\eta } \right|k\left( x \right){{\left( {1 + {{\left| {{\nabla _H}u\left( {x{e^{s{X_{{i_0}}}}}} \right)} \right|}^2}} \right)}^{\frac{{q - 1}}{2}}}\left| {{D_{{X_{{i_0}}},s,1}}u} \right|dx} \nonumber \\
\le   &  c{\left\| k \right\|_{{L^m}\left( \Omega  \right)}}\left( {\int_{B\left( {{x_0},\frac{r}{{{2^l}}}} \right)} {{{\left( {1 + {{\left| {{\nabla _H}u} \right|}^2}} \right)}^{\frac{{pm\left( {q - 1} \right)}}{{2\left[ {m\left( {p - 1} \right) - p} \right]}}}}dx}  + 1} \right)\nonumber \\
&+ c{\left\| k \right\|_{{L^m}\left( \Omega  \right)}}\int_{B\left( {{x_0},\frac{r}{{{2^{l+1}}}}} \right)} {{{\left| {{D_{{X_{{i_0}}},s,1}}u} \right|}^p}dx} .
\end{align}

Similarly to the treatment of ${J_{21}}$, we have by \eqref{eq14}, Young's inequality and H\"{o}lder's inequality that
\begin{align}\label{eq418}
  {J_{71}} \le& \int_\Omega  {b\left( x \right){{\left( {1 + {{\left| {{\nabla _H}u} \right|}^2} + {{\left| {{\nabla _H}u\left( {x{e^{s{X_{{i_0}}}}}} \right)} \right|}^2}} \right)}^{\frac{{q - 2}}{2}}}\left| {{D_{{X_{{i_0}}},s,1}}{\nabla _H}u} \right|\left| {{D_{{X_{{i_0}}},s,1}}\eta } \right|{\eta ^3}\left| {{\nabla _H}u\left( {x{e^{s{X_{{i_0}}}}}} \right)} \right|dx} \nonumber \\
   &   + \int_\Omega  {b\left( x \right){{\left( {1 + {{\left| {{\nabla _H}u} \right|}^2} + {{\left| {{\nabla _H}u\left( {x{e^{s{X_{{i_0}}}}}} \right)} \right|}^2}} \right)}^{\frac{{q - 2}}{2}}}\left| {{D_{{X_{{i_0}}},s,1}}{\nabla _H}u} \right| \cdot }  \nonumber \\
   &\;\;\;\;\left| {{D_{{X_{{i_0}}},s,1}}\eta } \right|s\left| {\frac{{\eta {{\left( {x{e^{s{X_{{i_0}}}}}} \right)}^3} - {\eta ^3}}}{s}} \right|\left| {{\nabla _H}u\left( {x{e^{s{X_{{i_0}}}}}} \right)} \right|dx\nonumber \\
    \le &\varepsilon \int_\Omega  {{\eta ^4}a\left( x \right){{\left( {1 + {{\left| {{\nabla _H}u} \right|}^2} + {{\left| {{\nabla _H}u\left( {x{e^{s{X_{{i_0}}}}}} \right)} \right|}^2}} \right)}^{\frac{{p - 2}}{2}}}{{\left| {{D_{{X_{{i_0}}},s,1}}{\nabla _H}u} \right|}^2}dx} \nonumber \\
  & + {c_\varepsilon }\int_\Omega  {{\eta ^2}{{\left| {{D_{{X_{{i_0}}},s,1}}\eta } \right|}^2}\frac{{b{{\left( x \right)}^2}}}{{a\left( x \right)}}{{\left( {1 + {{\left| {{\nabla _H}u} \right|}^2} + {{\left| {{\nabla _H}u\left( {x{e^{s{X_{{i_0}}}}}} \right)} \right|}^2}} \right)}^{\frac{{2q - 2 - p}}{2}}}{{\left| {{\nabla _H}u\left( {x{e^{s{X_{{i_0}}}}}} \right)} \right|}^2}dx}  \nonumber \\
   & + c\int_\Omega  {{{\left| {{D_{{X_{{i_0}}},s,1}}\eta } \right|}^2}b\left( x \right){{\left( {1 + {{\left| {{\nabla _H}u} \right|}^2} + {{\left| {{\nabla _H}u\left( {x{e^{s{X_{{i_0}}}}}} \right)} \right|}^2}} \right)}^{\frac{{q - 1}}{2}}}\left| {{\nabla _H}u\left( {x{e^{s{X_{{i_0}}}}}} \right)} \right|dx} \nonumber \\
    \le &\varepsilon \int_\Omega  {{\eta ^4}a\left( x \right){{\left( {1 + {{\left| {{\nabla _H}u} \right|}^2} + {{\left| {{\nabla _H}u\left( {x{e^{s{X_{{i_0}}}}}} \right)} \right|}^2}} \right)}^{\frac{{p - 2}}{2}}}{{\left| {{D_{{X_{{i_0}}},s,1}}{\nabla _H}u} \right|}^2}dx} \nonumber \\
   & + {c_\varepsilon }{\left\| {{a^{ - 1}}} \right\|_{{L^\infty }\left( \Omega  \right)}}\left\| b \right\|_{{L^m}\left( \Omega  \right)}^2\left( {\int_{B\left( {{x_0},\frac{r}{{{2^l}}}} \right)} {{{\left( {1 + {{\left| {{\nabla _H}u} \right|}^2}} \right)}^{\frac{{m\left( {2q - p} \right)}}{{2\left( {m - 2} \right)}}}}dx}  + 1} \right)\nonumber \\
  &+ c{\left\| b \right\|_{{L^m}\left( \Omega  \right)}}\left( {\int_{B\left( {{x_0},\frac{r}{{{2^l}}}} \right)} {{{\left( {1 + {{\left| {{\nabla _H}u} \right|}^2}} \right)}^{\frac{{mq}}{{\left( {m - 1} \right)}}}}dx}  + 1} \right).
\end{align}
Applying \eqref{eq15}, H\"{o}lder's inequality and Young's inequality, it follows
\begin{align}\label{eq419}
  {J_{72}} & \le \int_\Omega  {\left| {{D_{{X_{{i_0}}},s,1}}\eta } \right|\eta {{\left( {x{e^{s{X_{{i_0}}}}}} \right)}^3}k\left( x \right){{\left( {1 + {{\left| {{\nabla _H}u\left( {x{e^{s{X_{{i_0}}}}}} \right)} \right|}^2}} \right)}^{\frac{{q - 1}}{2}}}\left| {{\nabla _H}u\left( {x{e^{s{X_{{i_0}}}}}} \right)} \right|dx}\nonumber \\
   &  \le c{\left\| k \right\|_{{L^m}\left( \Omega  \right)}}\left( {\int_{B\left( {{x_0},\frac{r}{{{2^l}}}} \right)} {{{\left( {1 + {{\left| {{\nabla _H}u} \right|}^2}} \right)}^{\frac{{mq}}{{\left( {m - 1} \right)}}}}dx}  + 1} \right).
\end{align}

The estimates of ${J_8},{J_9},{J_{10}}$ are similar to the estimate of ${J_7}$, and the results are the same.

Taking the test function in \eqref{eq212}
\[\varphi  = {D_{{X_{{i_0}}},s,1}}{D_{{X_{{i_0}}}, - s,1}}\left( {{\eta ^4}u} \right)\left( x \right),\]
we can get the estimates similarly to ${J_i}(i = 0,1, \cdots ,10)$, just needing to replace $x{e^{s{X_{{i_0}}}}}$ with $x{e^{ - s{X_{{i_0}}}}}$. Then we add those estimates, take $\varepsilon $ small enough, and then use \eqref{eq18} to get
\begin{align*}
   & \int_{B\left( {{x_0},\frac{r}{{{2^{l + 2}}}}} \right)} {{{\left( {1 + {{\left| {{\nabla _H}u} \right|}^2}
+ {{\left| {{\nabla _H}u\left( {x{e^{s{X_{{i_0}}}}}} \right)} \right|}^2}} \right)}^{\frac{{p - 2}}{2}}}{{\left| {{D_{{X_{{i_0}}},s,1}}{\nabla _H}u} \right|}^2}dx}\\
   &+ \int_{B\left( {{x_0},\frac{r}{{{2^{l + 2}}}}} \right)} {{{\left( {1 + {{\left| {{\nabla _H}u} \right|}^2}+ {{\left| {{\nabla _H}u\left( {x{e^{-s{X_{{i_0}}}}}} \right)} \right|}^2}} \right)}^{\frac{{p - 2}}{2}}}{{\left| {{D_{{X_{{i_0}}}, - s,1}}{\nabla _H}u} \right|}^2}dx}  \\
   \le&  {\left\| {{a^{ - 1}}} \right\|_{{L^\infty }\left( {B\left( {{x_0},r} \right)} \right)}}\int_\Omega  {a\left( x \right){\eta ^4}{{\left( {1 + {{\left| {{\nabla _H}u} \right|}^2} + {{\left| {{\nabla _H}u\left( {x{e^{s{X_{{i_0}}}}}} \right)} \right|}^2}} \right)}^{\frac{{p - 2}}{2}}}{{\left| {{D_{{X_{{i_0}}},s,1}}{\nabla _H}u} \right|}^2}dx} \\
   & + {\left\| {{a^{ - 1}}} \right\|_{{L^\infty }\left( {B\left( {{x_0},r} \right)} \right)}}\int_\Omega  {a\left( x \right){\eta ^4}{{\left( {1 + {{\left| {{\nabla _H}u} \right|}^2} + {{\left| {{\nabla _H}u\left( {x{e^{ - s{X_{{i_0}}}}}} \right)} \right|}^2}} \right)}^{\frac{{p - 2}}{2}}}{{\left| {{D_{{X_{{i_0}}}, - s,1}}{\nabla _H}u} \right|}^2}dx} \\
   \le& c{\left\| {{a^{ - 1}}} \right\|_{{L^\infty }\left( \Omega  \right)}}\left( {\left\| k \right\|_{{L^m}\left( \Omega  \right)}^2 + \left\| b \right\|_{{L^m}\left( \Omega  \right)}^2} \right)\left( {\int_{B\left( {{x_0},\frac{r}{{{2^l}}}} \right)} {{{\left( {1 + {{\left| {{\nabla _H}u} \right|}^2}} \right)}^{\frac{{m\left( {2q - p} \right)}}{{2\left( {m - 2} \right)}}}}dx}  + 1} \right)\\
   &+ c{\left\| {{a^{ - 1}}} \right\|_{{L^\infty }\left( \Omega  \right)}}\left\| b \right\|_{{L^m}\left( \Omega  \right)}^2\left( {\int_{B\left( {{x_0},\frac{r}{{{2^l}}}} \right)} {{{\left( {1 + {{\left| {{\nabla _H}u} \right|}^2}} \right)}^{\frac{{pm\left( {2q - 2 - p} \right)}}{{2\left[ {m\left( {p - 2} \right) - 2p} \right]}}}}dx}  + 1} \right)\\
   & + c\left( {{{\left\| b \right\|}_{{L^m}\left( \Omega  \right)}} + {{\left\| k \right\|}_{{L^m}\left( \Omega  \right)}}} \right)\left( {\int_{B\left( {{x_0},\frac{r}{{{2^l}}}} \right)} {{{\left( {1 + {{\left| {{\nabla _H}u} \right|}^2}} \right)}^{\frac{{pm\left( {q - 1} \right)}}{{2\left[ {m\left( {p - 1} \right) - p} \right]}}}}dx}  + 1} \right)\\
   & + c\left( {{{\left\| {{a^{ - 1}}} \right\|}_{{L^\infty }\left( \Omega  \right)}}\left\| b \right\|_{{L^m}\left( \Omega  \right)}^2 + {{\left\| k \right\|}_{{L^m}\left( \Omega  \right)}}} \right)\int_{B\left( {{x_0},\frac{r}{{{2^{l + 1}}}}} \right)} {{{\left| {Tu} \right|}^p}dx} \\
   &+ c\left( {{{\left\| {{a^{ - 1}}} \right\|}_{{L^\infty }\left( \Omega  \right)}}\left\| b \right\|_{{L^m}\left( \Omega  \right)}^2 + {{\left\| k \right\|}_{{L^m}\left( \Omega  \right)}} + {{\left\| b \right\|}_{{L^m}\left( \Omega  \right)}}} \right)\int_{B\left( {{x_0},\frac{r}{{{2^l}}}} \right)} {{{\left| u \right|}^p}dx} \\
   & + c\left( {{{\left\| k \right\|}_{{L^m}\left( \Omega  \right)}} + {{\left\| {{a^{ - 1}}} \right\|}_{{L^\infty }\left( \Omega  \right)}}\left\| b \right\|_{{L^m}\left( \Omega  \right)}^2} \right)\int_{B\left( {{x_0},\frac{r}{{{2^{l + 1}}}}} \right)} {{{\left| {{D_{{X_{{i_0}}},s,1}}u} \right|}^p}dx} \\
   & + c\left( {{{\left\| b \right\|}_{{L^m}\left( \Omega  \right)}} + {{\left\| k \right\|}_{{L^m}\left( \Omega  \right)}}} \right)\left( {\int_{B\left( {{x_0},\frac{r}{{{2^l}}}} \right)} {{{\left( {1 + {{\left| {{\nabla _H}u} \right|}^2}} \right)}^{\frac{{mq}}{{\left( {m - 1} \right)}}}}dx}  + 1} \right)\\
    \le& c{\left\| {{a^{ - 1}}} \right\|_{{L^\infty }\left( \Omega  \right)}}\left( {\left\| k \right\|_{{L^m}\left( \Omega  \right)}^2 + \left\| b \right\|_{{L^m}\left( \Omega  \right)}^2} \right)\left( {\int_{B\left( {{x_0},\frac{r}{{{2^l}}}} \right)} {{{\left( {1 + {{\left| {{\nabla _H}u} \right|}^2}} \right)}^{\frac{{m\left( {2q - p} \right)}}{{2\left( {m - 2} \right)}}}}dx}  + 1} \right)\\
    & + c{\left\| {{a^{ - 1}}} \right\|_{{L^\infty }\left( \Omega  \right)}}\left\| b \right\|_{{L^m}\left( \Omega  \right)}^2\left( {\int_{B\left( {{x_0},\frac{r}{{{2^l}}}} \right)} {{{\left( {1 + {{\left| {{\nabla _H}u} \right|}^2}} \right)}^{\frac{{pm\left( {2q - 2 - p} \right)}}{{2\left[ {m\left( {p - 2} \right) - 2p} \right]}}}}dx}  + 1} \right)\\
    & + c\left( {{{\left\| b \right\|}_{{L^m}\left( \Omega  \right)}} + {{\left\| k \right\|}_{{L^m}\left( \Omega  \right)}}} \right)\left( {\int_{B\left( {{x_0},\frac{r}{{{2^l}}}} \right)} {{{\left( {1 + {{\left| {{\nabla _H}u} \right|}^2}} \right)}^{\frac{{pm\left( {q - 1} \right)}}{{2\left[ {m\left( {p - 1} \right) - p} \right]}}}}dx}  + 1} \right)\\
    & + c{\kappa ^{l+2}}\left( {\int_{B\left( {{x_0},2r} \right)} {\left( {{{\left( {1 + {{\left| {{\nabla _H}u} \right|}^2}} \right)}^{\frac{{mp\left( {2q - 2-p} \right)}}{2[{m\left( {p - 2} \right) - 2p}]}}} + {{\left| u \right|}^p}} \right)dx}  + 1} \right)\\
    & + c\left( {{{\left\| k \right\|}_{{L^m}\left( \Omega  \right)}} + {{\left\| {{a^{ - 1}}} \right\|}_{{L^\infty }\left( \Omega  \right)}}\left\| b \right\|_{{L^m}\left( \Omega  \right)}^2} \right)\int_{B\left( {{x_0},\frac{r}{{{2^l}}}} \right)} {{{\left| {{\nabla _H}u} \right|}^p}dx} \\
    & + c\left( {{{\left\| b \right\|}_{{L^m}\left( \Omega  \right)}} + {{\left\| k \right\|}_{{L^m}\left( \Omega  \right)}}} \right)\left( {\int_{B\left( {{x_0},\frac{r}{{{2^l}}}} \right)} {{{\left( {1 + {{\left| {{\nabla _H}u} \right|}^2}} \right)}^{\frac{{mq}}{{\left( {m - 1} \right)}}}}dx}  + 1} \right)\\
    \le& c\kappa \left[ {\int_{B\left( {{x_0},\frac{r}{{{2^l}}}} \right)} {\left[ {{{\left( {1 + {{\left| {{\nabla _H}u} \right|}^2}} \right)}^{\frac{{m\left( {2q - p} \right)}}{{2\left( {m - 2} \right)}}}} + {{\left( {1 + {{\left| {{\nabla _H}u} \right|}^2}} \right)}^{\frac{{pm\left( {2q - 2 - p} \right)}}{{2\left[ {m\left( {p - 2} \right) - 2p} \right]}}}} + {{\left( {1 + {{\left| {{\nabla _H}u} \right|}^2}} \right)}^{\frac{{pm\left( {q - 1} \right)}}{{2\left[ {m\left( {p - 1} \right) - p} \right]}}}}} \right]dx}  + 1} \right]\\
    & + c{\kappa ^{l+2}}\left[ {\int_{B\left( {{x_0},2r} \right)} {\left( {{{\left( {1 + {{\left| {{\nabla _H}u} \right|}^2}} \right)}^{\frac{{mp\left( {2q - 2-p} \right)}}{2[{m\left( {p - 2} \right) - 2p}]}}} + {{\left| u \right|}^p}} \right)dx}  + 1} \right]\\
     \le &c{\kappa ^{l+2}}\left[ {\int_{B\left( {{x_0},2r} \right)} {\left[ {{{\left( {1 + {{\left| {{\nabla _H}u} \right|}^2}} \right)}^{\frac{{mp(2q - 2-p)}}{2[{m(p - 2) - 2p}]}}} + {{\left| u \right|}^p}} \right]dx}  + 1} \right].
\end{align*}

So far, similarly to the proof of Lemma 3.1 in \cite{MZZ07}, by above estimate, we prove \eqref{eq110} for ${i_0} \in \left\{ {1, \cdots ,2n} \right\}$.

For ${i_0} \in \left\{ {n + 1, \cdots ,2n} \right\},$ repeating the above discussion, we can also prove \eqref{eq110}.

Note that
\[\int_{B\left( {{x_0},\frac{r}{{{2^{l + 2}}}}} \right)} {{{\left| {\nabla _H^2u\left( x \right)} \right|}^2}dx}  \le \int_{B\left( {{x_0},\frac{r}{{{2^{l + 2}}}}} \right)} {{{\left( {1 + {{\left| {{\nabla _H}u\left( x \right)} \right|}^2}} \right)}^{\frac{{p - 2}}{2}}}{{\left| {\nabla _H^2u\left( x \right)} \right|}^2}dx} ,\]
and by combining \eqref{eq110}, we have $\nabla _H^2u \in L_{loc}^2\left( \Omega  \right),$ i.e. $u \in HW_{loc}^{2,2}\left( \Omega  \right)$.

\begin{remark}\label{Re41}

From the above proof process, we see that the statement in Theorem \ref{Th12} also holds for the weak solution $u \in H{W^{1,{\tilde q}}}\left( \Omega  \right)$
 to system \eqref{eq12}. In fact, when we consider system \eqref{eq12}, an additional term $-\int_\Omega {T u \cdot \varphi dx} $ will appear in \eqref{eq42}. Let us use the integration by parts and Lemma \ref{Le22} to deduce
 \begin{align*}
-\int_\Omega {T u \cdot \varphi dx} = &{\int _\Omega }{D_{{X_{{i_0}}},s,1}}\left( {Tu} \right) \cdot \left( {{D_{{X_{{i_0}}},s,1}}\left( {{\eta ^4}u} \right)} \right)dx\\
 =& {\int _\Omega }T\left( {{D_{{X_{{i_0}}},s,1}}u} \right) \cdot {\eta ^4}\left( {x{e^{s{X_{{i_0}}}}}} \right) \cdot {D_{{X_{{i_0}}},s,1}}udx \\
 & + {\int _\Omega }T\left( {{D_{{X_{{i_0}}},s,1}}u} \right) \cdot {\eta ^3}\left( {x{e^{s{X_{{i_0}}}}}} \right)u\left( x \right) \cdot {D_{{X_{{i_0}}},s,1}}\eta dx\\
 & + {\int _\Omega }T\left( {{D_{{X_{{i_0}}},s,1}}u} \right) \cdot {\eta ^2}\left( {x{e^{s{X_{{i_0}}}}}} \right)\eta \left( x \right)u\left( x \right) \cdot {D_{{X_{{i_0}}},s,1}}\eta dx\\
  &+ {\int _\Omega }T\left( {{D_{{X_{{i_0}}},s,1}}u} \right) \cdot \eta \left( {x{e^{s{X_{{i_0}}}}}} \right){\eta ^2}\left( x \right)u\left( x \right) \cdot {D_{{X_{{i_0}}},s,1}}\eta dx\\
  & + {\int _\Omega }T\left( {{D_{{X_{{i_0}}},s,1}}u} \right) \cdot {\eta ^3}\left( x \right)u\left( x \right) \cdot {D_{{X_{{i_0}}},s,1}}\eta dx\\
   =&  - \frac{1}{2}{\int _\Omega }{\left( {{D_{{X_{{i_0}}},s,1}}u} \right)^2} \cdot T\left( {{\eta ^4}\left( {x{e^{s{X_{{i_0}}}}}} \right)} \right)dx\\
& - {\int _\Omega }\left( {{D_{{X_{{i_0}}},s,1}}u} \right) \cdot T\left( {{\eta ^3}\left( {x{e^{s{X_{{i_0}}}}}} \right)u\left( x \right) \cdot {D_{{X_{{i_0}}},s,1}}\eta } \right)dx\\
&- {\int _\Omega }\left( {{D_{{X_{{i_0}}},s,1}}u} \right) \cdot T\left( {{\eta ^2}\left( {x{e^{s{X_{{i_0}}}}}} \right)\eta \left( x \right)u\left( x \right) \cdot {D_{{X_{{i_0}}},s,1}}\eta } \right)dx\\
& - {\int _\Omega }\left( {{D_{{X_{{i_0}}},s,1}}u} \right) \cdot T\left( {\eta \left( {x{e^{s{X_{{i_0}}}}}} \right){\eta ^2}\left( x \right)u\left( x \right) \cdot {D_{{X_{{i_0}}},s,1}}\eta } \right)dx\\
& - {\int _\Omega }\left( {{D_{{X_{{i_0}}},s,1}}u} \right) \cdot T\left( {{\eta ^3}\left( x \right)u\left( x \right) \cdot {D_{{X_{{i_0}}},s,1}}\eta } \right)dx\\
 \le& c\int_{B\left( {{x_0},\frac{r}{{{2^{l + 1}}}}} \right)} {\left( {{{\left| {{D_{{X_{{i_0}}},s,1}}u} \right|}^2} + {{\left| u \right|}^2}} \right)dx}  + c\int_{B\left( {{x_0},\frac{r}{{{2^{l + 1}}}}} \right)} {{{\left| {Tu} \right|}^2}dx} \\
   \le& c\int_{B\left( {{x_0},\frac{r}{{{2^{l + 1}}}}} \right)} {\left( {{{\left| {{D_{{X_{{i_0}}},s,1}}u} \right|}^p} + {{\left| u \right|}^p} + 1} \right)dx}  + c\int_{B\left( {{x_0},\frac{r}{{{2^{l + 1}}}}} \right)} {{{\left| {Tu} \right|}^p}dx}\\
   \leq& c{\kappa }^{l+1}\left[ {\int_{B\left( {{x_{0}},2r}\right) }{\left[ {{{%
\left( {1+{{\left\vert {{\nabla _{H}}u\left( x\right) }\right\vert }^{2}}}%
\right) }^{\frac{{mp(2q-2-p)}}{2[{m(p-2)-2p}]}}}+{{\left\vert u\right\vert }%
^{p}}}\right] dx}+1}\right] .
\end{align*}
\end{remark}
The remaining proof is consistent with the previous proof.

\section{Proof of Theorem \ref{Th13}}

\textbf{Proof of Theorem \ref{Th13}.} Let $\eta $ be a cut-off function between $B\left( {{x_0},\frac{r}{{{2^{l + 2}}}}} \right)$
 and $B\left( {{x_0},\frac{r}{{{2^{l + 1}}}}} \right)$ with $\left| {{\nabla _H}\eta } \right| \le c$ and $\left| {T\eta } \right| \le c$. Take the test function in \eqref{eq212}
\[\varphi  = {\Delta _{T, - s}}\left( {{\eta ^2}{\Delta _{T,s}}u\left( x \right)} \right),\]
then
\begin{equation}\label{eq51}
\int_\Omega  {\langle{\Delta _{T,s}}\left( {A\left( {x,{\nabla _H}u} \right)} \right),{\nabla _H}\left( {{\eta ^2}{\Delta _{T,s}}u\left( x \right)} \right)\rangle dx}  = 0.
\end{equation}
For the sake of simplicity, we denote
\[{\tilde {\rm A}_1}: = A\left( {x,{\nabla _H}u\left( {x{e^{sT}}} \right)} \right) - A\left( {x,{\nabla _H}u\left( x \right)} \right)\]
and
\[{\tilde {\rm A}_2}: = A\left( {x{e^{sT}},{\nabla _H}u\left( {x{e^{sT}}} \right)} \right) - A\left( {x,{\nabla _H}u\left( {x{e^{sT}}} \right)} \right),\]
and then
\[{\Delta _{T,s}}\left( {A\left( {x,{\nabla _H}u\left( x \right)} \right)} \right) = {\tilde {\rm A}_1} + {\tilde {\rm A}_2}.\]
Using it and
\[{\nabla _H}\left( {{\eta ^2}{\Delta _{T,s}}u\left( x \right)} \right) = {\eta ^2}{\Delta _{T,s}}{\nabla _H}u\left( x \right) + 2\eta {\nabla _H}\eta  \otimes {\Delta _{T,s}}u\left( x \right),\]
we see that \eqref{eq51} becomes
\begin{align}\label{eq52}
0 =& \int_\Omega  {{\eta ^2}\langle{{\tilde {\rm A}}_1}, {\Delta _{T,s}}{\nabla _H}u\left( x \right)\rangle dx}  + 2\int_\Omega  {\eta \langle{{\tilde {\rm A}}_1},  {{\nabla _H}\eta  \otimes {\Delta _{T,s}}u\left( x \right)} \rangle dx} \nonumber \\
   & + \int_\Omega  {{\eta ^2}\langle{{\tilde {\rm A}}_2}, {\Delta _{T,s}}{\nabla _H}u\left( x \right)\rangle dx}  + 2\int_\Omega  {\eta \langle{{\tilde {\rm A}}_2} ,  {{\nabla _H}\eta  \otimes {\Delta _{T,s}}u\left( x \right)} \rangle dx}   \nonumber \\
    = :&{\tilde I_1} + {\tilde I_2} + {\tilde I_3} + {\tilde I_4}.
\end{align}

Before estimating ${\tilde I_1},\;{\tilde I_2},\;{\tilde I_3}$ and ${\tilde I_4}$, we note the following estimates obtained from \eqref{eq13}, \eqref{eq14} and \eqref{eq15}:
\begin{equation}\label{eq53}
  \langle{\tilde {\rm A}_1} , {\Delta _{T,s}}{\nabla _H}u\left( x \right)\rangle \ge a\left( x \right){\left( {1 + {{\left| {{\nabla _H}u\left( x \right)} \right|}^2} + {{\left| {{\nabla _H}u\left( {x{e^{sT}}} \right)} \right|}^2}} \right)^{\frac{{p - 2}}{2}}}{\left| {{\Delta _{T,s}}{\nabla _H}u\left( x \right)} \right|^2},
\end{equation}
\begin{equation}\label{eq54}
  \left| {{{\tilde {\rm A}}_1}} \right| \le b\left( x \right){\left( {1 + {{\left| {{\nabla _H}u\left( x \right)} \right|}^2} + {{\left| {{\nabla _H}u\left( {x{e^{sT}}} \right)} \right|}^2}} \right)^{\frac{{q - 2}}{2}}}\left| {{\Delta _{T,s}}{\nabla _H}u\left( x \right)} \right|,
\end{equation}
\begin{equation}\label{eq55}
  \left| {{{\tilde {\rm A}}_2}} \right| \le sk\left( x \right){\left( {1 + {{\left| {{\nabla _H}u\left( {x{e^{sT}}} \right)} \right|}^2}} \right)^{\frac{{q - 1}}{2}}}.
\end{equation}
Now let us estimate ${\tilde I_1},\;{\tilde I_2},\;{\tilde I_3}$ and ${\tilde I_4}$. It derives from \eqref{eq53} that
\begin{align}\label{eq56}
   {\tilde I_1}& \ge \int_\Omega  {{\eta ^2}a\left( x \right){{\left( {1 + {{\left| {{\nabla _H}u} \right|}^2} + {{\left| {{\nabla _H}u\left( {x{e^{sT}}} \right)} \right|}^2}} \right)}^{\frac{{p - 2}}{2}}}{{\left| {{\Delta _{T,s}}{\nabla _H}u} \right|}^2}dx} \nonumber \\
   &  \ge c\int_\Omega  {{\eta ^2}{{\left( {1 + {{\left| {{\nabla _H}u} \right|}^2}+ {{\left| {{\nabla _H}u\left( {x{e^{sT}}} \right)} \right|}^2}} \right)}^{\frac{{p - 2}}{2}}}{{\left| {{\Delta _{T,s}}{\nabla _H}u} \right|}^2}dx} .
\end{align}
By \eqref{eq54} and Young's inequality, we have
\begin{align}\label{eq57}
   {\tilde I_2} \le& 2\int_\Omega  {\eta \left| {{\nabla _H}\eta } \right|b\left( x \right){{\left( {1 + {{\left| {{\nabla _H}u} \right|}^2} + {{\left| {{\nabla _H}u\left( {x{e^{sT}}} \right)} \right|}^2}} \right)}^{\frac{{q - 2}}{2}}}\left| {{\Delta _{T,s}}{\nabla _H}u} \right|\left| {{\Delta _{T,s}}u} \right|dx} \nonumber\\
  \le &  \varepsilon \int_\Omega  {{\eta ^2}{{\left( {1 + {{\left| {{\nabla _H}u} \right|}^2}+ {{\left| {{\nabla _H}u\left( {x{e^{sT}}} \right)} \right|}^2}} \right)}^{\frac{{p - 2}}{2}}}{{\left| {{\Delta _{T,s}}{\nabla _H}u} \right|}^2}dx} \nonumber \\
  &+ {c_\varepsilon }\int_\Omega  {{{\left| {{\nabla _H}\eta } \right|}^2}b{{\left( x \right)}^2}{{\left( {1 + {{\left| {{\nabla _H}u} \right|}^2} + {{\left| {{\nabla _H}u\left( {x{e^{sT}}} \right)} \right|}^2}} \right)}^{\frac{{2q - 2 - p}}{2}}}{{\left| {{\Delta _{T,s}}u} \right|}^2}dx} .
\end{align}
We use \eqref{eq55} and Young's inequality to get
\begin{align}\label{eq58}
   {\tilde I_3} \le& s\int_\Omega  {{\eta ^2}k\left( x \right){{\left( {1 + {{\left| {{\nabla _H}u\left( {x{e^{sT}}} \right)} \right|}^2}} \right)}^{\frac{{q - 1}}{2}}}\left| {{\Delta _{T,s}}{\nabla _H}u} \right|dx} \nonumber \\
  \le &  \varepsilon \int_\Omega  {{\eta ^2}{{\left( {1 + {{\left| {{\nabla _H}u} \right|}^2}+ {{\left| {{\nabla _H}u\left( {x{e^{sT}}} \right)} \right|}^2}} \right)}^{\frac{{p - 2}}{2}}}{{\left| {{\Delta _{T,s}}{\nabla _H}u} \right|}^2}dx}\nonumber \\
  &+ {c_\varepsilon }{s^2}\int_\Omega  {{\eta ^2}k{{\left( x \right)}^2}{{\left( {1 + {{\left| {{\nabla _H}u\left( {x{e^{sT}}} \right)} \right|}^2}} \right)}^{\frac{{2q - p}}{2}}}dx} .
\end{align}
It implies by \eqref{eq55} that
\begin{equation}\label{eq59}
  {\tilde I_4} \le 2s\int_\Omega  {\eta \left| {{\nabla _H}\eta } \right|k\left( x \right){{\left( {1 + {{\left| {{\nabla _H}u\left( {x{e^{sT}}} \right)} \right|}^2}} \right)}^{\frac{{q - 1}}{2}}}\left| {{\Delta _{T,s}}u} \right|dx} .
\end{equation}
Substituting \eqref{eq56}-\eqref{eq59} into \eqref{eq52} and taking $\varepsilon $ small enough, it gives
\begin{align}\label{eq510}
   & \int_\Omega  {{\eta ^2}{{\left( {1 + {{\left| {{\nabla _H}u} \right|}^2}+ {{\left| {{\nabla _H}u\left( {x{e^{sT}}} \right)} \right|}^2}} \right)}^{\frac{{p - 2}}{2}}}{{\left| {{\Delta _{T,s}}{\nabla _H}u} \right|}^2}dx}  \nonumber \\
   \le&  c\int_\Omega  {{{\left| {{\nabla _H}\eta } \right|}^2}b{{\left( x \right)}^2}{{\left( {1 + {{\left| {{\nabla _H}u} \right|}^2} + {{\left| {{\nabla _H}u\left( {x{e^{sT}}} \right)} \right|}^2}} \right)}^{\frac{{2q - 2 - p}}{2}}}{{\left| {{\Delta _{T,s}}u} \right|}^2}dx}\nonumber \\
   & + c{s^2}\int_\Omega  {{\eta ^2}k{{\left( x \right)}^2}{{\left( {1 + {{\left| {{\nabla _H}u\left( {x{e^{sT}}} \right)} \right|}^2}} \right)}^{\frac{{2q - p}}{2}}}dx} \nonumber \\
   & + cs\int_\Omega  {\eta \left| {{\nabla _H}\eta } \right|k\left( x \right){{\left( {1 + {{\left| {{\nabla _H}u\left( {x{e^{sT}}} \right)} \right|}^2}} \right)}^{\frac{{q - 1}}{2}}}\left| {{\Delta _{T,s}}u} \right|dx} ,
\end{align}
so we divide the two sides by ${s^2}$, and use H\"{o}lder's inequality, Young's inequality and \eqref{eq18} to gain
\begin{align}\label{eq511}
  & \int_\Omega  {{\eta ^2}{{\left( {1 + {{\left| {{\nabla _H}u} \right|}^2}+ {{\left| {{\nabla _H}u\left( {x{e^{sT}}} \right)} \right|}^2}} \right)}^{\frac{{p - 2}}{2}}}{{\left| {{D_{T,s,1}}{\nabla _H}u} \right|}^2}dx} \nonumber \\
   \le&  c\int_\Omega  {{{\left| {{\nabla _H}\eta } \right|}^2}b{{\left( x \right)}^2}{{\left( {1 + {{\left| {{\nabla _H}u} \right|}^2} + {{\left| {{\nabla _H}u\left( {x{e^{sT}}} \right)} \right|}^2}} \right)}^{\frac{{2q - 2 - p}}{2}}}{{\left| {{D_{T,s,1}}u} \right|}^2}dx} \nonumber \\
    &  + c\int_\Omega  {{\eta ^2}k{{\left( x \right)}^2}{{\left( {1 + {{\left| {{\nabla _H}u\left( {x{e^{sT}}} \right)} \right|}^2}} \right)}^{\frac{{2q - p}}{2}}}dx}  \nonumber \\
    & + c\int_\Omega  {\eta \left| {{\nabla _H}\eta } \right|k\left( x \right){{\left( {1 + {{\left| {{\nabla _H}u\left( {x{e^{sT}}} \right)} \right|}^2}} \right)}^{\frac{{q - 1}}{2}}}\left| {{D_{T,s,1}}u} \right|dx}  \nonumber \\
     \le& c\left\| b \right\|_{{L^m}\left( \Omega  \right)}^2{\left( {\int_{B\left( {{x_0},\frac{r}{{{2^{l{\rm{ + 1}}}}}}} \right)} {{{\left[ {{{\left( {1 + {{\left| {{\nabla _H}u} \right|}^2} + {{\left| {{\nabla _H}u\left( {x{e^{sT}}} \right)} \right|}^2}} \right)}^{\frac{{2q - 2 - p}}{2}}}{{\left| {{D_{T,s,1}}u} \right|}^2}} \right]}^{\frac{m}{{m - 2}}}}dx} } \right)^{\frac{{m - 2}}{m}}} \nonumber \\
       & + c\left\| k \right\|_{{L^m}\left( \Omega  \right)}^2{\left( {\int_{B\left( {{x_0},\frac{r}{{{2^l}}}} \right)} {{{\left( {1 + {{\left| {{\nabla _H}u} \right|}^2}} \right)}^{\frac{{m\left( {2q - p} \right)}}{{2\left( {m - 2} \right)}}}}dx} } \right)^{\frac{{m - 2}}{m}}}\nonumber \\
       & + c{\left\| k \right\|_{{L^m}\left( \Omega  \right)}}{\left( {\int_{B\left( {{x_0},\frac{r}{{{2^{l{\rm{ + 1}}}}}}} \right)} {{{\left( {{{\left( {1 + {{\left| {{\nabla _H}u\left( {x{e^{sT}}} \right)} \right|}^2}} \right)}^{\frac{{q - 1}}{2}}}\left| {{D_{T,s,1}}u} \right|} \right)}^{\frac{m}{{m - 1}}}}dx} } \right)^{\frac{{m - 1}}{m}}} \nonumber \\
         \le&  c\left\| b \right\|_{{L^m}\left( \Omega  \right)}^2\left( {\int_{B\left( {{x_0},\frac{r}{{{2^{l{\rm{ + 1}}}}}}} \right)} {{{\left[ {{{\left( {1 + {{\left| {{\nabla _H}u} \right|}^2} + {{\left| {{\nabla _H}u\left( {x{e^{sT}}} \right)} \right|}^2}} \right)}^{\frac{{2q - 2 - p}}{2}}}{{\left| {{D_{T,s,1}}u} \right|}^2}} \right]}^{\frac{m}{{m - 2}}}}dx}  + 1} \right)\nonumber \\
        &+ c\left\| k \right\|_{{L^m}\left( \Omega  \right)}^2\left( {\int_{B\left( {{x_0},\frac{r}{{{2^l}}}} \right)} {{{\left( {1 + {{\left| {{\nabla _H}u} \right|}^2}} \right)}^{\frac{{m\left( {2q - p} \right)}}{{2\left( {m - 2} \right)}}}}dx}  + 1} \right) \nonumber \\
         & + c{\left\| k \right\|_{{L^m}\left( \Omega  \right)}}\left( {\int_{B\left( {{x_0},\frac{r}{{{2^{l{\rm{ + 1}}}}}}} \right)} {{{\left( {{{\left( {1 + {{\left| {{\nabla _H}u\left( {x{e^{sT}}} \right)} \right|}^2}} \right)}^{\frac{{q - 1}}{2}}}\left| {{D_{T,s,1}}u} \right|} \right)}^{\frac{m}{{m - 1}}}}dx}  + 1} \right) \nonumber \\
          \le& c\left\| b \right\|_{{L^m}\left( \Omega  \right)}^2\left( {\int_{B\left( {{x_0},\frac{r}{{{2^{l{\rm{ + 1}}}}}}} \right)} {\left[ {{{\left( {1 + {{\left| {{\nabla _H}u} \right|}^2} + {{\left| {{\nabla _H}u\left( {x{e^{sT}}} \right)} \right|}^2}} \right)}^{\frac{{pm\left( {2q - 2 - p} \right)}}{{2\left[ {m\left( {p - 2} \right) - 2p} \right]}}}} + {{\left| {{D_{T,s,1}}u} \right|}^p}} \right]dx}  + 1} \right)  \nonumber \\
            & + c\left\| k \right\|_{{L^m}\left( \Omega  \right)}^2\left( {\int_{B\left( {{x_0},\frac{r}{{{2^l}}}} \right)} {{{\left( {1 + {{\left| {{\nabla _H}u} \right|}^2}} \right)}^{\frac{{m\left( {2q - p} \right)}}{{2\left( {m - 2} \right)}}}}dx}  + 1} \right)\nonumber \\
             & + c{\left\| k \right\|_{{L^m}\left( \Omega  \right)}}\left( {\int_{B\left( {{x_0},\frac{r}{{{2^{l{\rm{ + 1}}}}}}} \right)} {\left( {{{\left( {1 + {{\left| {{\nabla _H}u\left( {x{e^{sT}}} \right)} \right|}^2}} \right)}^{\frac{{pm\left( {q - 1} \right)}}{{2\left[ {m\left( {p - 1} \right) - p} \right]}}}} + {{\left| {{D_{T,s,1}}u} \right|}^p}} \right)dx}  + 1} \right)\nonumber \\
            \le& c\left\| b \right\|_{{L^m}\left( \Omega  \right)}^2\left( {\int_{B\left( {{x_0},\frac{r}{{{2^l}}}} \right)} {{{\left( {1 + {{\left| {{\nabla _H}u} \right|}^2}} \right)}^{\frac{{pm\left( {2q - 2 - p} \right)}}{{2\left[ {m\left( {p - 2} \right) - 2p} \right]}}}}dx}  + 1} \right)\nonumber \\
          & + c\kappa^{l+1} \left( {\left\| b \right\|_{{L^m}\left( \Omega  \right)}^2 + {{\left\| k \right\|}_{{L^m}\left( \Omega  \right)}}} \right)\left( {\int_{B\left( {{x_0},2r} \right)} {\left[ {{{\left( {1 + {{\left| {{\nabla _H}u} \right|}^2}} \right)}^{\frac{{mp\left( {2q - 2-p} \right)}}{2[{m\left( {p - 2} \right) - 2p}]}}} + {{\left| u \right|}^p}} \right]dx}  + 1} \right) \nonumber \\
          & + c\left\| k \right\|_{{L^m}\left( \Omega  \right)}^2\left( {\int_{B\left( {{x_0},\frac{r}{{{2^l}}}} \right)} {{{\left( {1 + {{\left| {{\nabla _H}u} \right|}^2}} \right)}^{\frac{{m\left( {2q - p} \right)}}{{2(m - 2)}}}}dx}  + 1} \right) \nonumber \\
          & + c{\left\| k \right\|_{{L^m}\left( \Omega  \right)}}\left( {\int_{B\left( {{x_0},\frac{r}{{{2^l}}}} \right)} {{{\left( {1 + {{\left| {{\nabla _H}u} \right|}^2}} \right)}^{\frac{{pm\left( {q - 1} \right)}}{{2\left[ {m\left( {p - 1} \right) - p} \right]}}}}dx}  + 1} \right) \nonumber \\
           \le &c{\kappa ^{l+2}}\left( {\int_{B\left( {{x_0},2r} \right)} {\left[ {{{\left( {1 + {{\left| {{\nabla _H}u} \right|}^2}} \right)}^{\frac{{mp\left( {2q - 2-p} \right)}}{2[{m\left( {p - 2} \right) - 2p}]}}} + {{\left| u \right|}^p}} \right]dx}  + 1} \right).
\end{align}
Therefore, similarly to the proof of Lemma 3.1 in \cite{MZZ07}, \eqref{eq111} holds from \eqref{eq511}. From \eqref{eq111} and
\[\int_{B\left( {{x_0},\frac{r}{{{2^{l + 2}}}}} \right)} {{{\left| {T{\nabla _H}u\left( x \right)} \right|}^2}dx}  \le \int_{B\left( {{x_0},\frac{r}{{{2^{l + 2}}}}} \right)} {{{\left( {1 + {{\left| {{\nabla _H}u\left( x \right)} \right|}^2}} \right)}^{\frac{{p - 2}}{2}}}{{\left| {T{\nabla _H}u\left( x \right)} \right|}^2}dx} ,\]
it follows $T{\nabla _H}u = {\nabla _H}Tu \in L_{loc}^2\left( \Omega  \right)$.

\begin{remark}\label{Re51}
Similarly to the discussion in Remark \ref{Re33}, we can obtain from the above proof process that Theorem \ref{Th13} holds for the weak solution $u \in H{W^{1,{\tilde q}}}\left( \Omega  \right)$ to system \eqref{eq12}.
\end{remark}



\end{document}